\documentclass{siamart220329}
\usepackage[utf8]{inputenc}
\usepackage{amsmath, bm}
\usepackage{amsfonts}
\usepackage{amssymb}
\usepackage{geometry}
\usepackage{graphicx}
\usepackage[export]{adjustbox}
\usepackage{hyperref}
\usepackage{framed}
\usepackage{subcaption}
\usepackage{soul}
\usepackage{multirow}
\usepackage{tabu}
\usepackage{mathtools}
\usepackage{subfiles}
\usepackage{listings}
\usepackage{todonotes}
\usepackage{caption}
\usepackage{csquotes}
\usepackage{bbm}
\usepackage{tikz}
\graphicspath{{images}}
\usepackage{cleveref} 
\usepackage{adjustbox}
\usepackage[most]{tcolorbox} 
\usepackage{bibunits}

\usepackage{lscape}

\usepackage{algpseudocode}

\newcommand{\mylabel}[1]{\label{#1} \label{LINE#1}}

\usepackage{thmtools}

\newsiamremark{example}{Example}
\newsiamremark{property}{Property}

\theoremstyle{remark}
\newtheorem{remark}{Remark}

\let\oldref\ref
\renewcommand{\ref}[1]{(\oldref{#1})}

\makeatletter
\newcommand*{\labeltext}[2]{%
  \begingroup
    \def\@currentlabel{#2}%
    \label{#1}%
  \endgroup
  #2%
}
\makeatother

\usepackage{enumitem,amssymb}
\newlist{todolist}{itemize}{2}
\setlist[todolist]{label=$\square$}
\usepackage{pifont}

\newcommand\expect[1]{
\mathbb{E} \left[ {#1} \right]
}

\newcommand{\norm}[1]{\left\|#1\right\|}

\usepackage{tikz}
\usepackage{float}
\usetikzlibrary{arrows,automata,positioning}

\usepackage{mathtools}

\newcommand\includegraphicsifexists[2][width=\linewidth]{\IfFileExists{#2}{\includegraphics[#1]{#2}}{}}

\NewDocumentCommand{\mycite}{o m}{\citeauthor{#2} \IfValueTF{#1}{\cite[#1]{#2}}{\cite{#2}}}

\title{Monte-Carlo/Moments micro-macro Parareal method for unimodal and bimodal scalar McKean-Vlasov SDEs
\thanks{This version was last modified on \today.
\funding{This project has received funding from the European High-Performance Computing Joint Undertaking (JU) under grant agreement No. 955701. The JU receives
support from the European Union’s Horizon 2020 research and innovation programme and Belgium, France, Germany, and Switzerland.}}
}

\headers{MC-moments Parareal for scalar McKean-Vlasov SDEs}{Ignace Bossuyt, Stefan Vandewalle, Giovanni Samaey}

\author{Ignace Bossuyt\thanks{Department  of  Computer  Science,  
KU  Leuven,  Celestijnenlaan 200A,  3001  Leuven,  Belgium (\email{ignace.bossuyt1@kuleuven.be}, \email{stefan.vandewalle@kuleuven.be}, \email{giovanni.samaey@kuleuven.be})}
\and
Stefan Vandewalle\footnotemark[2]
\and Giovanni Samaey\footnotemark[2]}

\begin{document}
\maketitle

\paragraph{Abstract}
We propose a micro-macro parallel-in-time Parareal method for scalar McKean-Vlasov stochastic differential equations (SDEs). 
In the algorithm, the fine Parareal propagator is a Monte Carlo simulation of an ensemble of particles, while an approximate ordinary differential equation (ODE) description of the mean and the variance of the particle distribution is used as a coarse Parareal propagator to achieve speedup.
We analyse the convergence behaviour of our method for a linear problem and provide numerical experiments indicating the convergence behavior of the algorithm on a set of examples.
We show, with numerical experiments, that convergence typically takes place in a low number of iterations, depending on the quality of the ODE predictor. 
For bimodal SDEs, we use multiple ODEs, each describing the mean and variance of the particle distribution in locally unimodal regions of the phase space.
For bimodal SDEs, we also develop a variant that converges faster by adaptively learning a model for the distribution of particles over different regions of the phase space (through a least-squares procedure).
The benefit of the proposed algorithm can be viewed through two lenses: (i) through the parallel-in-time lens, speedup is obtained through the use of a very cheap coarse integrator (an ODE moment model), and (ii) through the moment models lens, accuracy is iteratively gained through the use of parallel machinery as a corrector. 
In contrast to the isolated use of a moment model, the proposed method (iteratively) converges to the true distribution generated by the underlying McKean-Vlasov SDE.



\section{Introduction and motivation}
\label{subSECTION_simulation_methods_SDE}

In this paper, we introduce a parallel-in-time algorithm for the numerical simulation of the solution of McKean-Vlasov (mean-field) stochastic differential equations (SDEs). 
The introduction is structured as follows.
In \cref{subsection_mckean_vlasov} we briefly discuss McKean-Vlasov SDEs. 
In \cref{SECTION_simulation_methods_SDE} we discuss existing methods for their simulation and we explain why, with respect to current state-of-the art approaches, time parallelisation can offer advantages. 
Then, in \cref{subsection_PINT} we briefly give an overview of existing parallel-in-time methods.
In \cref{subsection_goal_of_paper}, we detail the goal and outline of the paper.

\subsection{Intro to McKean-Vlasov SDEs}
\label{subsection_mckean_vlasov}
This introduction is based on \cite{haji_ali_multilevel_2018} and \cite{reisinger_adaptive_2022}.
Let $X \in \mathbb R^d$ be a $d$-dimensional process on a time interval $t \in [0, T]$ and $\lambda(t)$ the marginal law (distribution) of $X$ at time $t$. 
Let $a \in \mathbb{R}^d$ be the drift coefficient, $b \in \mathbb{R}^{d \times n}$ a diffusion coefficient,  $W(t) \in \mathbb{R}^n$ an $n$-dimensional Brownian motion.
We considr McKean-Vlasov SDEs of the form 
\begin{equation}
\begin{aligned}
dX &= a(X,\lambda(t),t) dt + b(X, \lambda(t), t) dW, \\
\lambda(0) &\sim p_0.
\end{aligned}
\label{general_McKean_Vlasov}
\end{equation}

\sloppy
In practice, the simulation of \eqref{general_McKean_Vlasov} requires (i) a discretisation with a finite amount of $P$ particles, and (ii) a time-discretisation.
Let $\bar X = \{ X^{(p)} (t) \}_{p=1}^P$ be a particle ensemble  that is initially distributed with law $p_0$ (at time $t=0$) and $\lambda_P(x,t) = \frac{1}{P} \sum_{p=1}^P \delta_{X^{(p)}(t)}(dx)$ be its distribution function.
The evolution of each particle $X^{(p)}  \in \mathbb{R}^d$ obeys
\begin{equation}
\begin{aligned}
dX^{(p)} &= 
a \left( X^{(p)}, \lambda_P(t), t \right) dt 
+ b \left(
X^{(p)}, \lambda_P(t), t \right) dW^{(p)}.
\end{aligned}
\label{general_equation}
\end{equation}
The dependence of the coefficients $a$ and $b$ on the empirical measure $\lambda_P$ creates a coupling between all the particles. 
Indeed, the drift and diffusion coefficients of each particle are not only determined by the position of the particle itself, but also by the distribution of all the other particles (see, e.g., \cite{snitzman_1991} for an introduction to McKean-Vlasov SDEs).
In this paper we adopt the It\^{o} interpretation of \eqref{general_equation}. 

Particle systems of the form \eqref{general_equation} have been used to model, for instance, networks of neurons \cite{baladron_mean_field_2012, bossy_synchronization_2019}, 
the synchronisation of nonlinear oscillators \cite{kostur_nonequilibrium_2002}, 
and for the stochastic simulation of the (deterministic) Burgers equation \cite{bossy_stochastic_1997}. 
McKean-Vlasov SDEs also arise in data assimilation methods such as ensemble Kalman filtering and other interacting particle ensemble methods for data assimilation (see, for instance, \cite{del_moral_stability_2018}). 

\subsection{Simulation of McKean-Vlasov SDEs}
\label{SECTION_simulation_methods_SDE}
Various techniques for the simulation of \eqref{general_equation} exist. Deterministic methods include the Gauss-quadrature method in \cite{kloeden_gauss_quadrature_2017}.
In this work, we use a stochastic method:
\paragraph{Euler-Maruyama time-stepping method for SDEs}
The most basic numerical discretisation scheme for one ensemble of particles, obeying a McKean-Vlasov SDE, is the Euler-Maruyama (EM) scheme (see e.g.,  \cite{kloeden_platen_1999}). 
Let $\Delta t$ be a time step, and let the index $n$ refer to time $t_n = n \Delta t$
For all particles $p=1..P$ in the McKean-Vlasov SDE \eqref{general_equation}, the method can be written as follows:
\begin{equation}
X_{n+1}^{(p)} 
= 
a(X^{(p)}_n, \lambda_P (\bar X_n), t_n) 
\Delta t + 
b(X^{(p)}_n, \lambda_P (\bar X_n), t_n) 
\sqrt{\Delta t} \xi,
\label{Euler_Maruyama}
\end{equation}
where $\xi \in \mathbb R$ is a standard normally distributed variable, $\xi \sim \mathcal{N}(0,1).$ 
In \cite{reisinger_adaptive_2022}, an adaptive variant of Euler-Maruyama has been developed.

\paragraph{Monte Carlo sampling methods for Euler-Maruyama}
Let $\Phi$ be a user-chosen function.
The expectation of a quantity of interest $\mathbb{E}[\Phi(\bar X)]$ can be estimated by computing the sample average of Q different $P$-particle ensembles:
\begin{equation}
M_{Q,P}(t) 
= 
\mathbb E_Q \left[ \mathbb E_P [\bar X^{(q)}] \right]
=
\frac{1}{Q} \sum_{q=1}^Q
\frac{1}{P} \sum_{p=1}^P \Phi(X^{(q,p)}(t)),
\label{equation_MC_estimator}
\end{equation}

\begin{remark}[About expectations]
In \eqref{equation_MC_estimator}, the expectation $\mathbb E_P$ goes over the particles in one ensemble.
In practice, the expectation $\mathbb E$ over all ensembles is approximated by $\mathbb E \approx \mathbb E_Q$ where $Q$ is a finite number of ensembles.
In the sequel, we will always work with finite $P$-particle ensembles, but we will not complicate the notation by writing $\mathbb E_P$ (for the mean) and $\mathbb V_P$ (for the variance), instead we will just write $\mathbb E$ and $\mathbb V$.
For classical SDEs, $P$ can be safely put to 1 and then $Q$ refers to the number of independent particles in the Monte Carlo method.
\end{remark}


The Monte Carlo simulation of McKean-Vlasov SDEs can be computationally expensive. 
Moment ODEs have been presented in \cite[p. 139]{jazwinski_stochastic_1970} as a cheap alternative for SDEs without mean-field coupling. 
Moment ODEs have been used, for instance, to model the stochastic spiking of neural networks in \cite{rodriguez_statistical_1996}. 
The solution to these ODEs requires no sampling, and therefore is much cheaper than a stochastic particle simulation, but these ODEs contain a model error.

\subsection{Parallelisation and parallel-in-time methods}
\label{subsection_PINT}
The simulation of SDEs without mean-field interaction can easily be parallelised over all stochastic realisations.
Simulations with mean-field coupling, however, are less trivially parallelisable. In \cite[Figure 5]{haji_ali_multilevel_2018}, it is pointed out that there can still be a large amount of work in multilevel Monte Carlo (MLMC) and multi-index Monte Carlo schemes for McKean-Vlasov SDEs that can not be parallelised. 

To increase parallelism in the simulation of interacting particles, parallel-in-time methods can be envisaged.
In parallel-in-time methods, the time domain is divided in different slices, on which simulation can be performed in parallel.
For a history of time-parallel methods, see, e.g., \cite{Gander2015} and \cite{ong_applications_2020}. 
For deterministic models, various methods have been proposed, including the Parareal algorithm \cite{lions_resolution_2001}, MGRIT \cite{falgout_parallel_2014}, and PFASST \cite{emmett_toward_2012}. 
The Parareal algorithm uses an expensive accurate time propagator, 
which is applied in parallel over all time slices, to correct the result of a coarse but approximate simulation method,
which is applied sequentially over the time domain. 
The Parareal algorithm was analysed for linear ODEs and PDEs in \cite{gander_analysis_2007}.

The micro-macro Parareal algorithm is a generalisation of the original Parareal algorithm, that allows to use a coarse propagator that acts on a reduced state variable, instead of possibly high-dimensional and multiscale original coordinates \cite{blouza_parallel_2010}, \cite{Legoll2013}. 
The method is designed for multiscale (stiff) systems, where the coarse model is a cheaper non-stiff reduced model, but with model error. 

\subsection{Goal of this paper and related work}
\label{subsection_goal_of_paper}

We propose a new micro-macro Parareal method for scalar McKean-Vlasov SDEs that uses a Monte Carlo discretisation of the interacting particle system \eqref{general_equation}, using the Euler-Maruyama method as a fine Parareal solver and a low-dimensional moment model ODE as a coarse Parareal solver. 
The key advantage is that the coarse model is very cheap to simulate.
Then, we also build a variant of MC-moments Parareal that adaptively improves the coarse model as the iterations progress.
As a quantity of interest, we consider the expected value of a function $\Phi$ of the particle ensemble $\mathbb E_Q[\mathbb E_P[\Phi(\bar X)]] \approx \mathbb E[\Phi(\bar X)]$ (weak numerical approximation).

\paragraph{Existing Parareal algorithms for SDEs}
\label{motivation_for_weak_error}
We now discuss related work.
In \cite{Legoll2020}, the micro-macro Parareal algorithm is applied to SDEs with scale separation, where the coarse propagator is a finite volume discretisation and the fine propagator is a Monte Carlo simulation. 
In \cite{dabaghi_hybrid_2023}, a hybrid Parareal method is proposed that couples a Monte Carlo simulation on the fine level with a Galerkin scheme on the coarse level.
In \cite{thalhammer_MLMC_space_time_multigrid,
thalhammer_phd_thesis}, the multilevel Monte Carlo (MLMC) method is combined with Parareal.
The idea in those papers is to wrap a MLMC loop around a space-time multigrid solver for non-interacting particle systems. 
In \cite{engblom_parallel_2009}, a  Parareal method has been developed where the coarse solver is a reaction rate equation (ODE), and the fine scale solver is a stochastic simulation. 
In \cite{legoll_adaptive_Parareal}, an adaptive Parareal algorithm is developed for the strong approximation of a very long trajectory, in the context of molecular dynamics simulations.

The MC-moments Parareal algorithm also works for classical SDEs without mean-field interaction. 
However, on a massively parallel machine (allowing parallelisation over the then independent samples), the MC-moments Parareal algorithm would be useful only if the computation of one particle path takes longer than desired. When the number of samples is known a priori, it may be more useful to parallelise over the stochasticity rather than over time since the former could be done without a need for iterations. See also \cite{bal_parallelization_2003}.

\paragraph{Learning-based Parareal algorithms}
The idea of online learning in Parareal has been pursued in other contexts.
In \cite{gander_analysis_2008} a Krylov-enhanced Parareal for linear ODEs is analysed where the coarse solver is learned from the results of the fine propagator.
In \cite{Chen2014}, the reduced model, that is used as a coarse Parareal propagator, is updated after each parallel sweep of the fine propagator, using model reduction techniques.
In \cite{pentland_stochastic_2022} and \cite{pamela_neural-parareal_2025}, the Parareal correction term $\mathcal F_n - \mathcal C_n$ and the coarse propagator $\mathcal C_n$ are iteratively trained by learning a Gaussian process and a neural (coarse) operator, respectively.
We propose a cheap online learning procedure, based on a simple model for a part of the fine dynamics (i.e., only for the particle fractions in different regions where particles behave locally unimodally).
The learning step that we propose is, in addition, quite cheap.

The remainder of the paper is organised as follows. 
The fine integrator in the proposed MC-moments Parareal method uses the Euler-Maruyama discretisation.
In \cref{SETION_moment_equations}, we outline the moment ODEs that are used as a coarse Parareal integrator.  
They (approximately) model the mean and variance of the particle ensemble in each disjoint region of the phase space where the SDE locally behaves unimodally.
Then, in \cref{SECTION_new_algorithm}, we introduce the proposed MC-moments Parareal algorithm, using the moment models from \cref{SETION_moment_equations}. 
The algorithm is analysed on a very simple linear equation in \cref{SECTION_analysis}. 
We present numerical experiments 
in \cref{SECTION_numerical_exp} for unimodal and bimodal scalar McKean-Vlasov SDEs.
\Cref{SECTION_end} presents a conclusion and proposes some future research directions.

\section{Moment ODEs for scalar McKean-Vlasov SDEs}
\label{SETION_moment_equations}
In this section, we discuss the construction of moment models that approximate the dynamics of the mean and variance of McKean-Vlasov SDEs.

\paragraph{General expression}
An approximate moment ODE for the mean $M \approx \mathbb E[ \bar X ]$ and variance $\Sigma \approx \mathbb V [ \bar X ]$ can be written, generalised to McKean-Vlasov SDEs from \cite[equation (9.2)]{saarkkaa_solin_applied_SDEs_2019}:
\begin{equation}
\begin{aligned}
\frac{dM}{dt} 
	&= \expect{a(x,\lambda_P,t)}, \\
\frac{d\Sigma}{dt}
	&= 2 \expect{a(x,\lambda_P,t)} 
	+ 
	\expect{b(x,\lambda_P,t)^2}.
\end{aligned}
\label{general_expectation}
\end{equation}
In theory, the number of particles $P$ can, but needs not, go to infinity.
This equation is exact, but not practical, since the computation of the expectations requires the knowledge of the distribution of the particles. 
Different approximations are now possible.
Below we discuss two techniques: one technique is based on Gaussian models and the other is based on Taylor series.

In general, we desire that a moment model satisfies the following abstract requirements:
\begin{itemize}
\item The initial condition must be consistent with the given initial condition of the SDE. 

\item For SDEs that permit an invariant distribution, the moment model should correctly capture (the statistical moments of) the invariant distribution.
\end{itemize}

\paragraph{Section overview}
In this section we first  review existing moment models for classical SDEs (\cref{subsection_moment_models_classical_ODEs}) and we extend them to McKean-Vlasov SDEs. 
We then propose a technique that is also applicable to multimodal SDEs, a situation where the distribution is typically not concentrated around its mean point (see \cref{subsection_moment_model_classical_SDE_multimodel}). 
The proposed moment models are exact for linear McKean-Vlasov SDEs (see \cref{lemma_exactness_moment_model}), and approximate for nonlinear systems.

\subsection{Moment ODEs for classical and McKean-Vlasov SDEs}
\mylabel{subsection_moment_models_classical_ODEs}

We consider a special case of the McKean-Vlasov SDEs \eqref{general_equation}, namely where the mean-field effect $\lambda_P$ only enters the SDE via the expected value of a function $f$ of the particles:
\begin{equation}
dX^{(p)} = a(X^{(p)},\mathbb{E}[f(X^{(p)})],t)dt + b(X^{(p)}, \mathbb{E}[f(X^{(p)})],t)dW^{(p)}.
\label{SDE_class_1}
\end{equation}
Various techniques to approximate the evolution of $M$ and $\Sigma$ have been proposed. 
Here we specifically zoom in on two methods, and then briefly compare them and review some other existing methods.

\paragraph{Gaussian-assumed density approximations}
Let $p_{\mathcal N}(x-\mu, \Sigma)$ denote a Gaussian (normal) distribution function for $x \in \mathbb R$ with mean $\mu$ and variance $\Sigma$.
In a Gaussian-assumed density approximation, the distribution of particles under the expectation in equation \eqref{general_expectation} is assumed to be Gaussian:
\begin{equation}
\begin{aligned}
\frac{dM}{dt} 
&\approx
\int a(x, \lambda_P, t) p_{\mathcal N}(x-M,\Sigma) dx, \\
\frac{d \Sigma}{dt} 
&\approx 
2 \int a(x, \lambda_P, t) p_{\mathcal N}(x-M,\Sigma) dx 
+ 
\int b(x, \lambda_P, t) p_{\mathcal N}(x-M,\Sigma) dx.
\end{aligned}
\label{Gaussian_model_McK_V_SDE}
\end{equation}
Numerically, an integration rule (sigma-point method) can be applied to approximate the Gaussian integrals in \eqref{Gaussian_model_McK_V_SDE}, see 
\cite[equation (3)]{julier_new_2000}, \cite[equation (23)]{arasaratnam_cubature_2009}.

In \cite{archambeau07a}, a Gaussian process approximation of stochastic differential equations is proposed.
In the context of machine learning, a Gaussian model is used while doing inference using SDEs in \cite{Solin_NEURIPS}.

\paragraph{A technique based on Taylor series expansion of SDE coefficients} 
We here generalize a technique from \cite{rodriguez_statistical_1996} which was proposed to approximate classical SDEs (in \cite{rodriguez_statistical_1996} multivariate SDEs are considered).
Let $a_X$ and $b_X$ denote the derivatives of $a$ and $b$ with respect to their first argument, and $b_{XX}$ denotes the second derivative of $b$ with respect to its first argument.
This system of ODEs can be used as a moment approximation:
\begin{equation}
\begin{aligned}
\frac{dM}{dt} 
&\approx 
a(M,f(M),t) 
	+ \frac{1}{2} a_{XX}(M,f(M),t)
	+ \frac{1}{2} b_{XX}(M,f(M),t) \Sigma,
\quad & M(0) = \mathbb{E}[\bar X(0)], \\
\frac{d\Sigma}{dt} 
&\approx
\left[ 2a_X(M,f(M),t) + b_X(M,f(M),t)^2 \right]\Sigma + b(M,f(M),t)^2, 
\quad & \Sigma(0) = \mathbb{V}[\bar X(0)]. 
\end{aligned}
\label{moment_model_class_1}
\end{equation}
This ODE is given in \cite{rodriguez_statistical_1996} and its derivation is based on a combination of It\^{o}'s lemma with a Taylor expansion of the drift and diffusion coefficients around the mean.
The model was improved by adding the last term in the evolution of the mean. 
For more information, see \cref{appendix_derivation_moment_model}.

\paragraph{Comparison, and other techniques}
It is not always clear a priori which moment model is best suitable for which SDE. 
In any case, the Taylor-based model is exact for linear McKean-Vlasov SDEs, regardless of the initial condition\footnote{The exactness of the Taylor series-based moment model is not strictly limited to linear McKean-Vlasov SDEs. It is also exact, for instance, for the polynomial drift model from \cite[equation (14)]{son_doan_mean_square_2015}.}; for the proof see \cref{proof_exactness_linear_McK_V_SDE}.

\begin{lemma}[Exactness of the Taylor-based moment model for linear SDEs]
\label{lemma_exactness_moment_model}
For linear McKean-Vlasov SDEs of the form
\begin{equation}
\begin{aligned}
dX^{(p)} 
&= 
\left( A(t)X^{(p)} + A_E (t) \mathbb{E}[\bar X] + A_0(t)\right) dt 
+
\left(B(t)X^{(p)} + B_E (t)\mathbb{E}[\bar X] + B_0(t) \right) dW^{(p)},
\end{aligned}
\end{equation}
and with $X(0) \sim p_0$, 
the moment equations \eqref{moment_model_class_1} are an exact description of its mean and variance.
\end{lemma}

An overview of various other methods to approximate the evolution of $M$ and $\Sigma$, and possibly also the full particle distribution function, is given in the book \cite[Chapter 9]{saarkkaa_solin_applied_SDEs_2019}. 
See also \cite{sukys_momentclosurejl_2021} for practical details and implementations of various moment models.
In \cite{zagli_dimension_2023}, a moment model has been developed for McKean-Vlasov SDEs based on cumulants of the particle distribution.

\subsection{Moment ODEs for classical SDEs with multimodal distributions}
\mylabel{subsection_moment_model_classical_SDE_multimodel}

\paragraph{Gaussian mixture models with multimodal particle distributions}
For SDEs with multimodal distributions, one can use a Gaussian mixture (a weighted combination of Gaussians) instead of one single Gaussian.
This has been done in the context of Kalman filtering in \cite{alspach_nonlinear_1972}, \cite{terejanu_adaptive_2011}.
A random initialisation of the Gaussian samples/particles can be used in the hope to scan the entire stationary phase space. 
This has been done in the context of the estimation of stationary distributions in \cite{lambert_chewi_bach_bonnabel_rigollet} and for the approximation of (classical) SDEs in \cite{li_numerical_2021}.

A fundamental difficulty for these moment models is the  fulfillment of the desire to capture well the evolution of the SDE over the whole time domain, including (i) the initial condition, and (ii) the stationary distribution (if this exists). 
This requires a robust methods for the weights to change over the coarse of time.

As far as we know, however, there exists no mature Gaussian mixture method for the approximation of general McKean-Vlasov SDEs where the weights of the mixture, corresponding to different modes, are iteratively updated (in the spirit of \cite{alspach_nonlinear_1972}). 
The creation of such a method was deemed out of the scope of this work.

\paragraph{Moment model based on Taylor series}Unimodal (Taylor-based) moment models may be nonphysical if the distribution function of the particles is not concentrated and symmetric around its mean \cite{rodriguez_statistical_1996}, see also \cite{schnoerr_validity_2014}.
Bimodal SDEs are especially challenging.
Suppose that the regions of attraction $\mathcal D_i$ of an SDE are known, for all $1 \leq i \leq I$ with $I$ be the total number of regions of attraction (this assumption is not essential, see also \cref{remark_regions_of_attraction}).
Let $\mathbb{E}_{\mathcal D_i}$ and $\mathbb{V}_{\mathcal D_i}$ denote operators that compute the mean and variance of particles residing in $\mathcal D_i$.
Then a new technique that we propose consists in creating a distinct moment model ($M_i$, $\Sigma_i$) for each non-overlapping region of attraction of the phase space. For each $1 \leq i \leq I$, we have
\begin{equation}
\begin{aligned}
\frac{dM_i}{dt} 
&= 
F(M_i, \Sigma_i)
\qquad & 
M_i(0) = \mathbb{E}_{\mathcal D_i}[X(0)], \\
\frac{d\Sigma_i}{dt}
&= 
G(M_i, \Sigma_i)
\qquad & 
\Sigma_i(0) = \mathbb{V}_{\mathcal D_i}[X(0)],
\end{aligned}
\label{moment_model_Taylor_multimodal}
\end{equation}
where $F$ and $G$ refer to either the Gaussian unimodal model \eqref{Gaussian_model_McK_V_SDE} or the Taylor-based model \eqref{moment_model_class_1}.
Let $\mathcal P_{\mathcal D_i}$ denote the fraction of the particles residing in region $\mathcal D_i$, such that $\sum_{i=1}^I \mathcal P_{\mathcal D_i} = 1$.
The mean and variance of the whole particle ensemble is given by \cite[equations 1.20 and 1.21]{fruhwirth-schnatter_finite_2006}:
\begin{equation}
M 
= 
\sum_{i=1}^I \mathcal P_{\mathcal D_i} M_i,
\qquad
\Sigma^2  
=
\sum_{i=1}^I (M_i^2 + \Sigma_i^2) \mathcal P_{\mathcal D_i} - M^2. 
\label{mixture_of_distributions_mean_variance}
\end{equation} 

Here we do not propose a model for the time evolution of $\mathcal P_{\mathcal D_i}$, instead we refer to \cref{subsection_learning_based_Parareal}.
If the number of modes equals one, then the multimodal moment model \eqref{moment_model_Taylor_multimodal} coincides with the unimodal moment model \eqref{Gaussian_model_McK_V_SDE} or \eqref{moment_model_class_1}.

\begin{property}
\label{property_of_multimodal_linear_moment_model}
For linear McKean-Vlasov SDEs, one can either use the unimodal ODE approximation \eqref{moment_model_class_1}, or apply \eqref{mixture_of_distributions_mean_variance} on the local means and variances of the multimodal moment model \eqref{moment_model_Taylor_multimodal} with constant weights (particle fractions). 
As a result of the superposition principle and linearity of the SDE (and its moment ODEs), these descriptions are equivalent. This is proved in \cref{nog_een_bewijsken}.
\end{property}

\begin{remark}[Computation of the regions of attraction]
\label{remark_regions_of_attraction}
For SDEs with a stationary distribution, the boundaries of these regions (separatrices) can be computed, for instance, based on an analysis of the invariant distribution, see for instance \cite{tang_potential_2017} or \cite{najafi_fast_2016} or using a Gaussian mixture \cite{lambert_chewi_bach_bonnabel_rigollet}.
\end{remark}

\section{The MC-moments Parareal algorithm}
\label{SECTION_new_algorithm}
In this section, we formulate our main contribution, the MC-moments Parareal method. 
The aim of the Parareal method is to parallelise the simulation of initial value problems (IVPs) of the form $du/dt = f(u,t)$ with $u(0) = u_0$. 
Let $u_n$ be the solution to the IVP at time $t=n \Delta t$, and let $\mathcal F_n: \mathbb R^d \rightarrow \mathbb R^d$ be a time-stepping method such that $u_{n+1} = \mathcal F_n (u_n)$.
First, we discuss the Parareal algorithm.

\subsection{Background: the Parareal algorithm}
The two ingredients for Parareal are (i) a fine propagator $\mathcal F_n$, which is accurate but computationally expensive, and (ii) a coarse propagator $\mathcal C_n$ which is cheaper but less accurate. 
Iterations of the Parareal algorithm \cite{lions_resolution_2001} can be written as follows: the initialisation equals $u^k_0 = u_0$ for all $k \geq 0$, and then for all $k \geq 0$ and $n \geq 0$:
\begin{equation} 
\begin{aligned}
u_{n+1}^0 &= \mathcal C_n(u^0_n), \\
u_{n+1}^{k+1} 
&= \mathcal{C}_n (u_n^{k+1}) 
+ \mathcal{F}_n (u_n^k) 
- \mathcal{C}_n (u_n^k).
\end{aligned}
\label{original_parareal}
\end{equation}
Here, $u^k_n$ is the approximation at time point $n$ and at iteration. 
The fine propagator $\mathcal F_n$ is used in parallel over all time slices; the coarse propagator $\mathcal C_n$ is applied sequentially in each iteration.

\subsection{Background: micro-macro Parareal}
In micro-macro Parareal \cite{Legoll2013}, the coarse propagator does not act on the original state variable $u \in \mathbb R^d$ (micro state), but on a reduced version $\rho \in \mathbb R^r$ (macro state). 
The restriction operator $\mathcal R: \mathbb R^d \rightarrow \mathbb R^r$ extracts macro information from a micro state ($\rho = \mathcal{R}(u)$).
The lifting operator $\mathcal L: \mathbb R^r \rightarrow \mathbb R^d$ provides a micro state $u$ that is consistent with a given macro state $\rho$ ($u = \mathcal{L}(\rho)$).
The matching operator $\mathcal M: (\mathbb R^r, \, \mathbb R^d) \rightarrow \mathbb R^d$ produces a micro state $u$ that is consistent with a given macro state $\rho$, based on prior information about the micro state $\hat u$ ($u = \mathcal{M}(\rho, \hat u)$).

The iterations of micro-macro Parareal are defined as follows:
For the initialisation, we have $u^k_0 = u_0$ and $\rho^k_0 = \mathcal R (u_0)$, for all $k \geq 0$.
Then, for all $n \geq 1$ and $k \geq 1$ we have
\begin{equation}
\begin{aligned} 
\rho_{n+1}^{0}  &= \mathcal{C}_n (\rho_n^{0}), \qquad 
u_{n+1}^{0} &= \mathcal{L}(\rho_{n+1}^{0}),
\end{aligned}
\label{mM_parareal_0}
\end{equation}
\begin{equation}
\begin{aligned} 
\rho_{n+1}^{k+1} 
&= \mathcal{C}_n (\rho_n^{k+1}) 
+ \mathcal{R} ( \mathcal{F}_n (u_n^k))
- \mathcal{C}_n( \rho_n^k), \\     
u_{n+1}^{k+1} &= \mathcal{M}(\rho_{n+1}^{k+1}, \mathcal{F}_n (u_n^k)).
\end{aligned}
\label{mM_parareal_other_k}
\end{equation}
Classical Parareal corresponds to the case $\mathcal{R}(x) =x$ and $\mathcal{L}(x,y) = \mathcal{M}(x,y) = x$ \cite[Remark 10]{Legoll2013}. 

\begin{property}
\label{properties_micro_macro_Parareal}
We give two properties of micro-macro Parareal.
\begin{itemize}
\item (Consistency) If the coupling operators are chosen such that $\mathcal{M}(\mathcal{R}u,u) = u$ for all $u$, then it holds that $\rho_n^k = \mathcal{R}u_n^k$ for all $k$ and $n$ (see \cite[equation (3.18)]{Legoll2013}). 
\item (Exactness property) If the matching operator satisfies $\mathcal M (\mathcal R u, u) = u$, for any $u \in \mathbb{R}^d$ and $\mathcal R(\mathcal M(\rho,u)) = U$ for any $\rho$ and $u$, then the micro-macro Parareal algorithm satisfies the exactness property, namely $u^k_n = u_n$ for all $k \geq n$, \cite[Theorem 12]{Legoll2013} 
\footnote{The details of what happens in the zeroth iteration for $n \geq 1$ are irrelevant for the exactness property.}.
\end{itemize}
\end{property}

\subsection{Generalisation of micro-macro Parareal}
In this subsection we provide a framework to formulate and analyse an extension of the MC-moments Parareal algorithm.
Let  $\mathcal I: (\mathbb R^r, \mathbb R^d, \mathbb R^r) \rightarrow (\mathbb R^d, \mathbb R^r)$ be an iterator function.
The generalised algorithm reads, for all $k \geq 0$ and $n \geq 0$:
\begin{equation}
\begin{aligned} 
\left( \rho_{n+1}^{0}, u_{n+1}^{0} \right)  &= 
\left(
\mathcal{C}_n (\rho_n^{0}), \,
\mathcal{L}(
\mathcal{C}_n (\rho_n^{0})
) \right) \\
(\rho_{n+1}^{k+1} , u_{n+1}^{k+1})
&= 
\mathcal I \left(
\mathcal{C}_n (\rho_n^{k+1}), \, 
\mathcal{F}_n (u_n^k), \,
\mathcal{C}_n( \rho_n^k)
\right).
\end{aligned}
\label{mM_parareal_other_k_NEW}
\end{equation}
\begin{property}
\label{exactness_property_of_generalised_Parareal}
\label{testlineref}
If the iterator function satisfies $\mathcal I(a,b,a) = (\mathcal R(b), b)$ for all $a$ and $b$, and for all $n \leq k$, then
MC-moments Parareal satisfies
the exactness property $u^k_n = u_n$ whenever $k \geq n$ and also obeys micro-macro consistency, that is $\rho^k_n = \mathcal R (u^k_n)$ for all $k \geq 0$ and $n \geq 0$\footnote{Here as well, the lifting operator in the zeroth iteration is irrelvant for the finite-exactness property.}.
\end{property}
\begin{proof}
We use the technique of proof by induction (over $n$).
\begin{itemize}
\item For $n=1$ it holds that, for all $k \geq 1$, $\left( \rho_1^k, u^k_1 \right) 
= \mathcal I( 
\mathcal C_0 \rho_0^k, 
\mathcal F_0 u^{k-1}_0, 
\mathcal C_0 \rho^{k-1}_0) 
= (\mathcal R 
\mathcal F_0 u_0, 
\mathcal F_0 u_0)
= (\mathcal R u_1, u_1)$.

\item Then the induction step. 
Suppose that $u^k_n = u_n$ and $\rho^k_n = \mathcal R( u^k_n) = \rho_n $ are exact for all $k \geq n$. 
Thus, it also holds that $u^{k+1}_n = u_n$ and $\rho^{k+1}_n = \mathcal R( u^{k+1}_n) = \rho_n$.
Then, 
$(\rho_{n+1}^{k+1} , u_{n+1}^{k+1})
=
\mathcal I \left(
\mathcal{C}_n (\rho_n^{k+1}) 
, \mathcal{F}_n (u_n^k)
, \mathcal{C}_n( \rho_n^k)
\right) =
\left(
\mathcal{C}_n (\rho_n) 
, u_{n+1}
, \mathcal{C}_n( \rho_n)
\right)
= (\mathcal R(u_{n+1}), u_{n+1}).
$
Thus $\rho^{k+1}_{n+1}= \rho_{n+1}$ and $u^{k+1}_{n+1} = u_{n+1}$ are also exact for all $k \geq n$. 
\end{itemize}
\end{proof}

\begin{remark}[micro-macro Parareal revisited]
Micro-macro Parareal is a special case of the generalised method \eqref{mM_parareal_other_k_NEW} with the iterator function
\begin{equation}
\mathcal I 
\left(
\mathcal{C}_n (\rho_n^{k+1}) 
, \mathcal{F}_n (u_n^k)
, \mathcal{C}_n( \rho_n^k)
\right) 
=
\left( \rho^{k+1}, \mathcal M \left(
\mathcal{C}_n (\rho_n^{k+1}) 
+ \mathcal{R} ( \mathcal{F}_n (u_n^k))
- \mathcal{C}_n( \rho_n^k), 
\mathcal{F}_n (u_n^k) \right) 
\right)
\end{equation}
This is not the same iterator as the iterator defined in \cite{Legoll2020}, where the iterator is considered separately from the matching operator. 
Here, the iterator is defined more generally.
\end{remark}

\subsection{The MC-moments Parareal algorithm}
\label{section_about_MC_moments_Parareal}
In the design of a fast micro-macro Parareal algorithm for McKean-Vlasov SDEs, we require (i) the exactness property of Parareal, and (ii) a low computational cost of the coarse propagator and of the coupling operators.
The basic idea behind the MC-moments algorithm is to use a moment model from \cref{SETION_moment_equations} as a coarse propagator in Parareal, while the fine propagator is a particle approximation of the SDE. 

Before we present the actual MC-moments Parareal algorithm, we  first define two helper functions $\mathcal S$ and $\mathcal T$.
The function $\mathcal T$ transforms a particle ensemble $u$ to another particle ensemble with desired mean $\mu$ and variance $\sigma$. It is based on the transformation given in \cite[Section 4.3]{caflisch_monte_1998}:
\begin{equation}
\mathcal T \left( 
\left[ \mu, \, \sigma \right],  
u \right) 
= 
\begin{cases}
\sqrt{\sigma} \xi + \mu
	& \qquad \mathrm{if \, \,} \mathbb V[u] = 0, \\
\sqrt{ \frac{ \sigma}{\mathbb V[u]} }
\left(
u - \mathbb E [u]
\right) + \mu,
	& \qquad \mathrm{else},
\end{cases}
\label{multimodal_M}
\end{equation}
where $\xi$ are normally distributed random variables. The function $\mathcal S$ equals $\mathcal S (u) 
= 
\left[
\mathbb E[u] 
,\mathbb V [u] 
\right].$
\begin{lemma}[Consistency of the operators $\mathcal S$ and $\mathcal T$]
\label{consistency_of_operators_S_and_T}
It holds that (i)
$\mathcal T (\mathcal S u, u) = u$ for any $u$; and that (ii) $\mathcal S (\mathcal T(\rho, u)) = \rho$ for any $u$ and $\rho$ with a positive variance in $\rho$.
\end{lemma}

In the following definitions, we formulate the MC-moments Parareal algorithm for unimodal and multimodal McKean-Vlasov SDEs.
\begin{tcolorbox}[breakable,sharp corners,colback=black!1]
\begin{definition}[MC-moments Parareal for unimodal scalar McKean-Vlasov SDEs]
\label{definition_operators_unimodal} \, \\
MC-moments Parareal is extended Parareal \eqref{mM_parareal_other_k_NEW} with these choices:
\begin{itemize}
\item The micro state $u^k_n$ is an ensemble of $P$ particles.

\item The fine solver $\mathcal F_n$ is an Euler-Maruyama discretisation \cref{Euler_Maruyama}.

\item The macro state equals $\rho^k_n  = (M^k_n, \Sigma^k_n)$, where $M^k_n \approx \mathbb E[u^k_n]$ and $\Sigma^k_n \approx \mathbb V[u^k_n]$.

\item The coarse solver $\mathcal C_n$ is a Taylor series-based moment moment model \eqref{moment_model_class_1} or a Gaussian model \eqref{Gaussian_model_McK_V_SDE}. 

\item The restriction operator $\mathcal R$ is defined as
$
\mathcal R(u) 
= 
\mathcal S(u)
= 
\left[ \mathbb E_P[u], \mathbb V_P[u] \right].
$

\item Let $\alpha = a+\mathcal R(b)-c$ (the addition and substraction are understood componentwise) and let $\alpha_{\mathbb V}$ be the second (variance) component of $\alpha$, then the iterator function is defined as
\begin{equation}
\mathcal I
\left( 
a, b, c
\right)
 = \begin{cases}
(\alpha, \, \mathcal T(\alpha, b)), \qquad 
	& \mathrm{if \, \,} \alpha_{\mathbb V} \geq 0, \\
(\mathcal R (b), \, b) \qquad 
	& \mathrm{\text{else.}} \\
\end{cases}
\label{iterator_function_unimodal}
\end{equation}

\item The lifting operator $\mathcal L$ is a transformation of the initial condition: 
$
\mathcal L(\rho) = \mathcal{T}(\rho, u_0).
\label{lifting_operator_unimodal}
$
\end{itemize}
\end{definition}
\end{tcolorbox}

\noindent
\begin{tcolorbox}[breakable,sharp corners,colback=black!2]
\begin{definition}[MC-moments Parareal for multimodal scalar McKean-Vlasov SDEs]
\label{MC_moments_Parareal_for_multimodal}
MC-moments Parareal is extended Parareal \eqref{mM_parareal_other_k_NEW} with these choices:
\begin{itemize}
\item The micro state $u^k_n$ is an ensemble of $P$ particles. 

\item The fine solver $\mathcal F_n$ is an Euler-Maruyama discretisation \cref{Euler_Maruyama}.

\item The macro state equals $\rho^k_n  = (\mathcal P_{\mathcal D_i}, M^k_{i,n}, \Sigma^k_{i,n})_{i=1}^I$, where $M^k_{i,n}$ and $\Sigma^k_{i,n}$ are approximations of the mean and variance of the particle ensemble $u^k_n$ in region $\mathcal D_i$ and $\mathcal P_{\mathcal D_i}$ is an approximation of the fraction of particles residing in region $\mathcal D_i$.

\item Let $\mathbb{E}_{\mathcal D_i}$ and $\mathbb{V}_{\mathcal D_i}$ denote the mean and variance of a subset of those particles that reside in the region $\mathcal D_i$ 
and let $\mathcal P_{\mathcal D_i}$ denote the (empirical) probability that particles reside in the domain $\mathcal D_i$, then the restriction operator equals
\begin{equation}
\mathcal R (X) 
:= 
\left[
\begin{matrix} 
\mathcal R_1 \\
\hdots  \\
\mathcal R_I \\
\end{matrix}
\right]
=
\left[
\begin{matrix} 
\mathcal P_{\mathcal D_1}(X)
	& \mathbb E_{\mathcal D_1}[X] 
	& \mathbb V_{\mathcal D_1} [X]  \\
\hdots  \\
\mathcal P_{\mathcal D_I}(X)
	& \mathbb E_{\mathcal D_I}[X]
	& \mathbb V_{\mathcal D_I} [X]
\end{matrix}
\right].
\label{multimodal_R}
\end{equation}

\item The coarse propagator $\mathcal C_n$ is, for each region $\mathcal D_i$
\begin{itemize}
\item on $M^k_{i,n}$ and $\Sigma^k_{i,n}$ the system of ODEs \eqref{moment_model_Taylor_multimodal}  is applied
\item on $\mathcal P_{ \mathcal{D}_i,n}^k$, no update is applied
\end{itemize}

\item Let $\alpha = a + \mathcal R_i(b) - c$ (the addition and substraction are understood componentwise) and let $\alpha_{\mathbb V}$ be the second (variance) component of $\alpha$.
If $\mathcal P_{i}(b) = 0$ and $n > k$, let $\alpha = a$.
Then, for each region of attraction $\mathcal D_i$, the iterator function $\mathcal I_i$ is defined as, 
\begin{equation}
\mathcal I_i
\left( 
a, b, c
\right)
 = \begin{cases}
(\alpha, \, \mathcal T(\alpha, b)) 
	\qquad 
	& \alpha_{\mathbb V} \geq 0 \mathrm{\, and \, \,} P_{\mathcal D_i}(\mathcal T (\alpha, b)) = P_{\mathcal D_i}(b), 
 	\\
(\mathcal R (b), \, b) \qquad 
	& \mathrm{\text{else, or if \,}} P_{\mathcal D_i}(b) = 0. \\
\end{cases}
\label{iterator_function_multimodal}
\end{equation}
Let $\mathcal R_P$ denote an operator that computes the particle fractions
of a particle ensemble.
The particle fractions are updated as 
$\mathcal P_{\mathcal D_i, n+1}^{k+1} 
= 
\mathcal R_P \left(
\mathcal F_n \left( u^k_n \right)
\right)$.

\item Lifting operator $\mathcal L$: for each region $\mathcal D_i$, move particles:
$
\mathcal L(\rho_i) = \mathcal T(\rho_i, u(0)_{\mathcal D_i})
$.
\end{itemize}
The global mean and variance are approximated by using \eqref{mixture_of_distributions_mean_variance} on $M_{i,n}^k$ and $\Sigma_{i,n}^k$.
\end{definition}
\end{tcolorbox}

For unimodal SDEs, the two cases in the definition of the iterator are required to avoid negative variance, which would be nonphysical (and equation \eqref{multimodal_M} is undefined in such cases).

\begin{property}[Exactness property of unimodal and multimodal MC-moments Parareal]
\label{exactness_property_MC_moments}
If the same random numbers are used in each iteration, MC-moments Parareal satisfies the exactness properties $u^k_n = u_n$ and $\rho^k_n = \mathcal R u_n$ for all $k \geq n$ (the particle ensemble converges exactly to the sequential simulation).
The strong and weak approximation errors are zero for all $k \geq n$.
\end{property}
\begin{proof}
The unimodal and multimodal iterator functions \eqref{iterator_function_unimodal} and \eqref{iterator_function_multimodal} satisfy $\mathcal I(a,b,a) = \left( \mathcal R(b), b \right)$ for all $a,b$ and all $n \leq k$. 
Thus, as a result of \cref{exactness_property_of_generalised_Parareal}, the exactness property holds.
\end{proof}

\begin{figure}[h]
\centering
\captionsetup{width=.9\linewidth}
\includegraphics[width=0.5\textwidth]{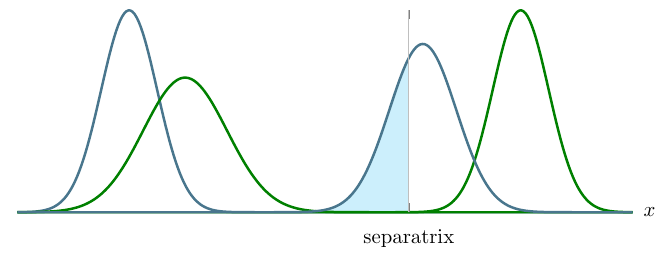}
\caption{Illustration of a potential pitfall with particle transformations. 
The regions of attraction are separated by a separatrix. 
The green distributions are the locally unimodal distributions in each region of attraction before matching. 
The cyan ones result after application of the matching operator, but this changes the particle fractions.}
\label{Illustration_particle_fractions_changed_after_matching}
\end{figure}

\begin{remark}[About the update of the particle fractions]
\label{remark_particle_fractions_without_coarse}
In \cref{MC_moments_Parareal_for_multimodal}, the
iterator function only carries out the transformation $\mathcal T$ if $P_{\mathcal D_i} (\mathcal T (\alpha, b)) = P_{\mathcal D_i} (b)$.
That is, particles are only moved if they stay inside the same region of attraction
before and after their transformation.
In other words, the particle fractions are not
affected by the iterator $\mathcal I_i$. 
Moreover, it holds that
$\mathcal P_{\mathcal D_i, n+1}^{k+1} = 
\mathcal R_P \left( 
\mathcal F_n \left( u^k_n \right)
\right) $.
The particle fractions are thus not affected by
the coarse Parareal propagator, but only by the fine Parareal propagator. 
We
now explain the reason for this design choice.

Let us consider the situation in \cref{Illustration_particle_fractions_changed_after_matching}, where a particle transformation $\mathcal T$ affects the particle fractions.
The transformation did not succeed 
to give the locally unimodal subensembles the desired (local) means and variances.
One could, instead, devise another (iterative) method that redistributes particles until the local subensembles satisfy the desired means and variances.
Yet, in the scope of this work we decided not to do so. 
Instead, we only perform matching in a region $\mathcal D_i$ if the transformation step does not modify the particle fractions in any region. 
Otherwise, we only accept the effect of the fine propagator.

Thus, the MC-moments Parareal algorithm leaves the particle fractions untouched, unless for the effect of the fine propagator $\mathcal F_n$.
The particle fractions converge as in Parareal without coarse propagator. For more information about this algorithm, see \cite[lemma 5]{bossuyt_new} (Dahlquist equation) or \cite[theorem 1]{gander_parareal_2024} (parabolic equations).
\end{remark}

\subsection{Other possible choices, high-level comparison} 
\mylabel{section_other_choices_comparison}
In \cref{MC_moments_Parareal_for_multimodal}, we choose the coarse propagator to be an ensemble of independent pairs of Taylor-based ODEs for the local means and variances in each nnoverlapping region of attraction for multimodal SDEs.
A model for the particle fractions is difficult to obtain because it requires information on the exit times (this usually requires sequential simulations).
In the next subsection, we present a learning-based variant (\cref{definition_MC_moments_Parareal_with_learning}) to tackle this challenge.

Another choice for the coarse propagator is a Gaussian mixture (GM).
A GM describes the evolution of a mixture of Gaussians that each have a weight $w_i$, mean $M_i$ and variance $\Sigma_i$. 
This choice, however, would require different micro-macro Parareal coupling operators. 
Then, the restriction operator maps a particle ensemble onto a mixture of Gaussians. 
This can be done, for instance, using methods like those described in \cite{kalai_disentangling_2012}.
The matching operator then transforms a particle ensemble such that its restriction (a GM fit) equals a desired state. 
A candidate method for this task consists of, 
for instance, 
computing a histogram (on a predefined grid) of the particle ensemble and of the Gaussian mixture, 
then applying the Parareal update in histogram space, 
and then somehow resampling the resulting histogram, see also \cite{Legoll2020}.
This procedure is, however, not straightforwardly generalisable to higher dimensions. 
The latter would require, for instance, the solution to an optimal transport problem, see \cite{Legoll2020}, and is more computationally expensive than the simple particle transformation $\mathcal T$.

In \cref{table_comparison_different_MC_Parareals} we list some possible alternatives and briefly indicate some (dis)advantages of each choice.

\begin{table}[h]
\centering
\captionsetup{width=.9\linewidth}
\begin{tabular}{c|c|c|c|c}
& \textit{coarse propagator}
	& \textit{restriction} 
	& \textit{matching} 
	& \textit{(dis)advantages} \\ \hline

1 &\eqref{Gaussian_model_McK_V_SDE},~\eqref{moment_model_class_1}
	& \eqref{multimodal_R} 
	& \eqref{multimodal_M} 
	& cheap \\

& & & & requires RoA \\
	
& & & & simple update \\

& \eqref{moment_model_Taylor_multimodal}
	&
	&
	& MPF nontrivial \\

\hline
2 
	& Gaussian mixture	
	& fitting a GM 
	& e.g. via histogram
	& coupling operators \\
	
& & & & more expensive \\
\end{tabular}
\caption{Comparison of moment models and potential options for corresponding restriction and matching operators. MPF=model for particle fractions. GM=Gaussian mixture. RoA=Regions of Attraction.}
\label{table_comparison_different_MC_Parareals}
\end{table}

\subsection{Extension: learning-based (MC-moments) Parareal}
\mylabel{subsection_learning_based_Parareal}
The convergence of the Parareal algorithm without coarse propagator may be very slow, depending on how fast the particle fractions converge (see \cref{remark_particle_fractions_without_coarse}). 
In this section, we develop a new MC-moments Parareal algorithm where the evolution of the particle fractions is learned as the Parareal iterations progress. 

\paragraph{Model for particle fractions}
In \cite[Chapter XI, equation 1.4]{kampen_stochastic_1981} a model is presented for the evolution of the particle fractions $\mathcal P_{\mathcal D_1}$ and $\mathcal P_{\mathcal D_2}$, with $i=1,2$, of a (classical) bimodal SDE:
\begin{equation}
\mathcal P_{\mathcal D_i}(t) 
= 
e^{
\alpha t
}
\left[ 
\mathcal P_{\mathcal D_i}(0) + 
\beta_i
\right]
- 
\beta_i,
\label{exponential_for_particle_fractions}
\end{equation}
where the parameters $\alpha$ and $\beta_i$ are related to exit times of particles from one mode to another.
More generally, let $\omega_{ij}$ be a parameter related to the exit times of particles from region $\mathcal D_i$ to $\mathcal D_j$.
In \cref{appendix_model_particle_fractions} we derive a generalisation of \eqref{exponential_for_particle_fractions} to multimodal SDEs:
\begin{equation}
\frac{d}{dt} 
\left[ 
\begin{matrix}
\mathcal P_{\mathcal D_1} \\
\mathcal P_{\mathcal D_2} \\
\vdots  \\
\mathcal P_{\mathcal D_I}
\end{matrix}
\right]
=
\underbrace{
\left[ 
\begin{matrix}
-\sum_{j \neq 1} \omega_{1j} 
	& \omega_{21} 
	& 
	& \hdots
	& \omega_{I1} \\ 
\omega_{12} 
	& -\sum_{j \neq 2} \omega_{2j} 
	& 
	& \hdots
	& \omega_{I2} \\ 
\hdots & \hdots & & & \vdots \\	
\omega_{1I} 
	& \omega_{2I}
	& 
	& \hdots
	& - \sum_{j \neq I} \omega_{Ij}  \\ 
\end{matrix}
\right]}_{A}
\left[ 
\begin{matrix}
\mathcal P_{\mathcal D_1} \\
\mathcal P_{\mathcal D_2} \\
\vdots  \\
\mathcal P_{\mathcal D_I}
\end{matrix}
\right],
\label{generalised_model_multimodal_PFs_repetition}
\end{equation}
with the initial condition $
\left[ 
\begin{matrix}
\mathcal P_{\mathcal D_1}(0) &
\mathcal P_{\mathcal D_2}(0) &
\hdots  &
\mathcal P_{\mathcal D_I}(0)
\end{matrix}
\right]
= 
\left[ 
\begin{matrix}
\mathcal P_{\mathcal D_1, 0} &
\mathcal P_{\mathcal D_2, 0} &
\hdots  &
\mathcal P_{\mathcal D_I, 0}
\end{matrix}
\right]
$.
In practice, however, this model is not practical since (an approximation of) the exit times usually requires sequential simulations, which we intend to avoid.

\paragraph{New algorithm}
We now present a new Parareal variant, fitting in the framework of generalised Parareal, where a model for the particle fractions is learned by solving a least squares problem, fitting the parameters $\theta = \{ \omega_{i,j} \}_{i=1, j\neq i}^J$ in the model for the evolution of the particle fractions 
to data collected from the fine propagator.

Let the operator $\mathcal Q$ be a nonlinear least-squares solver that takes as input time-series data of particle fractions, and as output the fitted parameters $\theta$, using the Levenberg-Marquardt algorithm implemented in. 
Let the operator $\mathcal H_i$ be defined such that $\mathcal P_{i,n+m} \approx \mathcal H_i (\mathcal P_{i,n},m,\theta)$\footnote{For a bimodal SDE, the operator $\mathcal H_i$ equals $\mathcal H_i(\mathcal P_{i,n}, m, \theta) = e^{
\alpha m \Delta t 
}
\left(
\mathcal P_n + 
\beta_i
\right)
- 
\beta_i$.}. 
Then, we define the learning-based MC-moments Parareal algorithm as follows:

\begin{tcolorbox}[breakable,sharp corners,colback=black!2]
\begin{definition}[MC-moments Parareal with learning]
\label{definition_MC_moments_Parareal_with_learning}
MC-moments Parareal with learning is the MC-moments Parareal algorithm with the following modifications: 
\begin{itemize}
\item At the end of each iteration, the learning operator $\mathcal Q$ updates the parameters: $\theta = \mathcal Q(\{ u^k_n \}_{n=1}^k)$.

\item For all $n \geq k$, all particle fractions $\hat{\mathcal P}_{i,n}^k$ are computed by extrapolation with the learned parameters: 
$\hat{\mathcal P}_{i,n}^k 
= 
\mathcal H_i (\mathcal P_{i,k}^k, k-n, \theta)$.

\item For all $n \geq k$, the global mean and variance are approximated by using
 \eqref{mixture_of_distributions_mean_variance} on $M_{i,n}^k$ and $\Sigma_{i,n}^k$ and the extrapolated $\hat{\mathcal P}_{i,n}^k$.
If, for $n+1 \geq k$, no particles in the micro variable $u^k_n$ are present in region of attraction $\mathcal D_i$, then the variables $M_{i,n+1}^{k+1}$ and $\Sigma_{i,n+1}^{k+1}$ are approximated by applying (sequentially) the coarse propagator: 
$(M_{i,n+1}^{k+1}, \Sigma_{i,n+1}^{k+1}) = \mathcal C_{\mathcal D_i} (M_{i,n}^{k+1}, \Sigma_{i,n}^{k+1})$.
\end{itemize}
\end{definition}
\end{tcolorbox}

\section{Analysis of MC-moments Parareal}
\label{SECTION_analysis}
In this section we provide a short convergence analysis of the proposed algorithm.

\subsection{Sources of approximation errors in the MC-moments Parareal algorithm}
\paragraph{Approximation errors for Euler-Maruyama discretisations}
Let $X_n^{(p)}$ be a path obtained through numerical simulation, and $X(t_n)^{(p)}$ be its exact solution. 
The approximation error of the Euler-Maruyama estimator $M_{Q,P,n}$, which uses Monte Carlo sampling on \eqref{Euler_Maruyama}, consists of a bias component (the time-stepping error that would arise with an infinite number of samples), and a statistical error (arising from the effects of using only a finite number of samples).

More precisely, 
\begin{equation}
\begin{aligned}
&\left| M_{Q,P,n} - M_{\infty,\infty}(t_n) \right| 
\leq 
\left| M_{Q,P,n} - M_{\infty,\infty,n} \right|
+ 
\left| M_{\infty,\infty,n} - M_{\infty,\infty}(t_n)  \right|.
\end{aligned}
\label{splitting_of_stochastic_error}
\end{equation}

Now assume an infinite number of particles, such that the second term in the right-hand side of \eqref{splitting_of_stochastic_error} vanishes.
The strong error $E_n$ at time $t_n$ concerns the expected error on individual realisations of the stochastic process.
Let $\Phi$ be a user-chosen function. 
The weak error $e_{\Phi,n}$, on the other hand, concerns the error on averaged quantities:
\begin{equation}
\begin{aligned}
E_n 
&=  
\max_{1 \leq p \leq P} \mathbb{E}[|\bar X_n  - \bar X(t_n)|], \\
e_{\Phi,n} 
&= 
| \mathbb{E}[\mathbb{E}_P[\Phi(\bar X_n)] - \mathbb{E}[\mathbb{E}_P [\Phi(\bar X(t_n))] ] ] |.
\end{aligned}
\label{definition_weak_strong_error}
\end{equation}

\paragraph{Choice of discretisation parameters and termination condition of the MC-moments algorithm}
In the proposed algorithm, there exist multiple sources of approximation errors: 
(i) the statistical error in $\mathcal F_n$ which unavoidably occurs in the stochastic, Monte Carlo approximation method, 
(ii) the time-discretisation error (bias) through the temporal discretisation in $\mathcal F_n$ and $\mathcal C_n$ and 
(iii) the model error introduced through the moment model in $\mathcal C_n$.
A proper implementation ensures that the statistical error is of the same order as the time-stepping error in $\mathcal F_n$. 
A termination condition for MC-moments Parareal then  stops the iteration procedure as soon as the model error-induced error from $\mathcal C_n$ on the Parareal iterates reaches this level. 
An actual implementation of this heuristic is left for future work.

\subsection{The effect of statistical error on the convergence}
In this section, we first present an error recursion for MC-moments Parareal.
The subsequent linear and superlinear convergence bounds follow directly from the analysis in \cite{gander_analysis_2007}. 
Here we explicitly analyse the case of affine ODEs while \cite{gander_analysis_2007} considered only linear ODEs. 

\begin{lemma}
\label{lemma_error_recursion_with_noise}
Let $F$ and $G$ be scalars, and let $\epsilon$ be a random variable with  $\expect{|\epsilon|} \leq \varepsilon$ and let $\epsilon^k_n$ denote a random draw from $\epsilon$ in the $k$-th Parareal iteration at timepoint $n$. 
Let $\mathcal C_n (u) = Gu + g_n(t)$, with $|G| < 1$, be a coarse propagator, and let $\mathcal F_n (u) = Fu + f_n(t) + \epsilon$ be a fine propagator for the classical Parareal algorithm \eqref{original_parareal}.
Let the reference solution $u$ satisfy $u(t_{n+1}) = \mathcal F_n(u(t_n))$ for all $n \geq 0$.
Then, the expected value of the approximation error, $\xi^k_n = \expect{|u^k_n - u_n|}$, satisfies
the recursion
\begin{equation}
\xi^{k+1}_{n+1} 
\leq
|F-G| \xi^k_n + |G| \xi^{k+1}_n  + \varepsilon.
\label{equation_error_recursion_with_noise}
\end{equation}
\end{lemma}
\begin{proof}
From the definition of the Parareal algorithm \eqref{original_parareal},
it holds that
\begin{equation}
\begin{aligned}
u^{k+1}_{n+1} 
&= \left( F u^k_{n} + f + \epsilon \right)
+ \left( Gu^{k+1}_{n} + g \right) 
- \left( Gu^k_n + g \right) \\
&= F u^k_{n} + \epsilon^{k+1}_{n+1}
+ Gu^{k+1}_{n} 
- Gu^k_n 
\end{aligned}
\end{equation}
For the approximation error, defined as $e^k_n = u^k_n - u_n$, it holds that
\begin{equation}
\begin{aligned}
e^{k+1}_n 
&= u^{k+1}_{n+1} - u_{n+1} \\
&= F u^k_{n} + \epsilon^{k+1}_{n+1} + Gu^{k+1}_{n} - Gu^k_n - u_{n+1} \\
&= F(u^k_{n} - u_n) + G(u^{k+1}_n - u_n) - G(u^k_n - u_n) + \epsilon^{k+1}_{n+1} \\
&= (F-G) e^k_n + Ge^{k+1}_n + \epsilon^{k+1}_{n+1}
\end{aligned} 
\label{equation_starting_point_first_proof}
\end{equation}
Now, taking the absolute value and then the expectation, leads to \eqref{equation_error_recursion_with_noise}.
\end{proof}
In practice, $\epsilon$ corresponds to a statistical error on the fine solver. 
The quantity $\mathbb E[|\epsilon|]$ can for instance be bounded through $\expect{|\epsilon|} \leq \sqrt{\expect{\epsilon^2}}$ (see \cite[equation (4.1)]{engblom_parallel_2009}).
Now we are in a position to formulate our main theoretical result, which is proven in \cref{proof_of_main_lemma}.

\begin{lemma}[Linear and superlinear bound for Parareal with noise]
\mylabel{lemma_what_happens_with_noise}
Let $\xi^k_{\mathrm{max}} = \max_n \expect{|u^k_n - u_n|}$.
Consider the same setting as \cref{lemma_error_recursion_with_noise}, then $\xi^k_{\mathrm{max}}$ satisfies a superlinear bound:
\begin{equation}
\xi^k_{\mathrm{max}}
\leq 
\frac{\left| F - G \right|^k}{k!}
\prod_{j=1}^k (N-j) 
\xi^0_{\mathrm{max}}
+
\varepsilon 
\frac{1-|G|^N}{1-|G|}
\sum_{j=0}^{k-1} 
|F-G|^{j} \binom{N-1}{j}
\label{superlinear_bound_with_noise}
\end{equation}
and a linear bound:
\begin{equation}
\xi^k_{\mathrm{max}}
\leq 
\left( \frac{\left| F - G \right|}{1 - | G|} \right) ^k
\xi^0_{\mathrm{max}}
+
\varepsilon 
\frac{1}{1-|G|}
\sum_{j=0}^{k-1} \left(\frac{|F-G|}{1-|G|} \right)^{j}.
\label{linear_bound_with_noise}
\end{equation}
\end{lemma}
A similar bound, derived in the context of a low-rank Parareal method, and by using the technique of generating functions, is given in 
\cite[Lemma 2, Theorem 3, Theorem 4]{carrel_low-rank_2023}.
In words, the exactness property is satisfied up to some statistical error.

\begin{example}
We consider a linear ODE where noise is added to the fine propagator with various levels of intensity.
We use the following operators.
\begin{itemize}
\item A fine propagator $\mathcal F(u) = Fu + 2 \beta \zeta$ where $\zeta$ is a uniformly distributed random variable, such that $\mathbb E[|\zeta|] = \beta$. 
\item A coarse propagator $\mathcal C(u) = Cu$.
\item The coupling operators are the same as for classical Parareal. In that case, letting $I$ denote the identity, $\mathcal R = I$ and $\mathcal M = I$. This corresponds to MC-moments Parareal in the case of a single-particle ensemble.
\end{itemize}
In \cref{figure_effect_statistical_error} we display the convergence history of Parareal, as well as the linear and superlinear bounds from \cref{lemma_what_happens_with_noise}.
The bounds capture the evolution of the error well as the iterations proceed.
\begin{figure}[h]
\centering
\includegraphics[width=0.5\textwidth]{{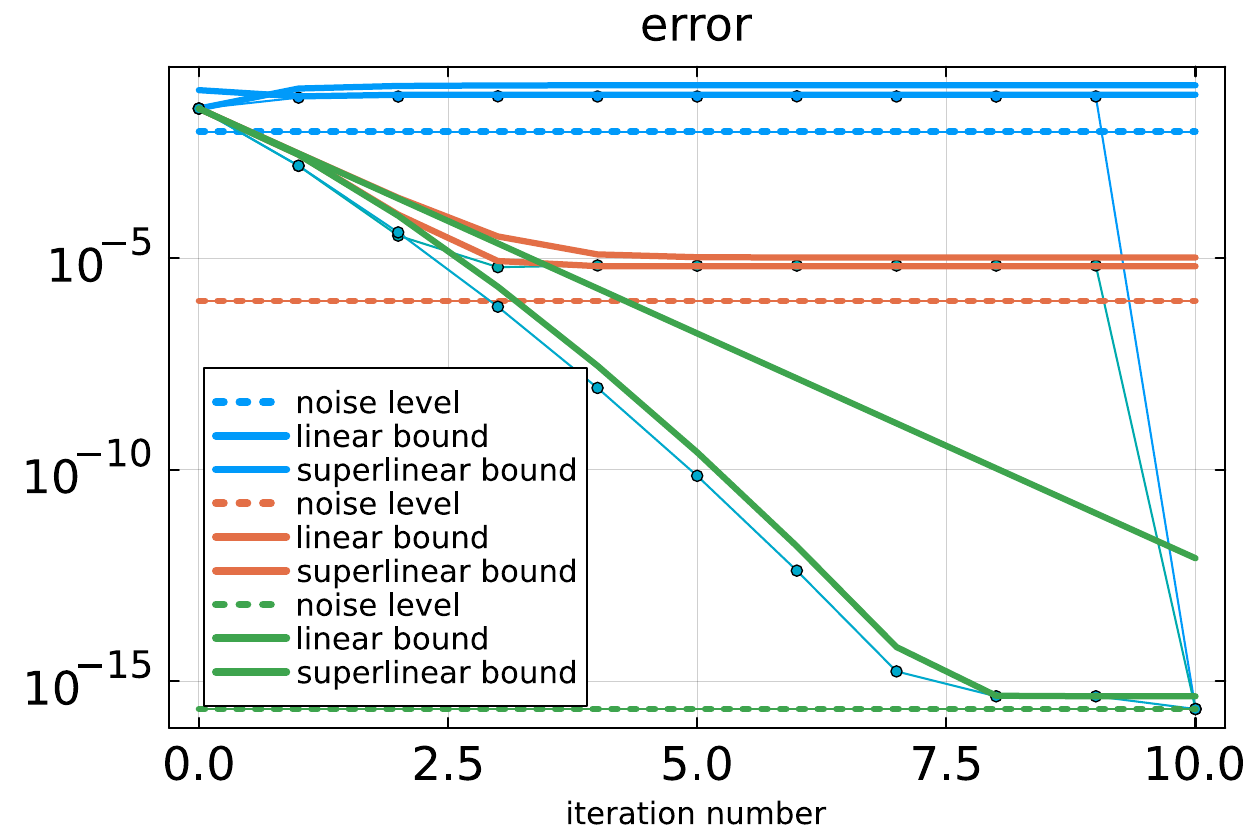}}
\caption{Illustration of the convergence bound for a simple test problem. The reference solution equals the sequential solution using the fine propagator $\mathcal F_n$ with $\beta=  0$. 
(Blue) $\beta = 10^{-2}$, (orange) $\beta= 10^{-6}$ and (green) $\beta = 0$.}
\label{figure_effect_statistical_error}
\end{figure}
\end{example}

\subsection{Convergence of micro-macro Parareal on the scalar Ornstein - Uhlenbeck SDE with model error}
\mylabel{subsecton_about_OU}

In this subsection, we study convergence of the MC-moments Parareal algorithm for a linear Ornstein-Uhlenbeck model with a perturbed coarse propagator.
In this section, we consider the Ornstein-Uhlenbeck SDE, which is  probably the simplest SDE, like the Dahlquist equation is for ODEs.

\begin{example}[label=exampleperturbedOU, name=Perturbed generalised Ornstein-Uhlenbeck SDE]
\label{exampleperturbedOU_label}
For the fine Parareal solver $\mathcal F_n$, we use an Euler-Maruyama discretisation of the Ornstein-Uhlenbeck SDE
\begin{equation}
dX^{(p)} = (aX^{(p)} + a_M \mathbb{E}[X])dt + B dW.
\label{MC1D_moments1D_fine_solver_OU_SDE}
\end{equation}
Let $\varepsilon_E$ and $\varepsilon_V$ be artificial perturbation parameters, then the coarse solver $\mathcal C_n$ consists of the ODE
\begin{equation}
\begin{aligned}
\frac{dM}{dt} &= (a+a_M)(1+\varepsilon_{M}) M, \\
\frac{d\Sigma}{dt} &= 2 a(1+\varepsilon_{M}) \Sigma + B^2(1+\varepsilon_{V})^2.
\end{aligned}
\label{MC1D_moments1D_coarse_solver_OU_ODE}
\end{equation}

We now consider the error on the mean and on the variance.
The mean and variance of the fine solver obey the ODEs
\begin{equation}
\begin{aligned}
\frac{dM_{\mathcal F}}{dt} 
&=
(a+ a_M) M_{\mathcal F} + \zeta_M, \\
\frac{d\Sigma_{\mathcal F}}{dt}
&=
2(a + a_E) \Sigma_{\mathcal F} + B^2.
\end{aligned}
\label{OU_evolution_mean_F}
\end{equation}
As $\varepsilon_{E}, \varepsilon_{V} \rightarrow 0$, the coarse model gets closer to the fine model.
Let $F$ be the result of one time step on the exact model for the mean and the variance of \eqref{OU_evolution_mean_F},
let $G$ be the result of one time step with the coarse solver \eqref{MC1D_moments1D_coarse_solver_OU_ODE}, and let $\epsilon$ be an upper bound on the statistical noise on the mean of the particle ensemble.
Then, \cref{lemma_what_happens_with_noise} can be used to bound the error of MC-moments Parareal.

For the numerical simulations, the chosen parameters are $\alpha	= -1$, $\beta = -0.5$, $\sigma = 1/100$ and $\{ X^{(p)}(0) \}_{p=1}^P = 100$, $\varepsilon_{E} = 1$ and $\varepsilon_{V} = 1$. 
The initial condition is chosen as $p_0(X) = \delta(X-100).$

In \Cref{fig_perturbed_OU_convergence}, we plot the convergence of the iterates with respect to the iteration number for various time horizons. 
We also plot the minimum of the linear and superlinear bounds from \cref{lemma_what_happens_with_noise}. For the variance, the bound overshoots the true error largely.

\begin{figure}[h]
\centering
\includegraphics[width=0.33\textwidth]{{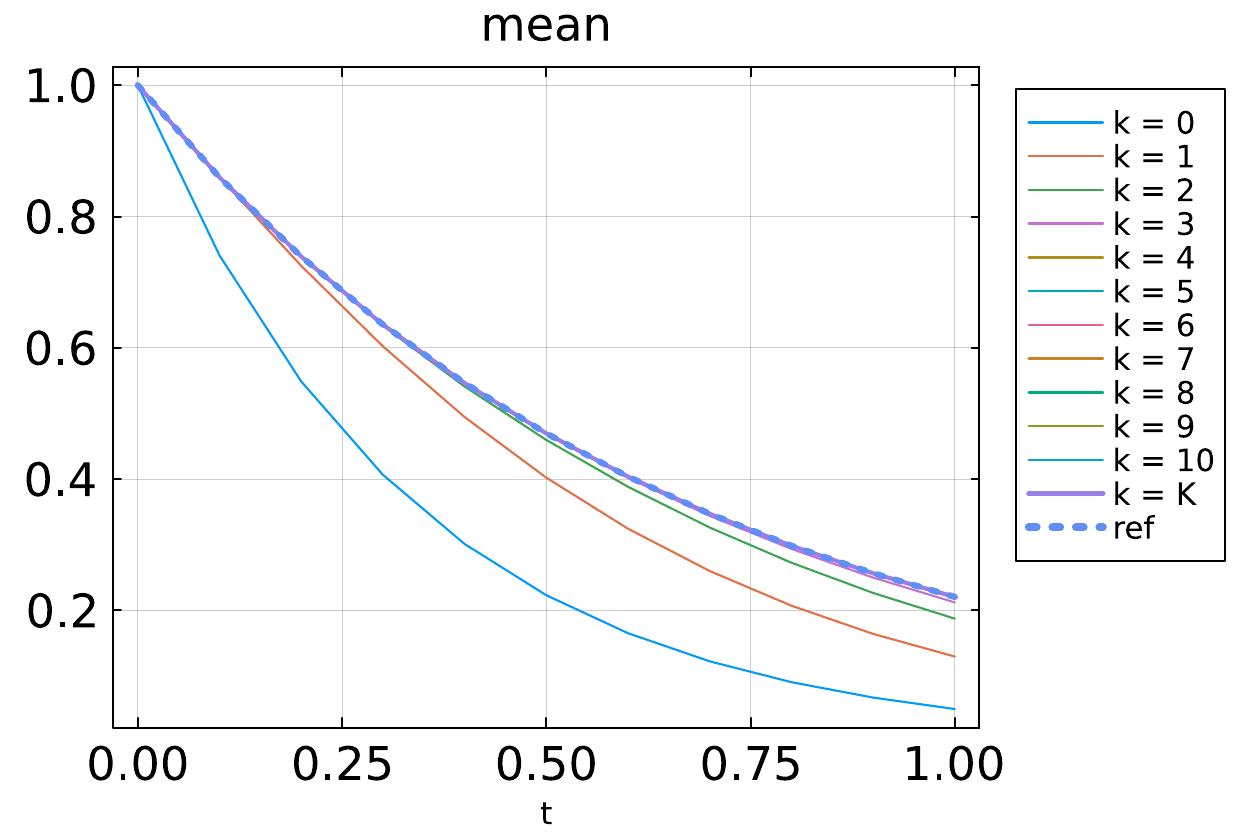}}
\includegraphics[width=0.33\textwidth]{{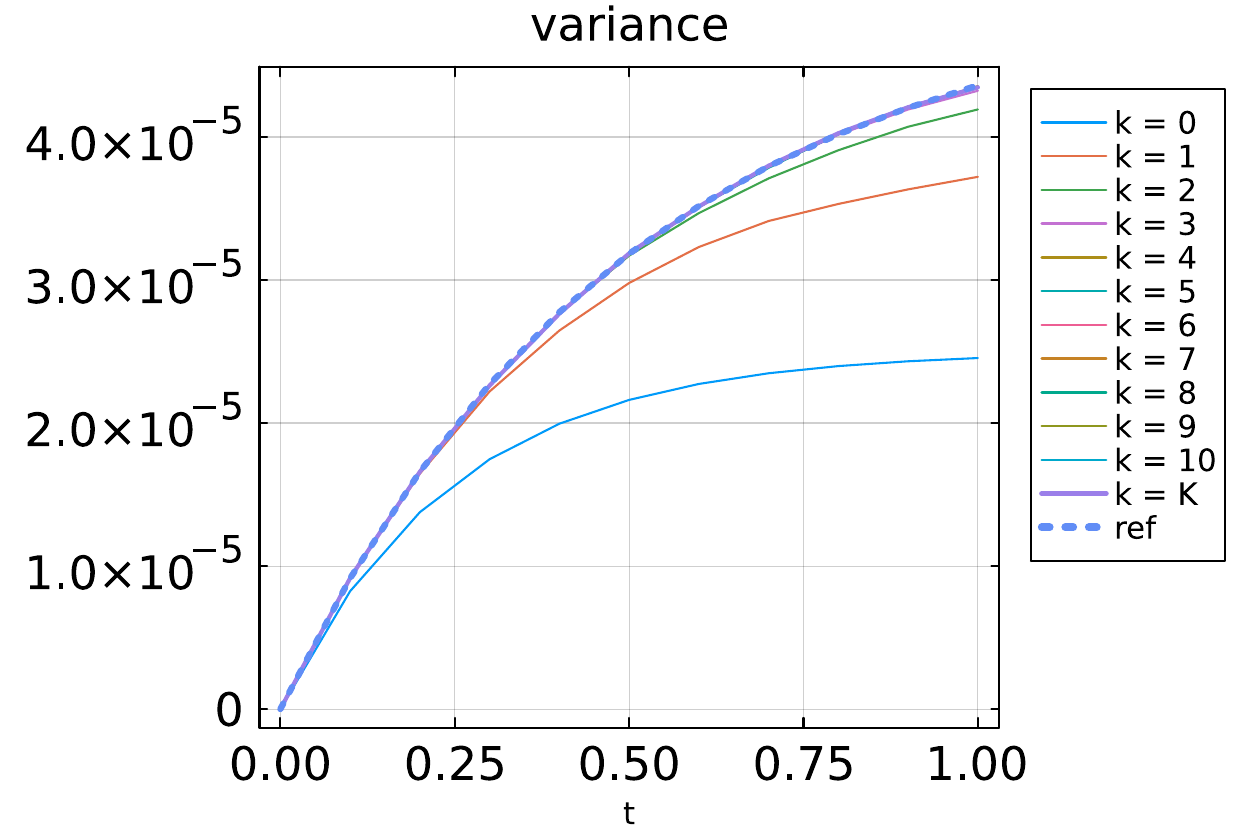}} \\
\includegraphics[width=0.33\textwidth]{{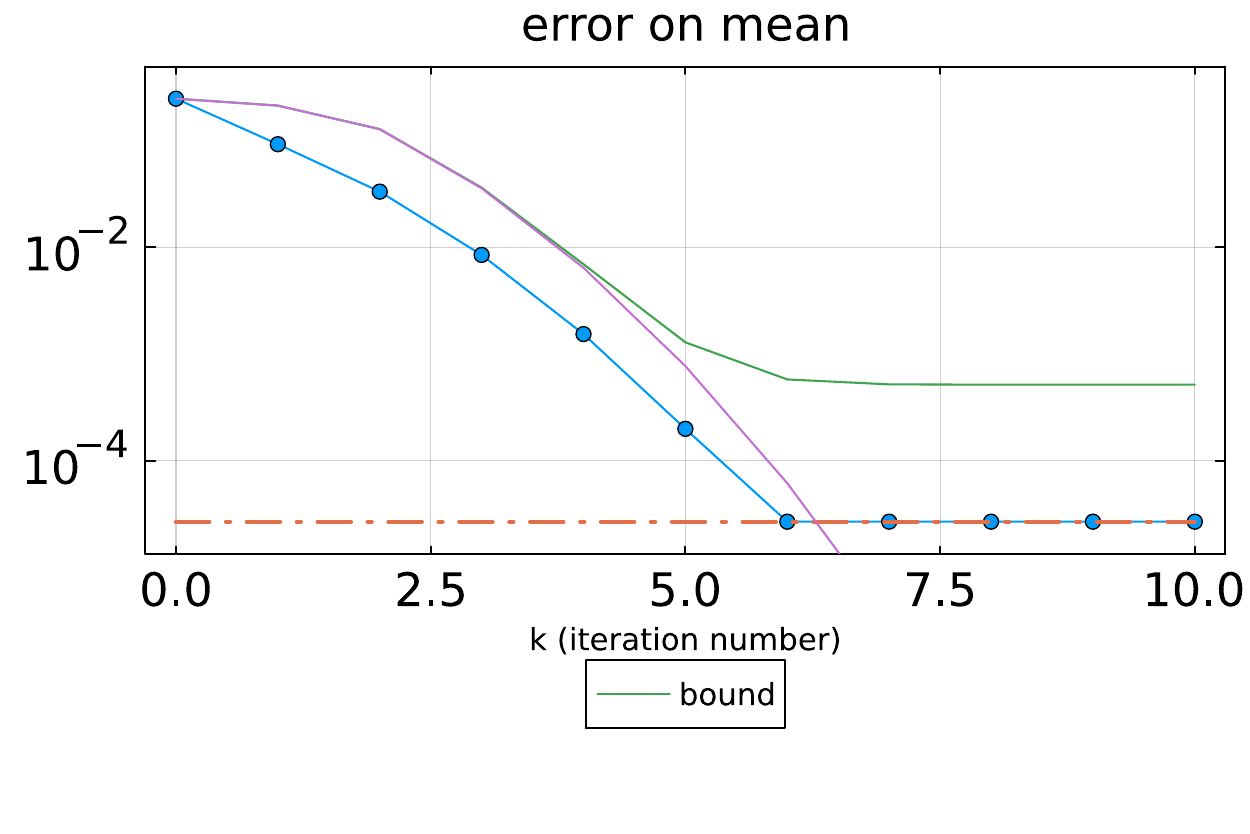}}
\includegraphics[width=0.33\textwidth]{{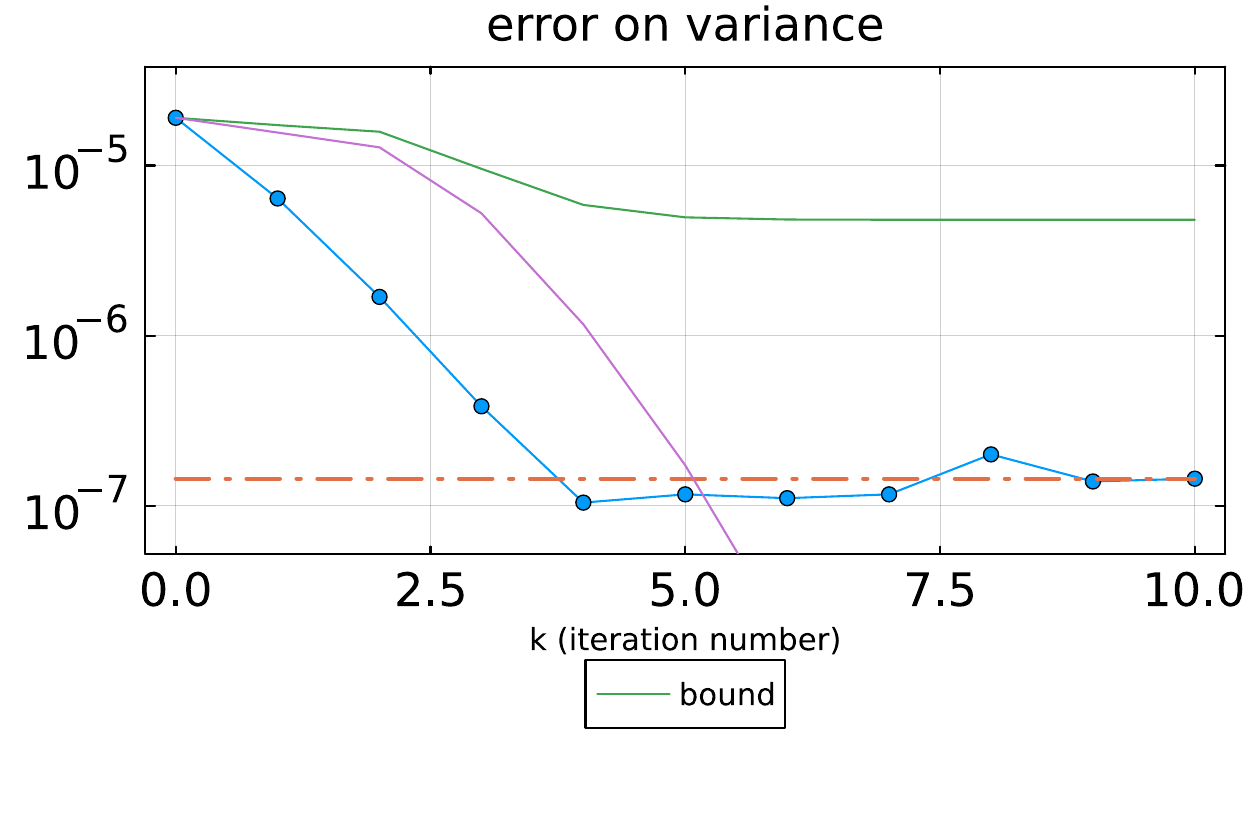}}
\caption{Convergence of (the error on) the mean and variance in the Parareal approximation of the Ornstein-Uhlenbeck SDE. 
We also plot the minimum of the linear and superlinear convergence bounds from \cref{lemma_what_happens_with_noise} with and without statistical noise included.}
\label{fig_perturbed_OU_convergence}
\end{figure}
\end{example}

\subsection{Classical Parareal, effect of time horizon}
\mylabel{Dahlquist_effect_time_horizon}
In this subsection, we study the effect of the end time $T$ on the convergence of Parareal, while keeping the number of subintervals $N$ constant.

\begin{example}
We use the following operators for classical Parareal:
\begin{itemize}
\item The fine solver is the exact solution of the Dahlquist equation $\exp(\Delta t \lambda)$ with $\lambda = -4$.

\item The coarse solver equals $\exp(\Delta t \lambda (1 + \epsilon))$ with model error $\epsilon = 0.1$.
\end{itemize}
In \cref{figure_Dahlquist_effect_time_horizon} (left) we show the convergence of Parareal for different time horizons $T$.
In the middle plot we show the sequential solution, computed with $\mathcal F_n$ as a function of time.
In the rightmost plot, we show the convergence of Parareal without coarse propagator (except in the zeroth iteration, where the coarse solver is applied sequentially).
If the time horizon is long, then even this variant converges relatively quickly.
\begin{figure}[h]
\centering
\includegraphics[width=0.33\textwidth]{{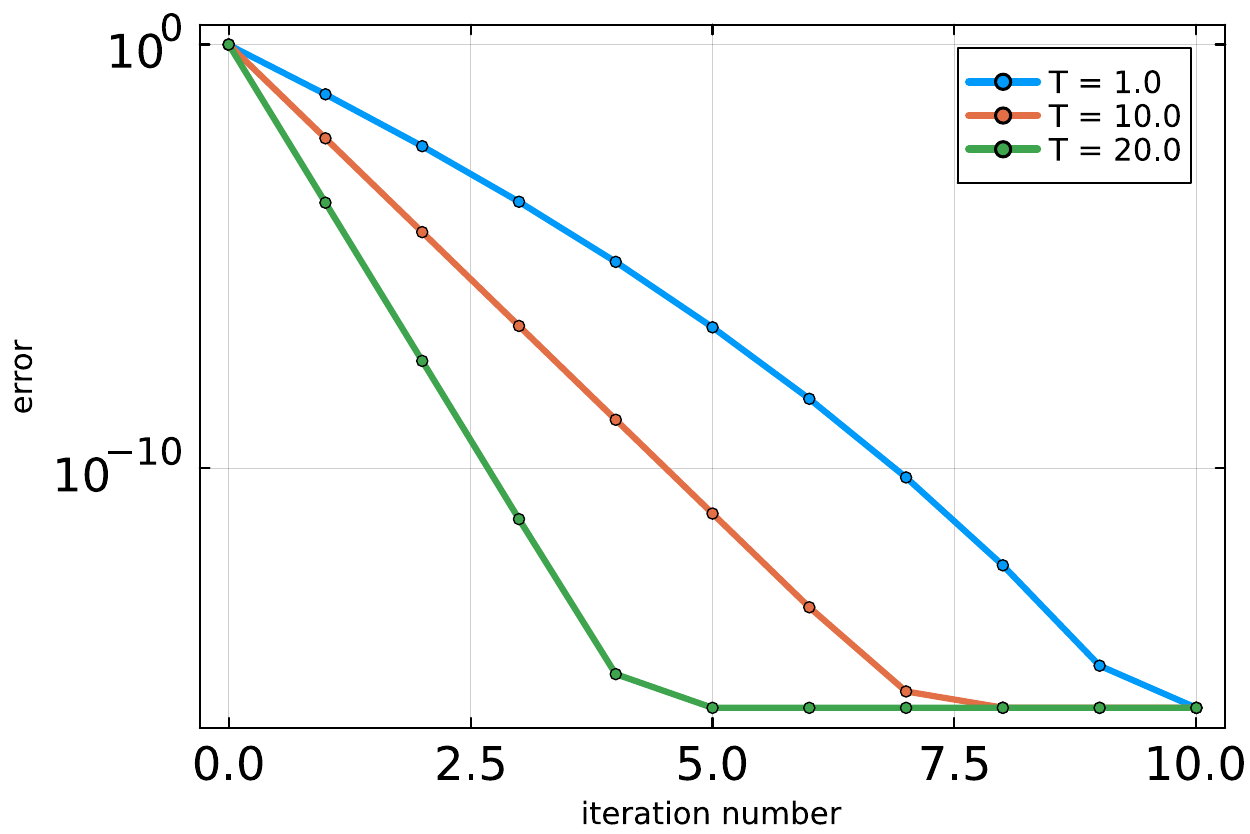}}
\includegraphics[width=0.33\textwidth]{{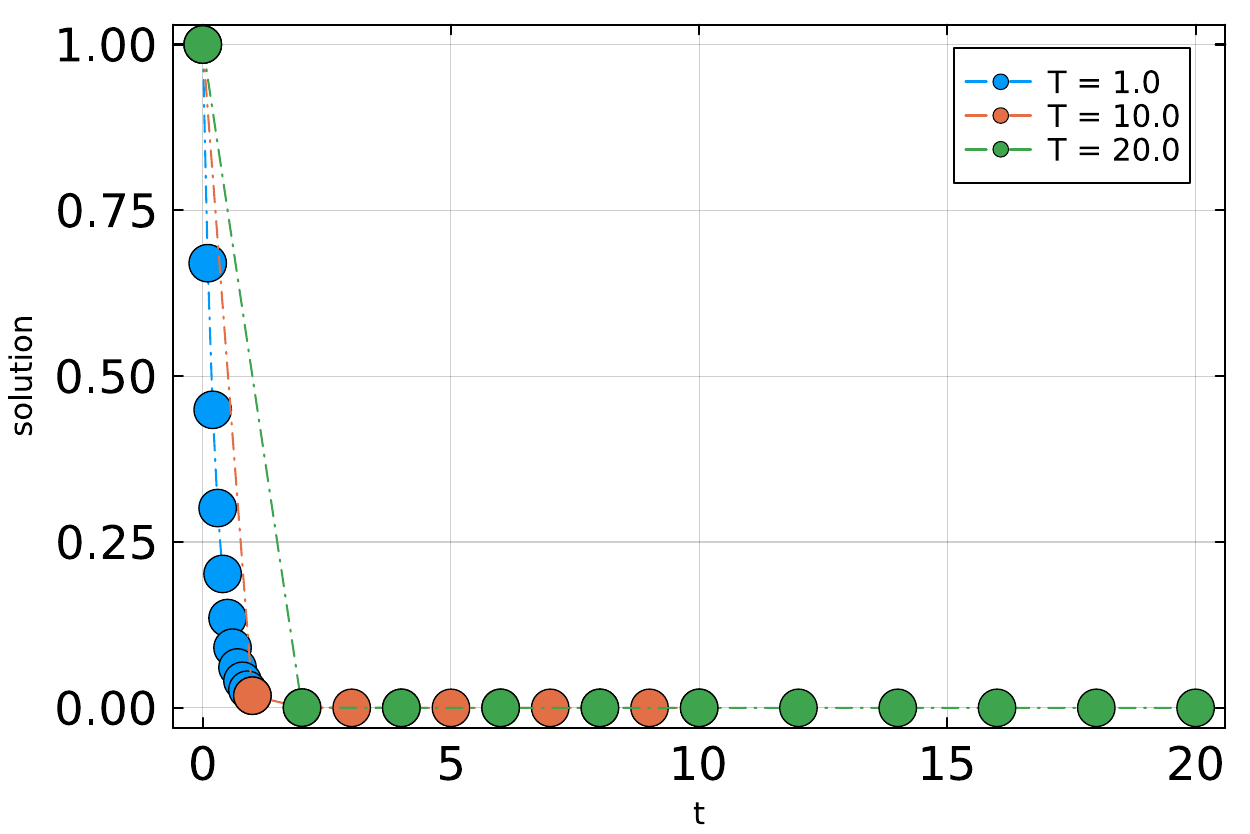}}%
\includegraphics[width=0.33\textwidth]{{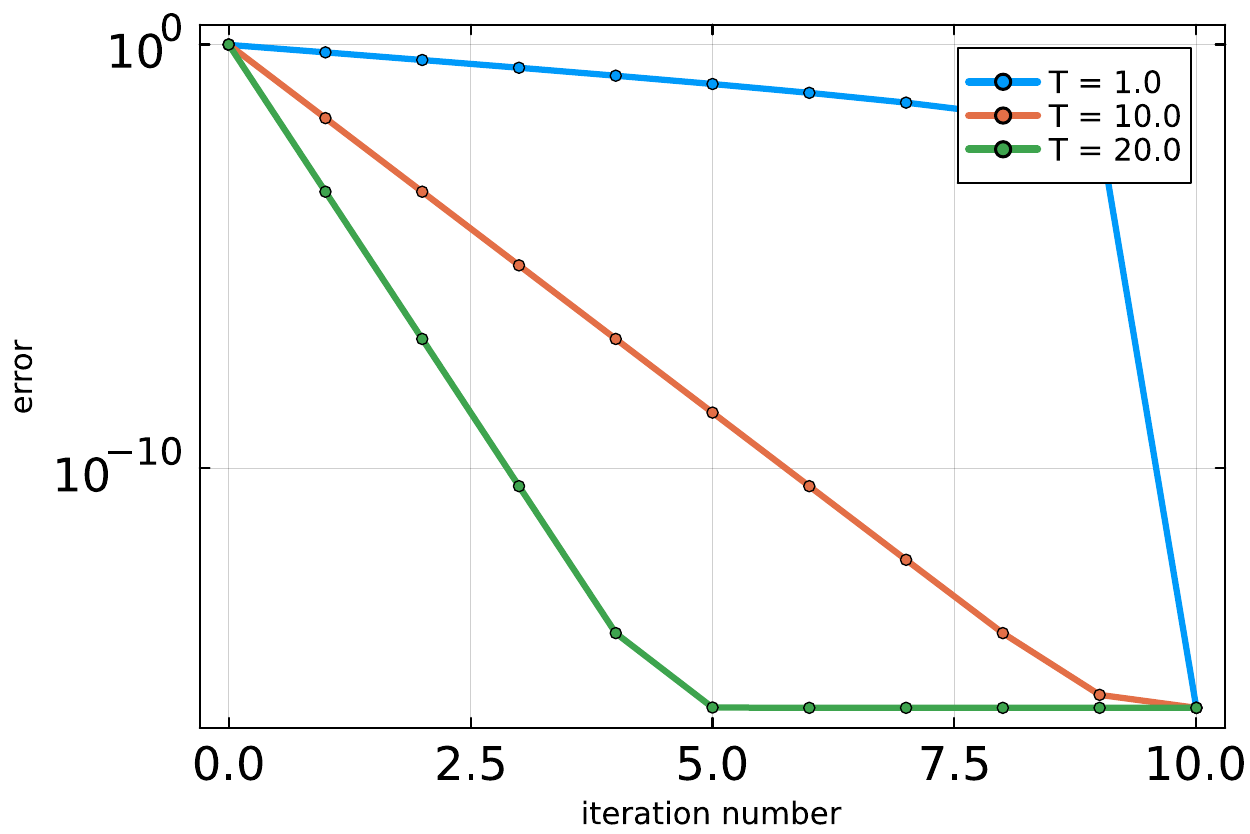}}
\caption{Convergence of classical Parareal with varying time horizon. (left) Convergence of Parareal: error as function of iteration number, (middle) Plot of the sequential solution for different time horizons, (right) Convergence of Parareal without coarse propagator.}
\label{figure_Dahlquist_effect_time_horizon}
\end{figure}
\end{example}

For simulations with longer time horizon, faster convergence takes place. Intuitively, this can be explained by observing that (i) the Dahlquist equation possesses a steady-state, and (ii) on longer time  horizons, the solution resides a larger fraction of the time in this steady state. The time-simulation of the steady solution is not very demanding, at least when the criterion is absolute error over the time domain.

\section{Numerical experiments}
\label{SECTION_numerical_exp}

In this section, we test our algorithm by means of numerical experiments. Some example setups are taken from \cite{kloeden_gauss_quadrature_2017} and also appear in  \cite{jourdain_2019}.

The code with numerical experiments (implemented sequentially) can be found in \cite{ignace_software_paper} as well as in the Supplementary Materials. 
For the MC discretisation of the McKean-Vlasov SDE, we implemented our own Euler-Maruyama code. 
For the Parareal approximations, the number of particles is $P = 10^5$. 
Reference solutions (based on the sequential application of the fine propagator) are computed with an ensemble of $10^6$ particles.
We numerically study the approximation error in the first two statistical moments on the Parareal approximations,
with respect to the statistical moments of the reference solution. These simulations are repeated 10 times in order to mitigate the effects of statistical outliers.
The moment ODEs are solved with the Tsitouras 5/4 Runge-Kutta method with an automatic stepsize controller, implemented in \cite{rackauckas_differentialequationsjl_2017}. 
For each ensemble, we use the same random numbers in each iteration. 

We use $N=10$ Parareal subintervals.
Let $\| \cdot \|_{n,\infty}$ denote the $\infty$-norm over all time points $n$.
We consider the maximum of the weak error on the mean and the variance:
\begin{equation}
\begin{aligned}
E_{\text{mean},k} 
&= 
\|  \mathbb{E}_P[\bar X_n^k] - \mathbb{E}_P[\bar X_n] \|_{n,\infty} 
\\ 
E_{\text{var},k} 
&= \|  \mathbb{V}_P[\bar X_n^k] - \mathbb{V}_P[\bar X_n] \|_{n,\infty}.
\end{aligned}
\label{DEFINITION_ALL_ERRORS}
\end{equation}

\paragraph{Overview of numerical experiments}
In \cref{subsection_unimodal_SDEs} we numerically study the convergence of the unimodal version of the MC-moments Parareal algorithm.
In \cref{Subsection_burgers} we consider the Burgers equation, which is unimodal but has a complex mean-field coupling.
In \cref{numerical_simulations_bimodal_SDEs} we present a detailed numerical convergence study of MC-moments Parareal for bimodal SDEs. 
In \cref{numerical_experiments_effect_time_horizon} we study the effect of the length of the time window on the convergence of MC-moments Parareal.

\subsection{Convergence of MC-moments Parareal for some unimodal SDEs}
\mylabel{subsection_unimodal_SDEs}
\begin{remark}[About MC-moment Parareal for linear McKean-Vlasov SDEs]
\label{MC_moments_Parareal_for_linear_SDEs}
Moment models for linear systems provide an exact description of the mean and variance.
One reason why MC simulations of linear systems are still interesting (and thus also the MC-moments Parareal algorithm) 
is that
a moment model alone is in general not able to describe the time evolution of the $\mathbb E_P[\Phi(X^{(p^)}(t))]$ if there is a nonlinear QoI $\Phi$.
A technique for computing an approximation to $\mathbb E[\Phi (X)]$ is to compute the empirical mean after the elementwise application of $\Phi$ on $\mathcal T([M, \, \Sigma], u_0)$\footnote{Alternatively, it is possible to estimate the quantity of interest $\mathbb E[\Phi(X)]$ as $\Phi(M) + \frac{\Sigma^2}{2} \Phi''(M)$ (see \cref{appendix_estimating_nonlinear_QoI}). Also, a Gaussian quadrature rule can be applied.}.
\end{remark}

\begin{example}[Linear Geometric Brownian Motion SDE with nonlinear QoI]
\mylabel{example_GBM_with_Parareal}
This example is taken from \cite{kloeden_gauss_quadrature_2017}.
It consists of a geometric Brownian motion with an extra interaction term, and a nonlinear QoI $\Phi(x) = \sin(x^2)$:
\begin{equation}
dX = (\alpha X + \beta \mathbb E_P[X] + a) dt + (\sigma X + \omega \mathbb E_P[X] + b)dW.
\end{equation}
with $a = -1/2$,  $\beta = 4/5$, $a= 1/3$, $\sigma  = \omega = \frac{1}{2}\sqrt{\frac{1}{2}}$ and $b=1/6$ and $p_0 = \mathcal N(1, 1/4)$.

In \cref{fig_GBM_convergence_nonlinear}, we show that
the two methods from \cref{MC_moments_Parareal_for_linear_SDEs} for the computation of a nonlinear QoI based on a moment model, produce very different results (in the left panel). 
We also show (in the middle panel) the convergence of MC-moments Parareal on the nonlinear QoI, for different coarse solvers: (i) the Taylor series-based SDE, and Gaussian approximations with (ii) two sigma-points and (iii) three sigma-points.
In the zeroth iteration, empirical mean of the QoI of the lifted ensemble computed.
In the right panel of \cref{fig_GBM_convergence_nonlinear}, we show a histogram of the particles at time $t=0$ and $t=T$.
The MC-moments Parareal variants with the three different coarse propagators all have a very similar convergence behavior.
 
\begin{figure}[h]
\centering
\includegraphics[width=0.33\textwidth]{{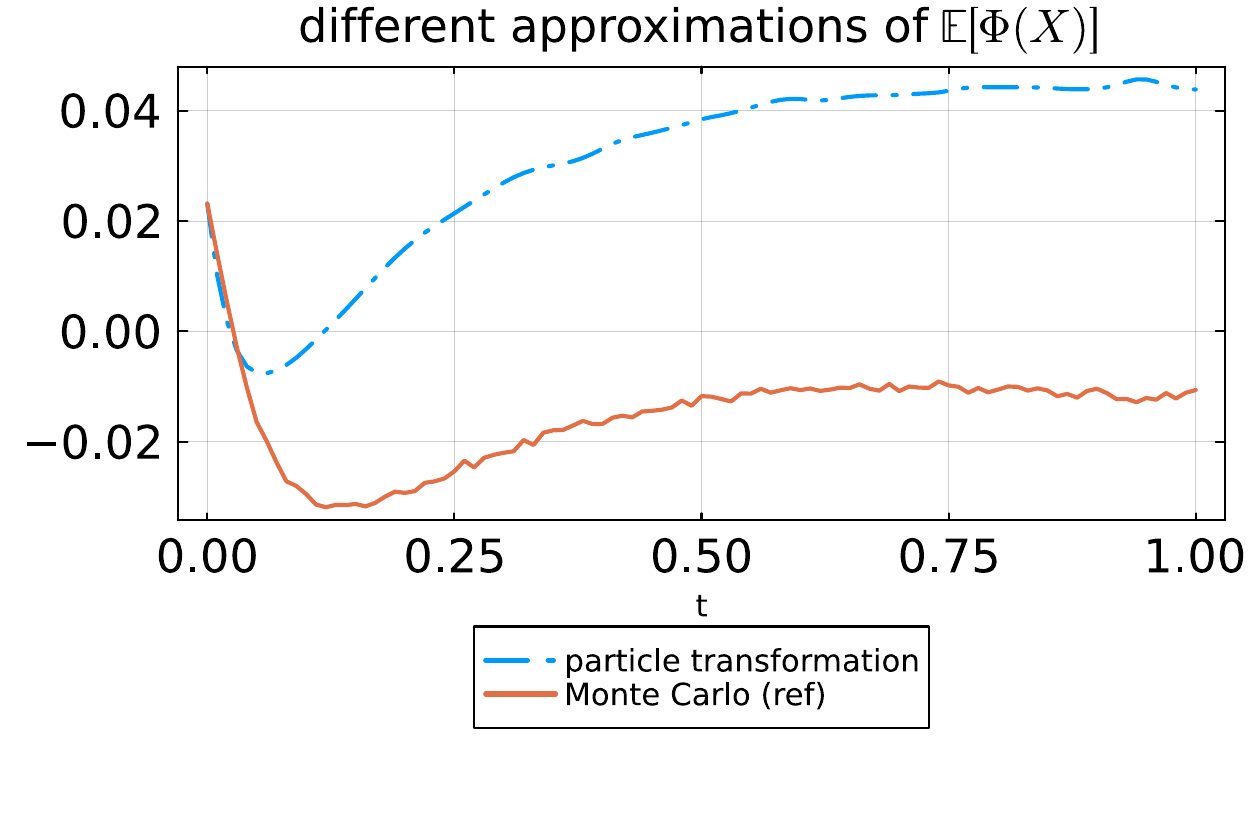}}%
\includegraphics[width=0.33\textwidth]{{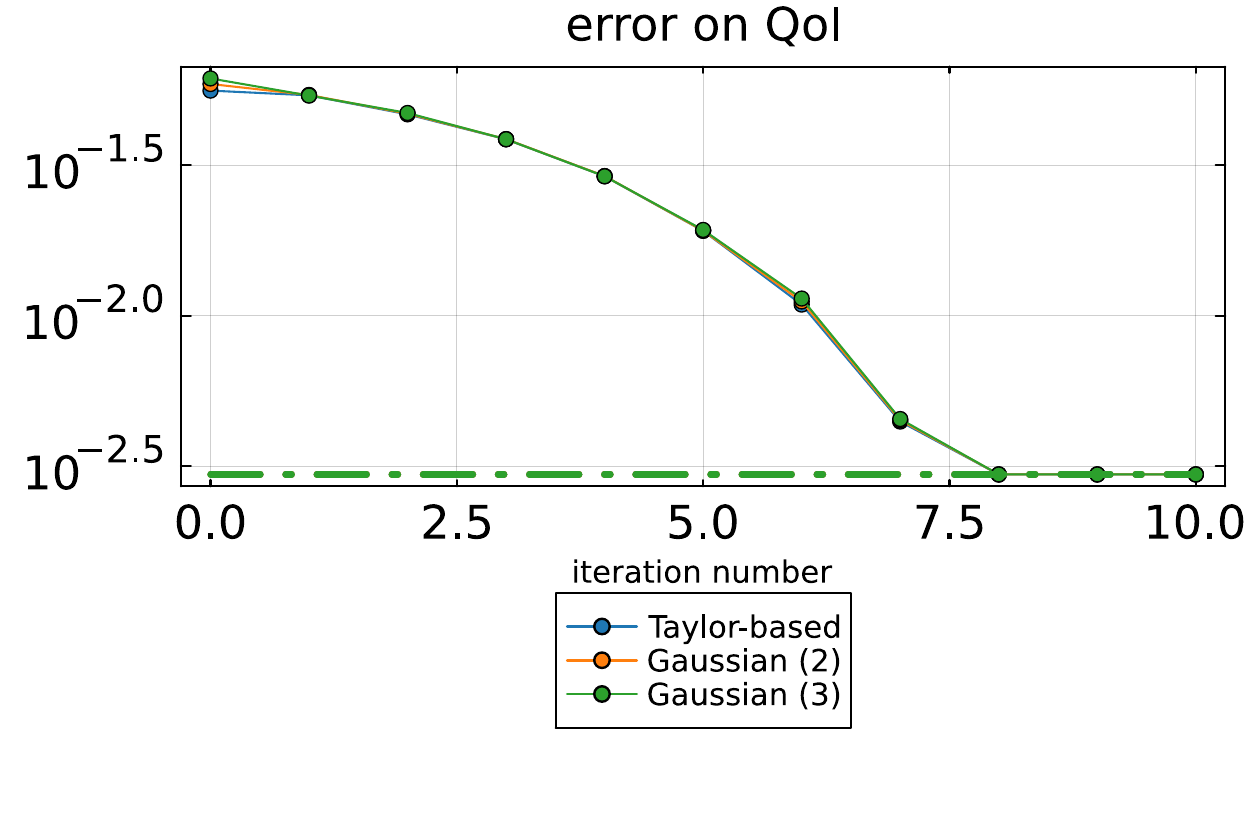}}%
\includegraphics[width=0.33\textwidth]{{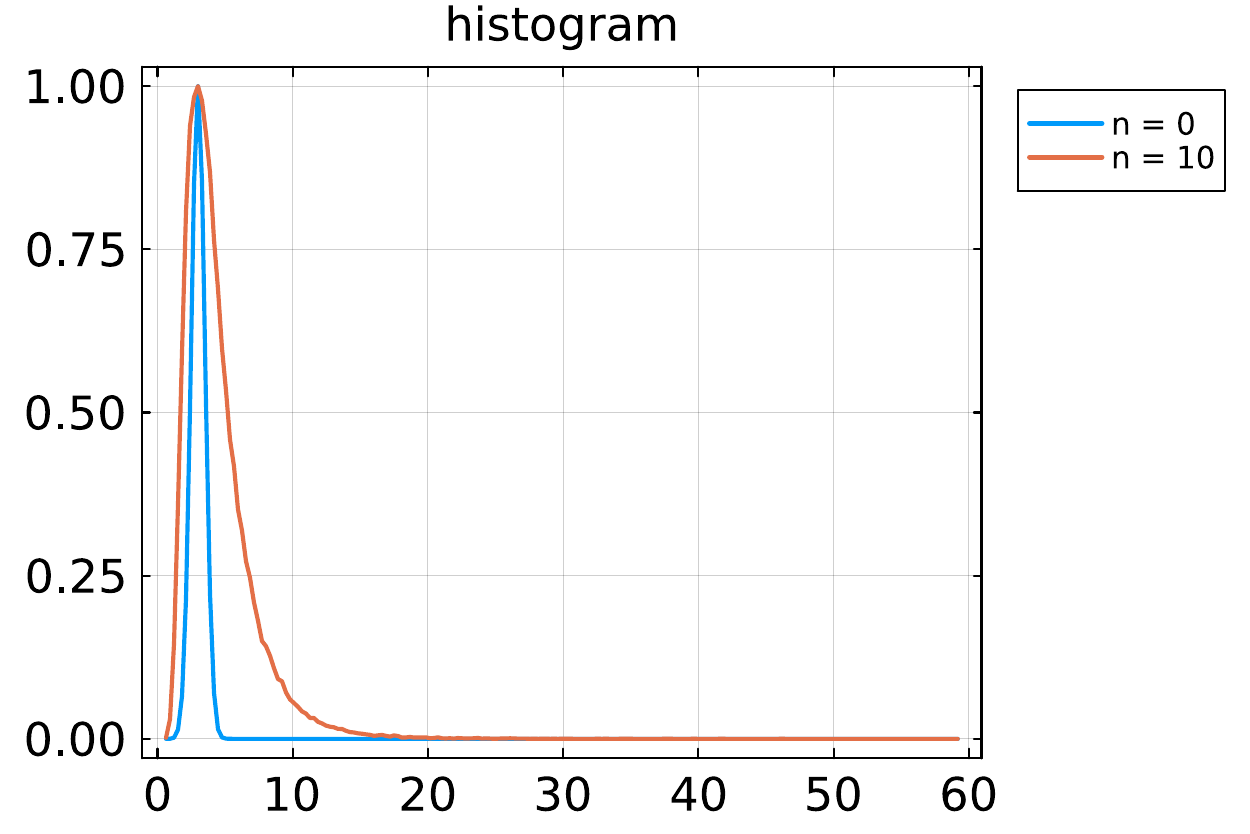}}
\caption{Modified geometric Brownian motion: (left) approximation of a nonlinear QoI using two approximation methods, (middle) MC-moments Parareal convergence of the weak error for the QoI with three variants, (right) histogram at the beginning and at the end of the time window.}
\label{fig_GBM_convergence_nonlinear}
\end{figure}
\end{example}

\begin{example}[Unimodal plane rotator]
\mylabel{example_plane_rotator}
We consider the unimodal plane rotator, taken from \cite{kloeden_gauss_quadrature_2017} and \cite{kostur_nonequilibrium_2002}:
\begin{equation}
\begin{aligned}
dX^{(p)} &=  \left[ \frac{K}{P} \sum_{q=1}^P \sin(X^{(q)} - X^{(p)}) - \sin(X^{(p)}) \right]dt + \sqrt{2 k_B T} dW ^{(p)}, 
\qquad 
X^{(p)}(0) &\sim \mathcal N \left( \frac{\pi}{4}, \frac{3 \pi}{4} \right),
\end{aligned}
\label{SDE_plane_rotator}
\end{equation}
where every particle $X^{(p)}$ is shifted modulo $2\pi$ into the interval $[0, 2\pi]$.
In \cite{kloeden_gauss_quadrature_2017} a shift modulo $2\pi$ is applied in each time-step. This corresponds to shifting at each timestep every particle $X^{(p)}$ modulo $2\pi$ into $[0, 2\pi]$ in the Euler-Maruyama algorithm. This creates a unimodal distribution of the solution particles. If no such shift is applied, the SDE \eqref{SDE_plane_rotator} possesses a multimodal particle distribution.  

The moment model \eqref{moment_model_class_1} equals 
\begin{equation}
\begin{aligned}
\frac{dM}{dt} 
&= 
-\sin(M) 
+ \frac{\Sigma}{2}\sin(M)  \\
\frac{d\Sigma}{dt},
&= 
-2 \left( -K - \cos(M) \right) \Sigma 
+ \sigma^2.
\end{aligned}
\label{plane_rotator_third_order}
\end{equation}

We choose the discretisation parameters $[0,T] = [0, 10]$ and $\Delta t = 1/1000$. 
We illustrate in \Cref{fig_plane_rotator_Taylor} 
the convergence of MC-moments Parareal, for four different coarse solvers:
\label{fig_plane_rotator_TaylorLINE}
\begin{itemize}
\item the moment model \eqref{plane_rotator_third_order} without the last term in the evolution of the mean

\item the moment model \eqref{plane_rotator_third_order}

\item a Gaussian approximation \eqref{Gaussian_model_McK_V_SDE} with 2 sigma-points 

\item a Gaussian approximation \eqref{Gaussian_model_McK_V_SDE} with 3 sigma-points 
\end{itemize}

\begin{figure}[h]
\includegraphics[width=0.32\textwidth]{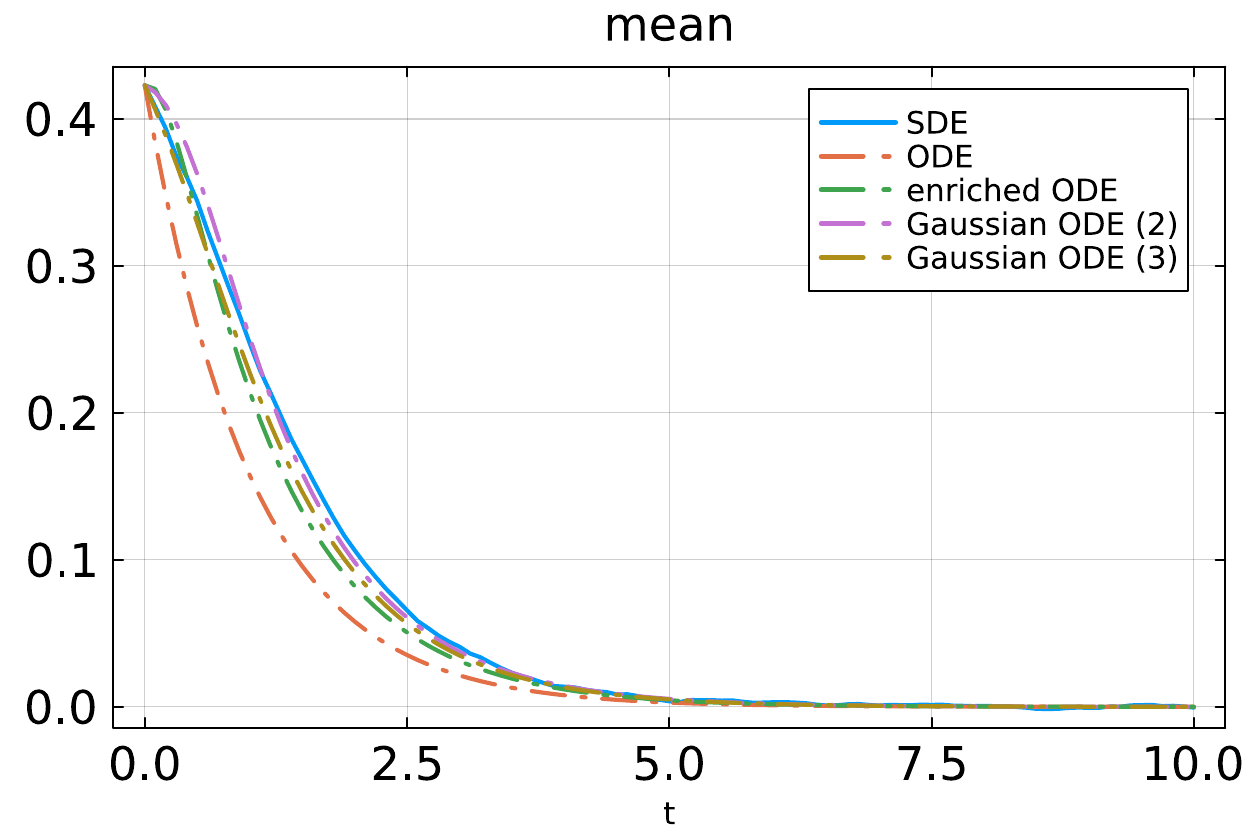}
\includegraphics[width=0.32\textwidth]{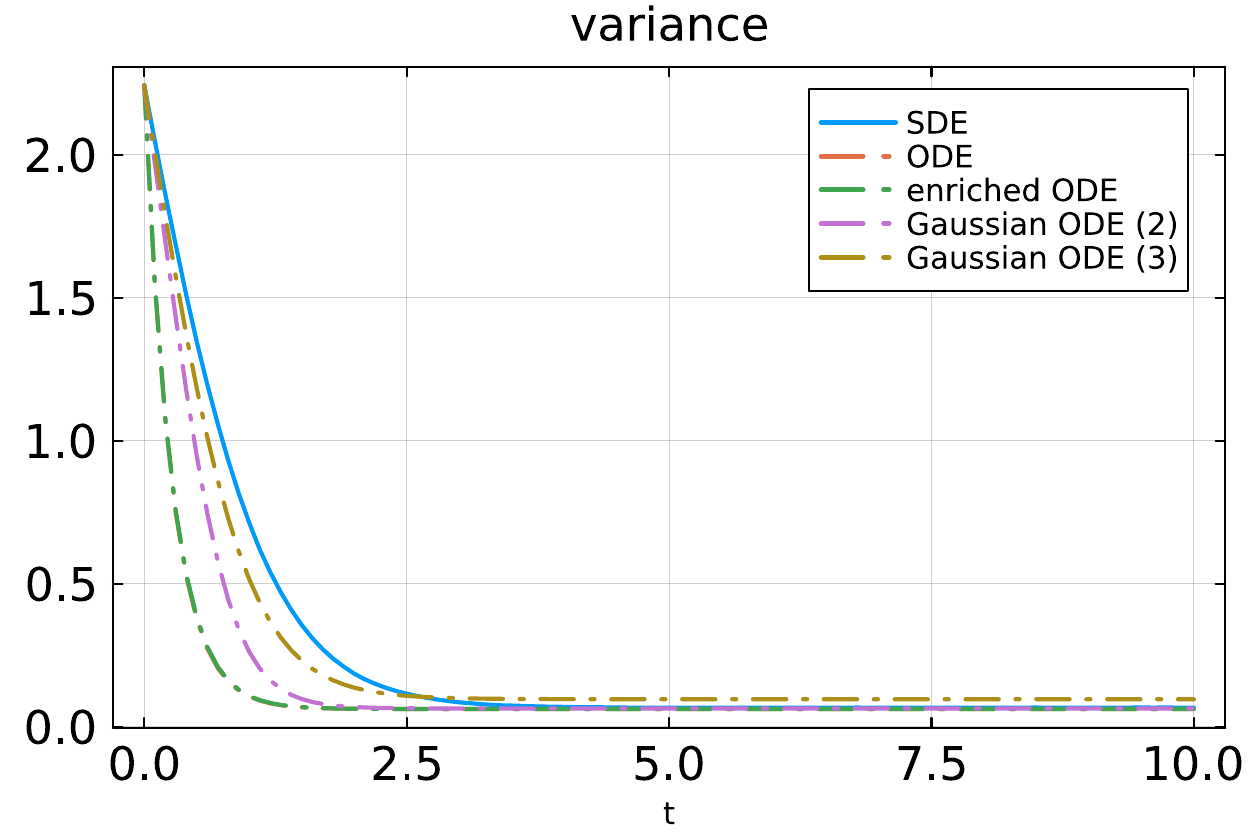}
\includegraphics[width=0.32\textwidth]{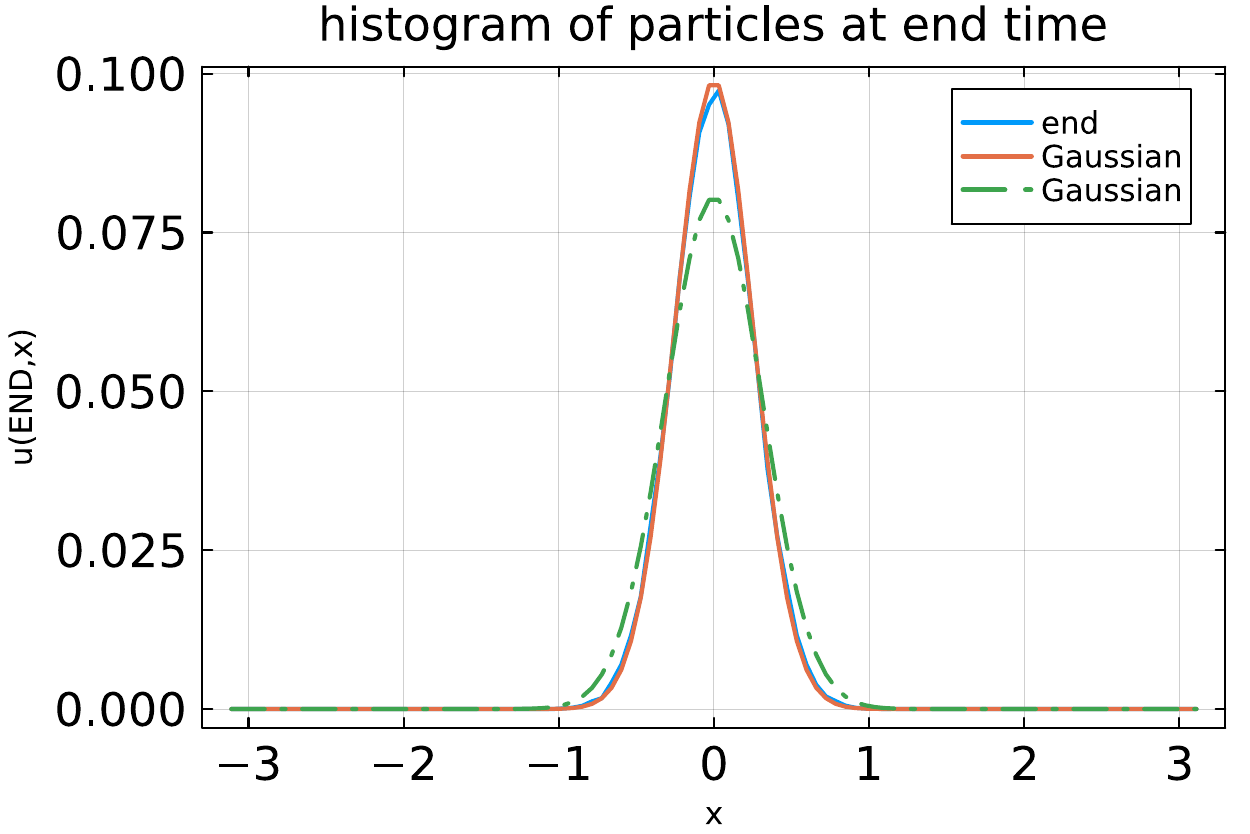}
\caption{Plane rotator: mean, variance and histogram. For the MC simulation of the McKean-Vlasov SDE, $10^5$ particles are used.}
\label{fig_plane_rotator_Taylor}
\end{figure}
In \Cref{fig_plane_rotator_WEAK} we plot the convergence of the weak error on the mean and variance, with respect to the iteration number.
The statistical error is reached in a few iterations.
For all the coarse solvers, the MC-moments Parareal method converges only in a few iterations to the level of the statistical error.
\label{fig_plane_rotator_WEAK_LINE}

\begin{figure}[h]
\centering
\includegraphics[width=0.33\textwidth]{{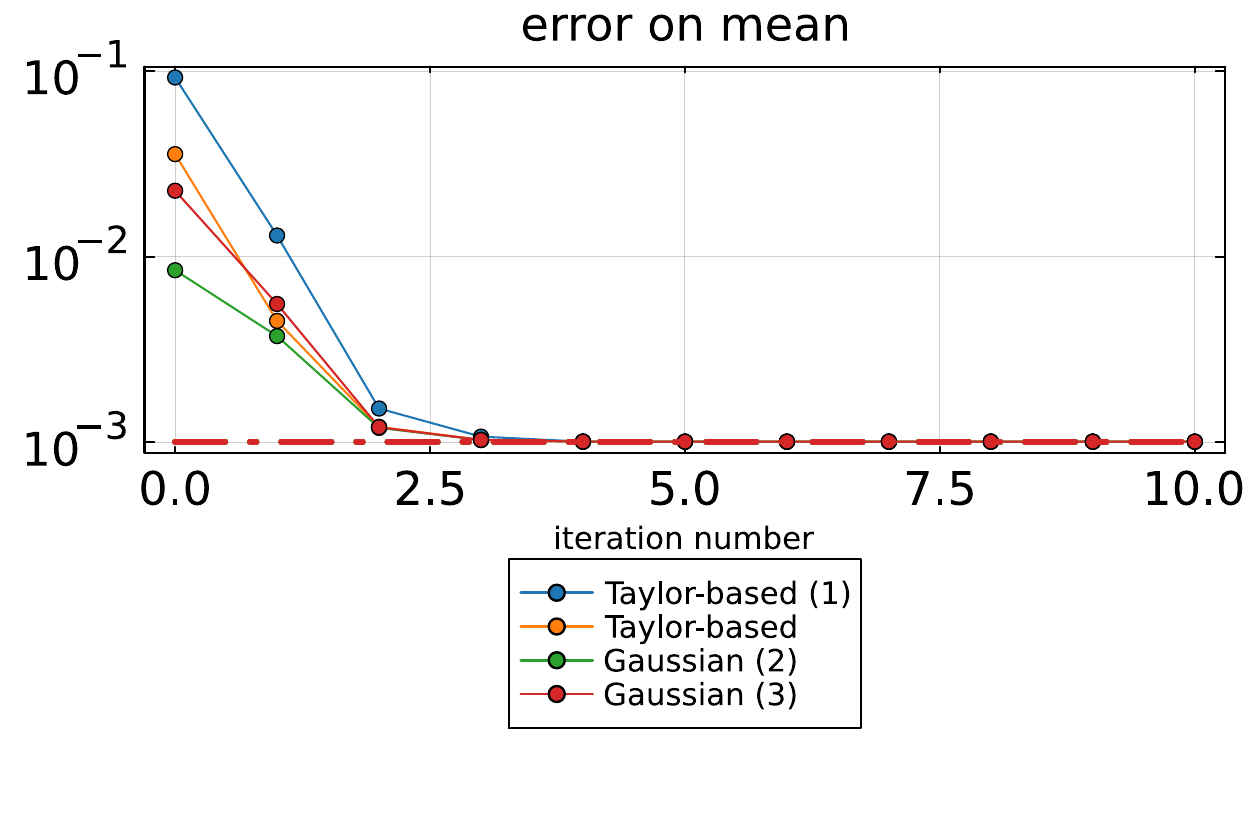}}%
\includegraphics[width=0.33\textwidth]{{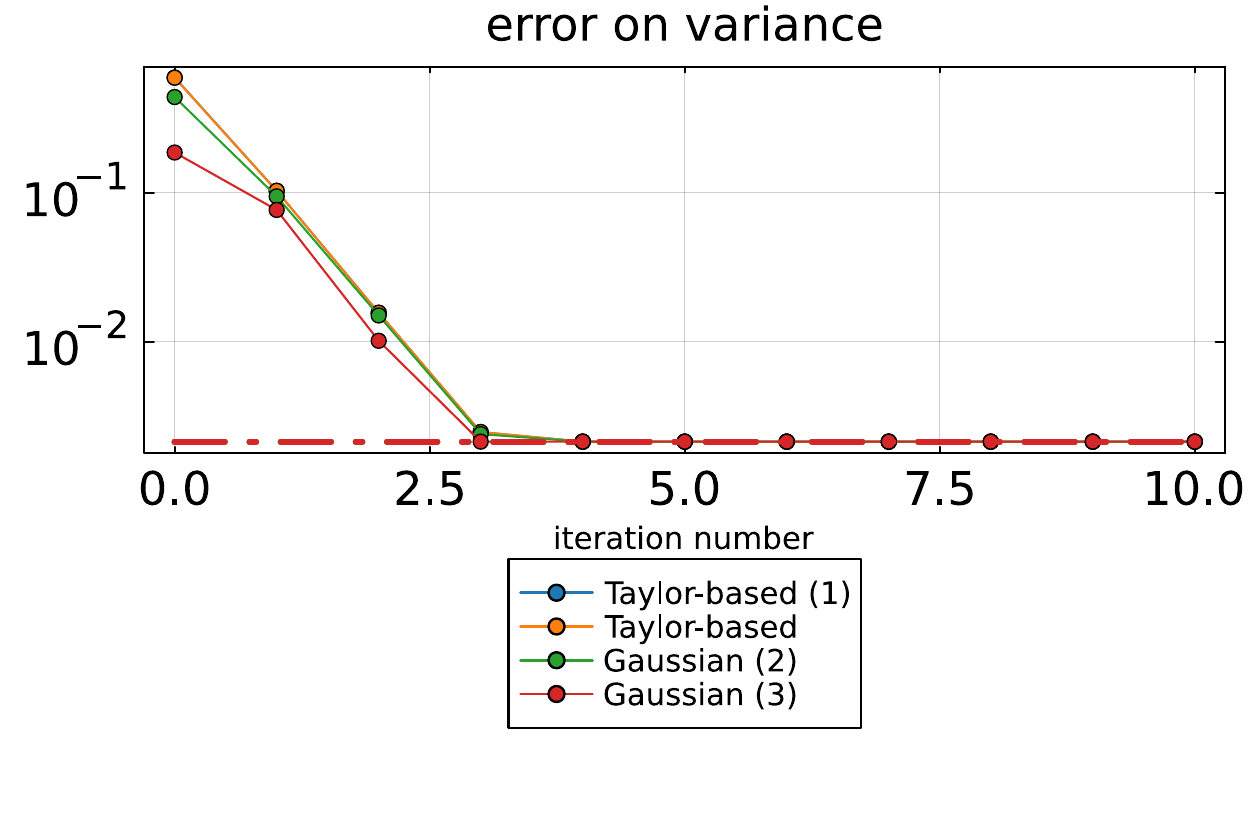}}%
\includegraphics[width=0.33\textwidth]{{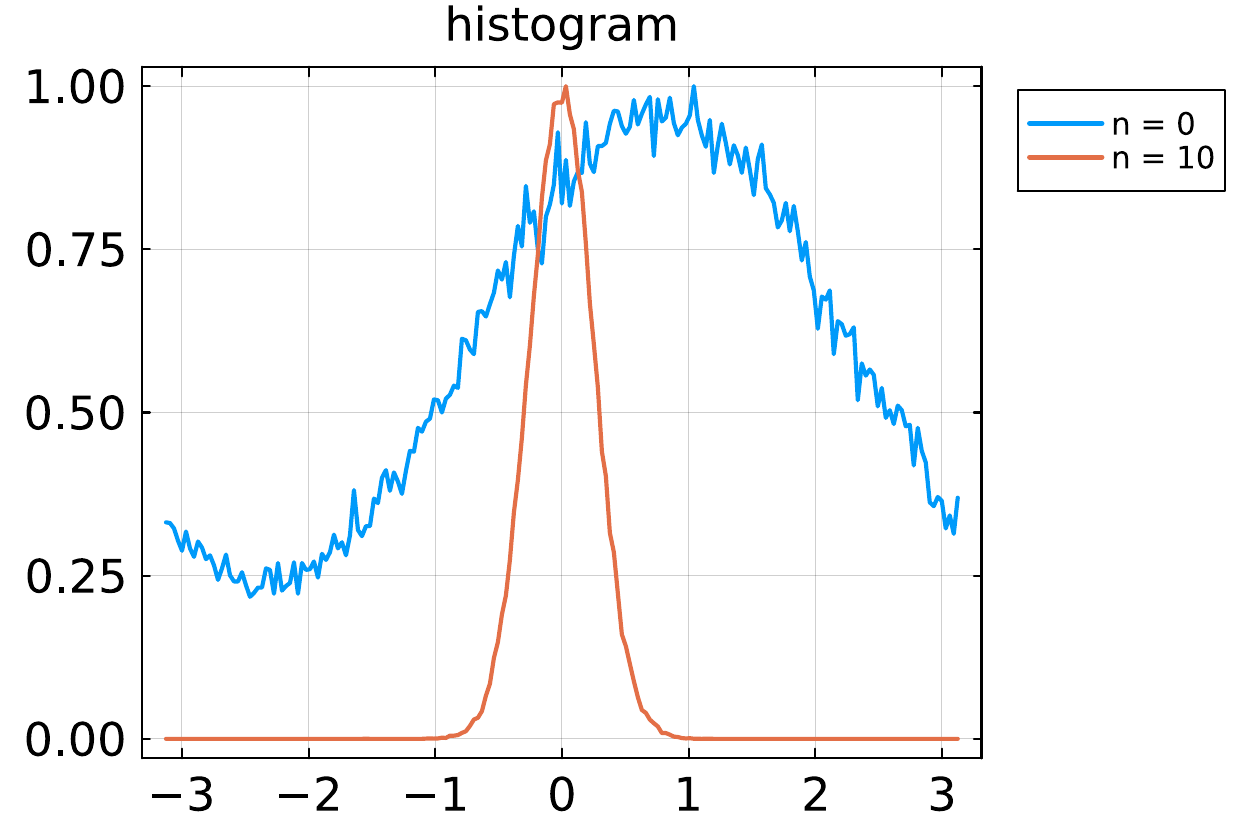}}
\includegraphics[width=0.33\textwidth]{{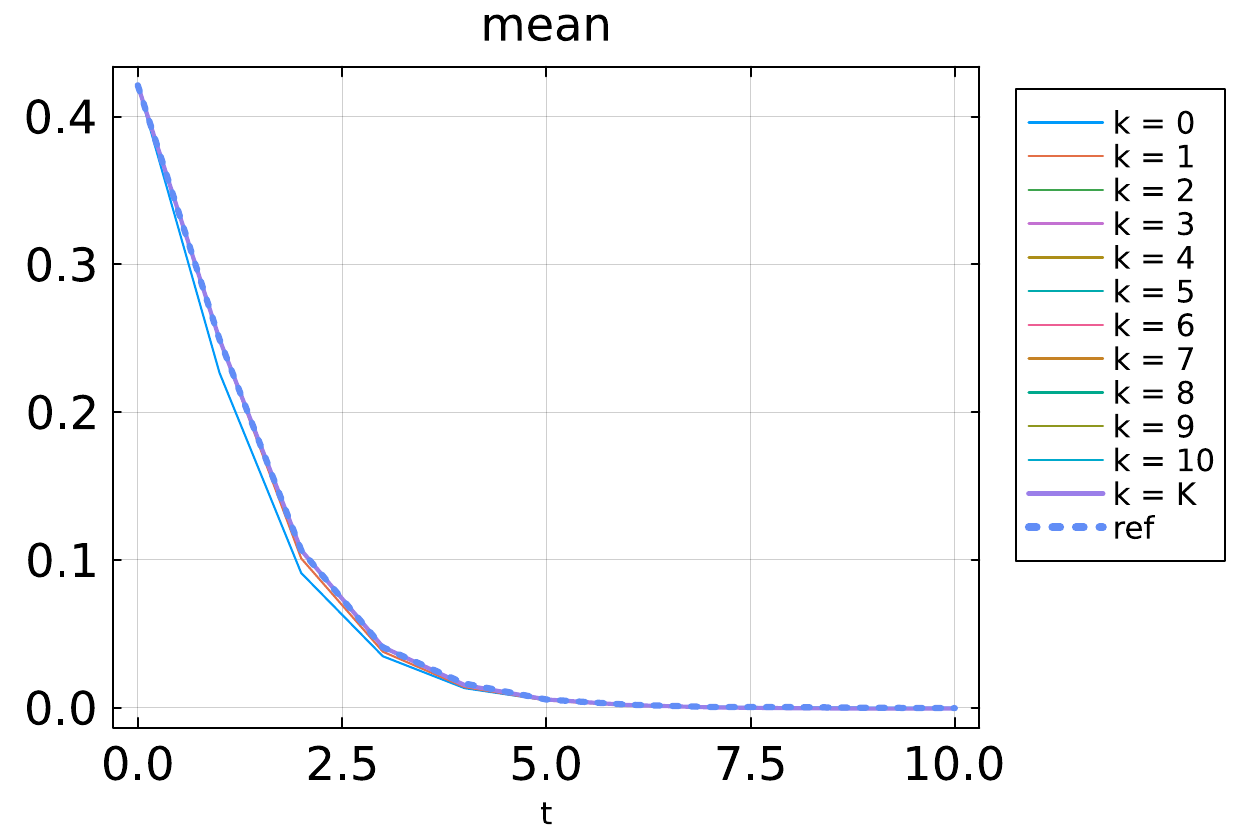}}%
\includegraphics[width=0.33\textwidth]{{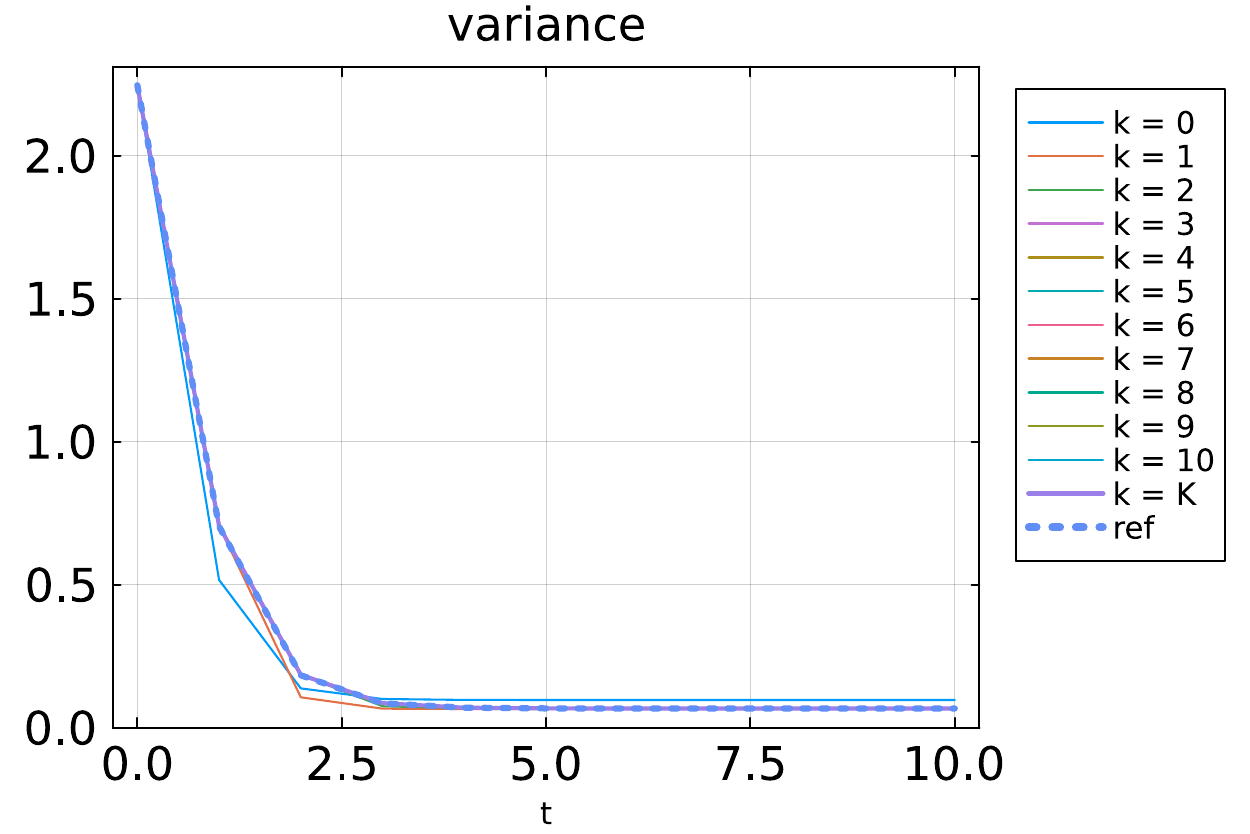}}%
\includegraphics[width=0.33\textwidth]{{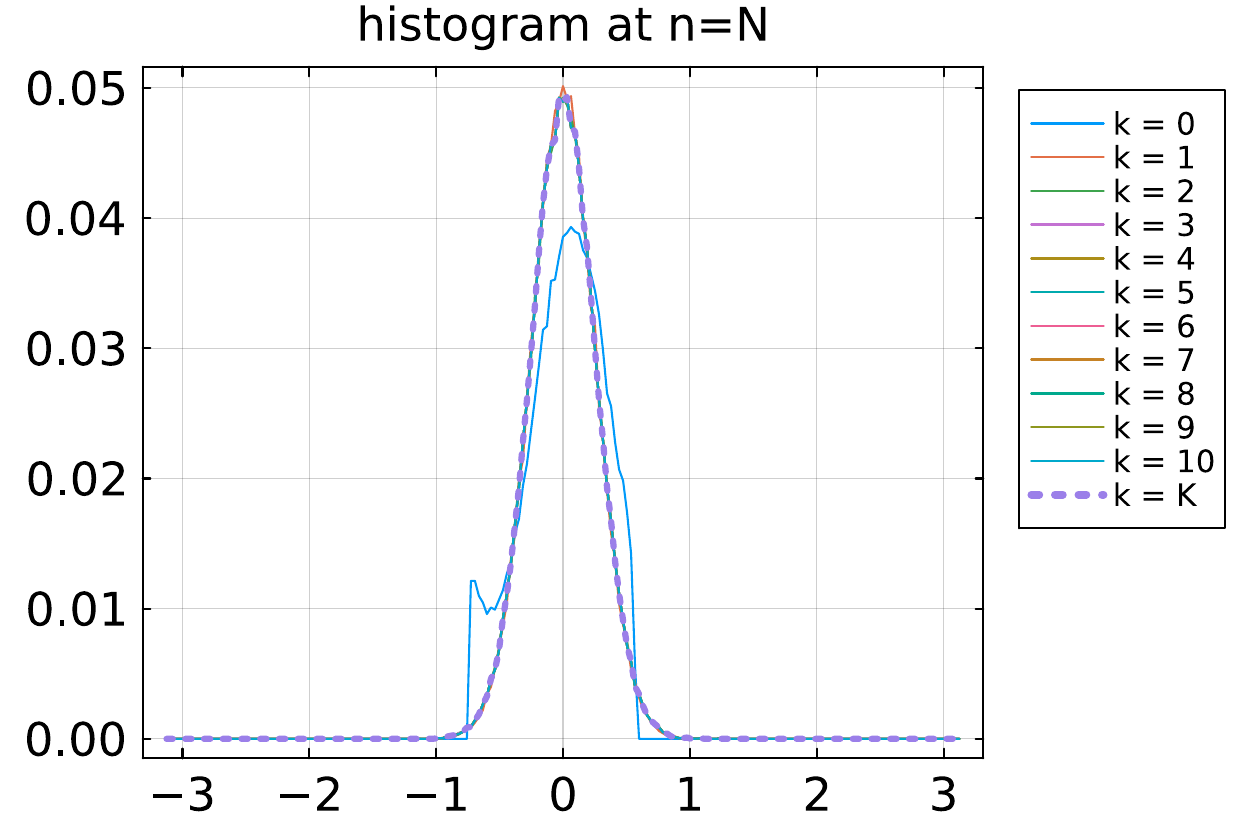}}
\caption{Plane rotator: convergence of weak error in MC-moments Parareal, with various coarse propagators. 
The mean, variance and histogram are shown on the edges of the Parareal subintervals for the Gaussian ODE with 3 sigma-points.}
\label{fig_plane_rotator_WEAK}
\end{figure}
\end{example} 

\subsection{McKean-Vlasov SDE with complex mean-field coupling}
\mylabel{Subsection_burgers}
\newcommand\exampleBurgers{viscous Burgers equation}
\begin{example}[label=exampleBurgers,name=Viscous Burgers equation]
\label{BURGERS_EQUATION_INTRO}
This example is taken from \cite{bossy_stochastic_1997}, and is also used as an example in \cite{kloeden_gauss_quadrature_2017}. 
This example is not directly applicable to our MC-moments Parareal algorithm because the drift coefficient does not contain the expected value of a function of the ensemble, but the empirical distribution function itself.
Nevertheless, we build a moment model based on a Gaussian assumption of the density function.

We consider the nonlinear hyperbolic partial differential equation
\begin{equation}
\begin{aligned}
\frac{\partial V}{\partial t} &= \frac{1}{2} \sigma^2 \frac{\partial^2 V}{\partial x^2} - V \frac{\partial V}{\partial x} 
\quad \text{with} \quad
V(0,x) &= 1 - H(x),
\end{aligned}
\label{Burgers_PDE}
\end{equation}
where $H(x)$ is a step function such that $H(x)=1 \, \forall x >0$, $H(0)=1/2$ and 0 elsewhere.
It is possible to associate with \eqref{Burgers_PDE} the McKean-Vlasov SDE \cite{bossy_stochastic_1997}
\begin{equation}
\begin{aligned}
dX^{(p)} &= \left[ \int_{\mathbb{R}} (1-H(X-y)) P(dy) \right] dt + \sigma dW 
\quad \text{with} \quad
X(0) &= 0.
\end{aligned}
\label{Burgers_SDE}
\end{equation}
Equation \eqref{Burgers_SDE} can be rewritten as
\begin{equation}
\begin{aligned}
dX^{(p)} &= (1 - \text{CDF}(X))dt + \sigma dW 
\quad \text{with} \quad
X(0) &= 0.
\end{aligned}
\label{Burgers_SDE_alternative_formulation}
\end{equation}

\paragraph{Derivation of the moment model}
Using the assumption that the particles are normally distributed (also called a Maxwellian distribution), an approximate moment model for \eqref{Burgers_SDE_alternative_formulation} reads
\begin{equation}
\begin{aligned}
\frac{dM}{dt} &= 1 - \left. \text{\text{CDF}} \right|_{(X=M)}  = \frac{1}{2}, \\
\frac{d\Sigma}{dt} &= - \frac{2}{\sqrt{2 \pi}} \sqrt{\Sigma} + \sigma^2.
\end{aligned}
\label{moment_model_Burgers}
\end{equation}
In the derivation, we used the fact that the integral of a normal distribution from $-\infty$ to its mean equals to $1/2$, and the fact that
\begin{equation}
\left. \frac{d(1-\text{CDF}(X))}{dX} \right|_{(X=m)} = -\text{PDF}(X=m) = - \frac{2}{\sqrt{2\pi \Sigma}}.
\end{equation}
 
The moment model for the mean of the particle ensemble is exact, since the solution is a traveling wave with speed $1/2$ \cite{bossy_stochastic_1997}.
In \cref{fig_Burgers_convergence} the convergence of the mean and the variance is plotted. 
The statistical error is reached in a few iterations.
\begin{figure}[H]
\centering
\includegraphics[width=0.33\textwidth]{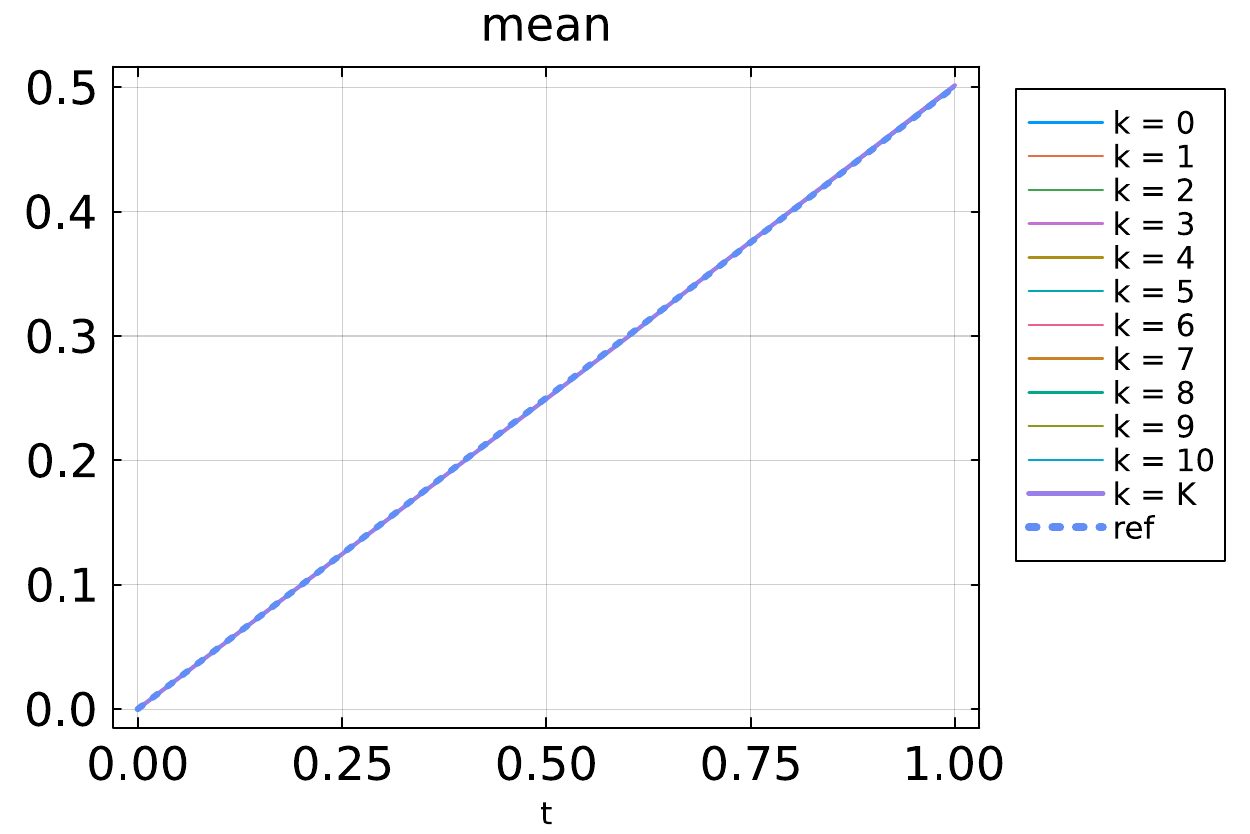}
\includegraphics[width=0.33\textwidth]{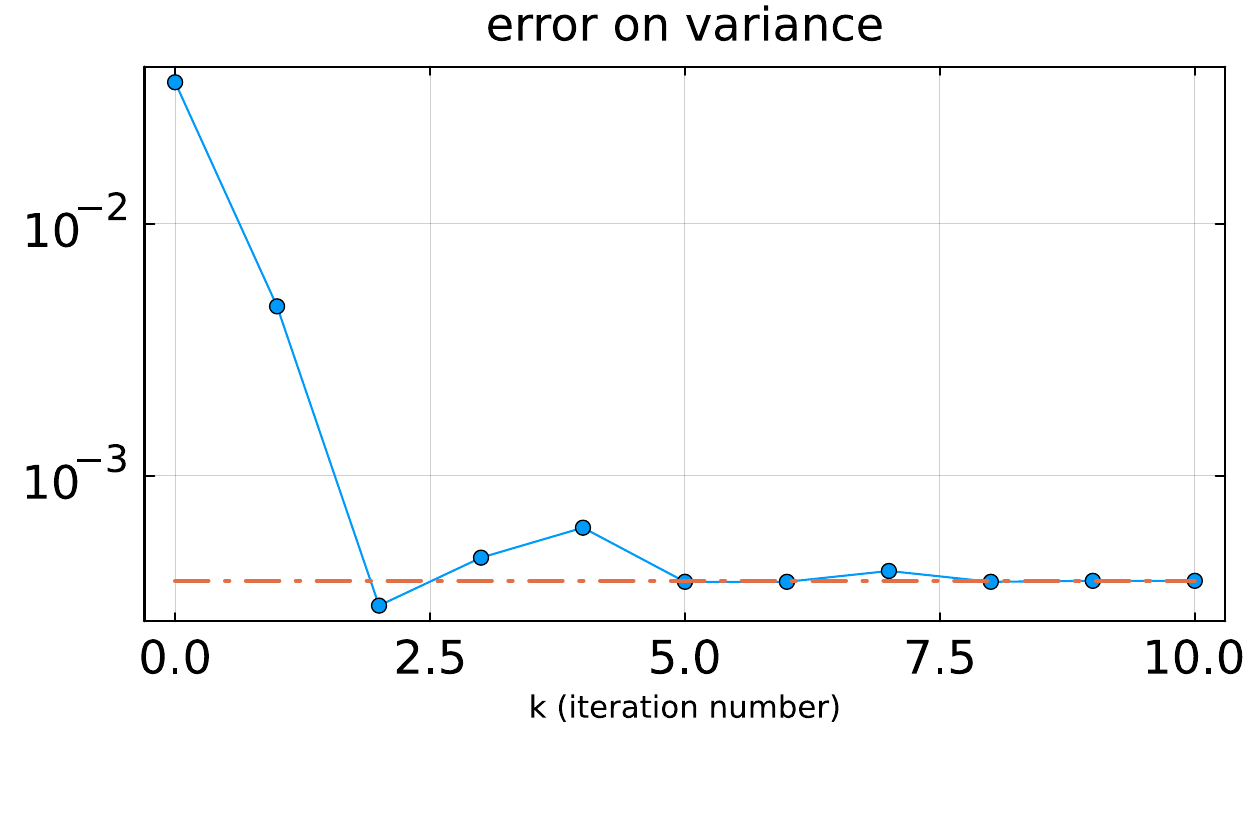}%
\includegraphics[width=0.33\textwidth]{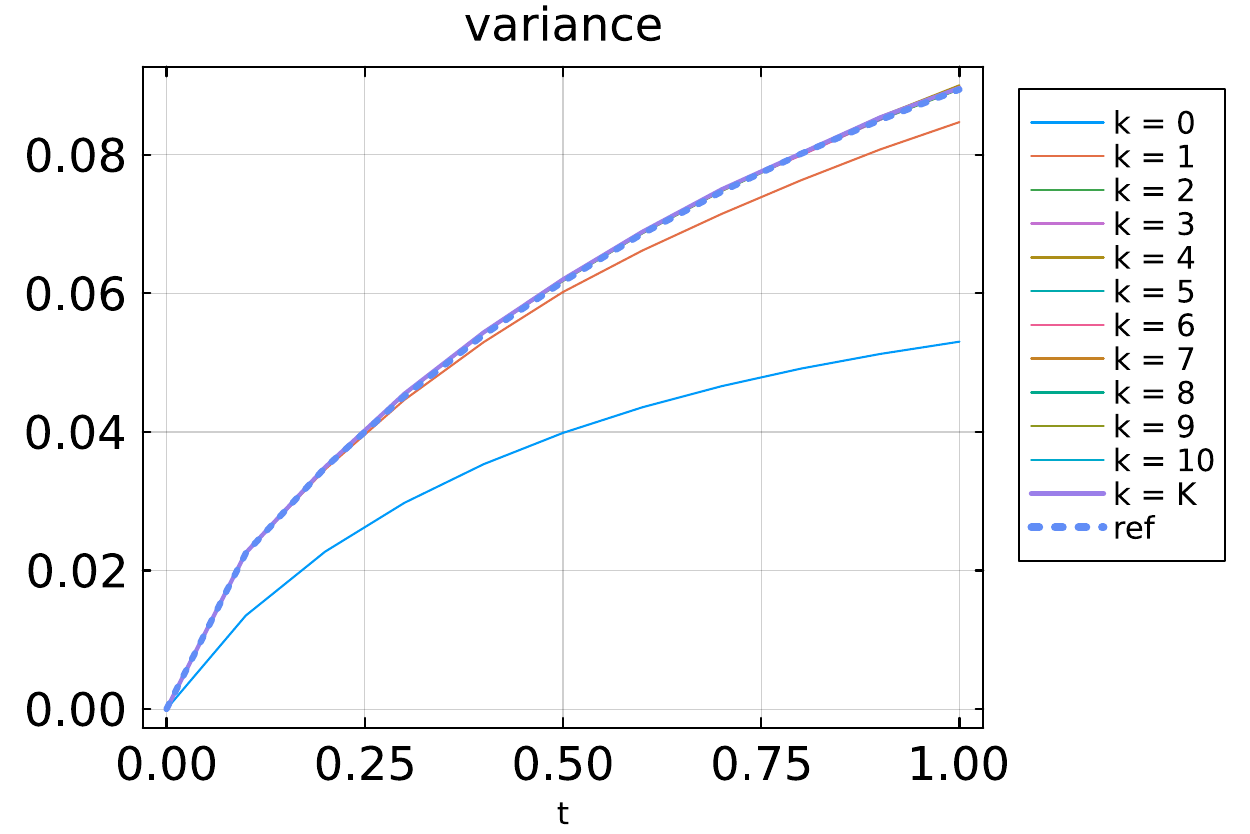}
\caption{Convergence of MC-moments Parareal for a stochastic particle discretisation of the Burgers equation. (left) time-plot of $M^k_n$ for different iteration numbers $k$, (middle) Convergence of the weak error on the mean, (right) time-plot of $\Sigma^k_n$ for different iteration numbers $k$}
\label{fig_Burgers_convergence}
\end{figure}

\begin{remark}[More general initial conditions]
\label{remark_about_weighted_particle_ensemble}
If in \eqref{Burgers_PDE} another initial condition was chosen, containing some negative particles (with $X^{(p)}<0$), the particle system \eqref{general_equation} does not suffice.
Instead a weighted particle ensemble should be used to describe the system's behaviour in a stochastic way, see \cite[Section 2.3]{bossy_stochastic_2005}. 
In the context of the MC-moments Parareal algorithm, such a weighted particle simulation, however, is left for future work.
We choose the parameters $[0,T] = [0, 20]$, $\sigma = \sqrt{0.2}$ and $\Delta t = 1/10$.
\end{remark}
\end{example}

\subsection{Bimodal McKean-Vlasov SDEs}
\mylabel{numerical_simulations_bimodal_SDEs}
In this subsection, we consider a McKean-Vlasov SDEs with multimodal particle distribution and constant diffusion coefficient: 
\begin{equation}
\begin{aligned}
dX^{(p)} &= -\left( 4\alpha \left({X^{(p)}}\right)^3 - 2\gamma X^{(p)} - \beta \mathbb E [X^{(p)}] + J \sqrt{\frac{a}{2b}} \right) dt + \sigma dW^{(p)},
\quad 
\bar X(0) \sim p_{\mathcal N} (m_0, \sigma_0).
\end{aligned}
\label{SDE_bimodal_first_time}
\end{equation}
Here, $\alpha$, $\gamma$, $\beta$ and $J$ are model parameters.
Without mean-field term $\beta \mathbb E[X^{(p)}]$, the drift term of SDE \eqref{SDE_bimodal_first_time} can be derived from the potential
\begin{equation}
\begin{aligned}
V(X) &= aX^4 - bX^2 + J\sqrt{\frac{a}{2b}}X.
\end{aligned}
\end{equation}
If $J=0$, the potential $V$ is symmetric around $X=0$. 
The parameter $\beta$ is the strength of the mean-field coupling and $m_0$ and $\sigma_0$ are the mean value and variance of the initial particle distribution.
For all simulations we take $J=0$, since we can bring in asymmetry by a proper choice of the end time and the initial condition.

\begin{figure}[h]
\begin{subfigure}[t]{.33\linewidth}
\includegraphics[width=1\textwidth]{{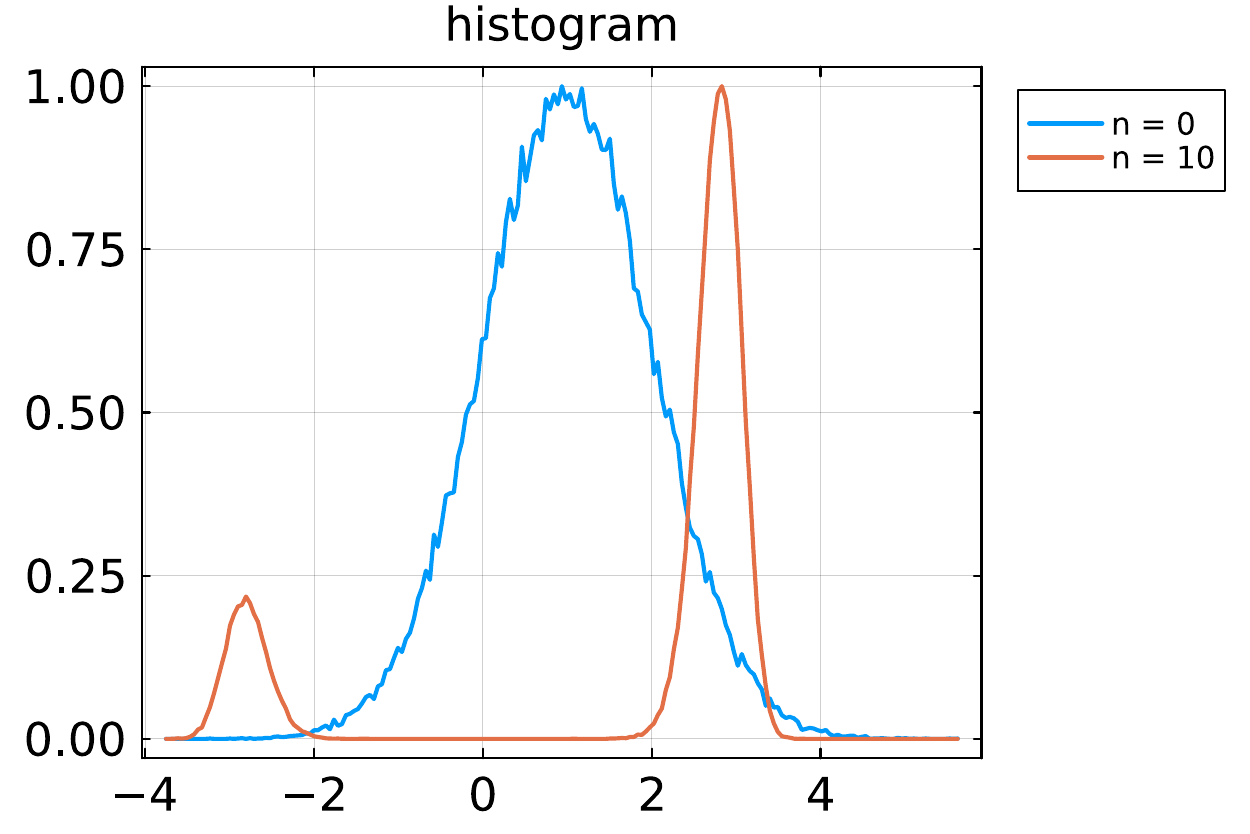}}%
	\caption{Setup 1}
  \end{subfigure}
  \begin{subfigure}[t]{.33\linewidth}
\includegraphics[width=1\textwidth]{{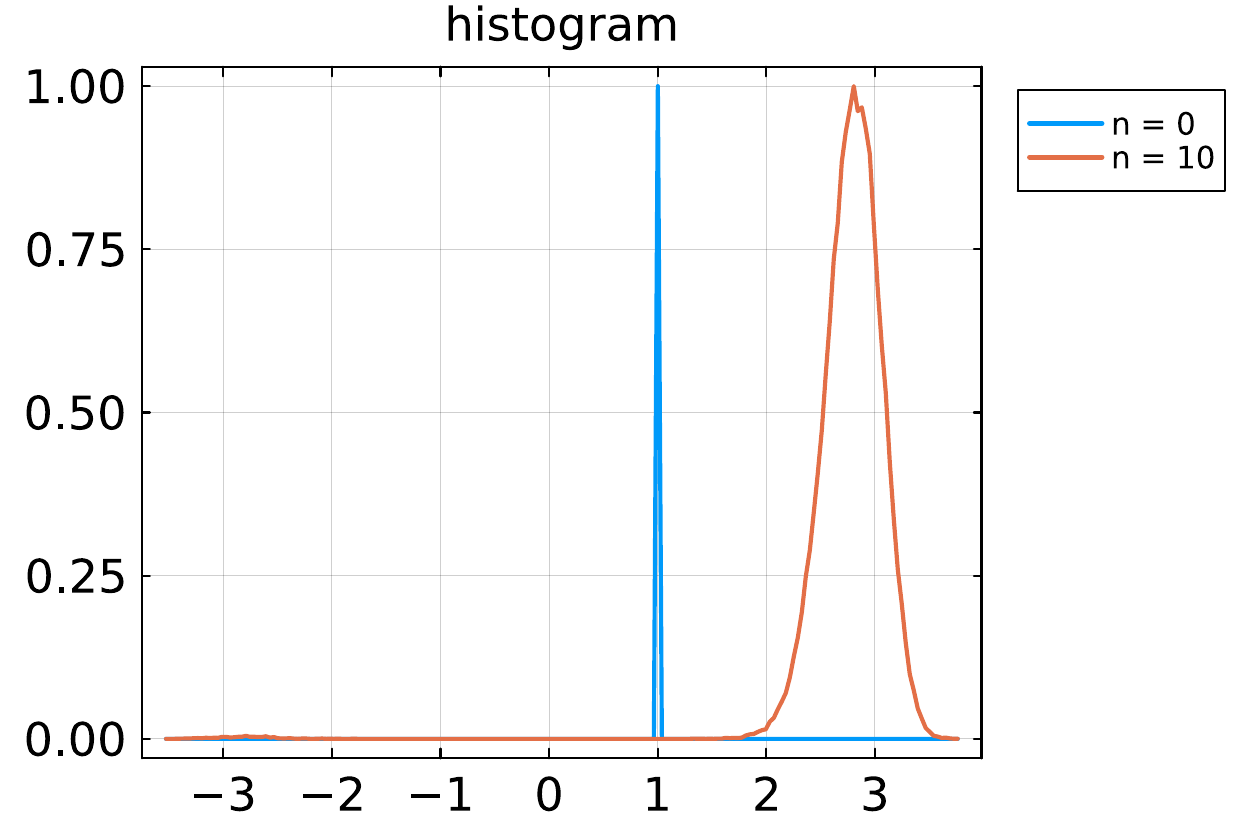}}%
	\caption{Setup 2}
  \end{subfigure}
  \begin{subfigure}[t]{.33\linewidth}
\includegraphics[width=1\textwidth]{{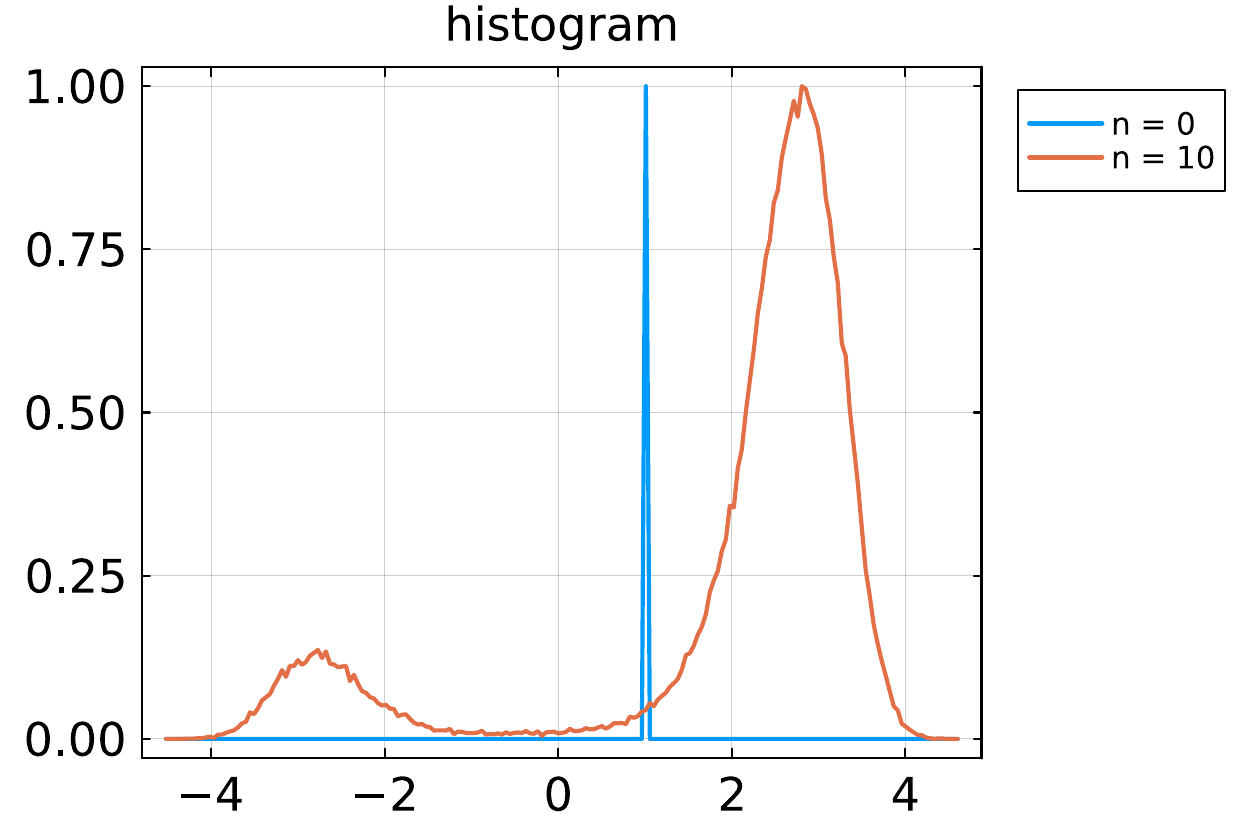}}%
	\caption{Setup 3}
  \end{subfigure}
    \begin{minipage}{.25\linewidth}
    \vspace*{-2cm}
  \caption{Overview of test systems for the bimodal SDE. 
  The histograms are normalised such that their maximal value equals one.}
  \label{fig_overview_all_tests}
    \end{minipage}
    \hspace*{0.06\linewidth}
    \begin{subfigure}[t]{.33\linewidth}
\includegraphics[width=1\textwidth]{{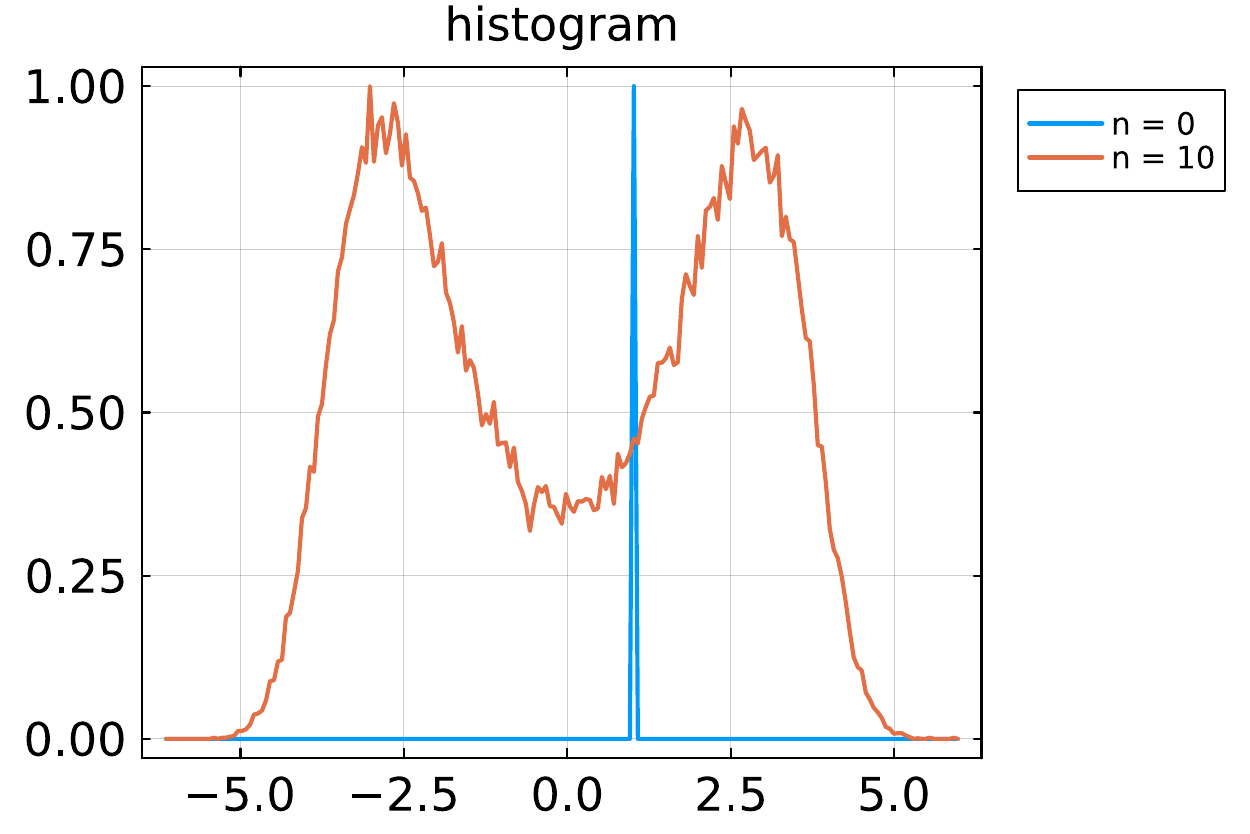}}%
	\caption{Setup 4}
  \end{subfigure}
  \begin{subfigure}[t]{.33\linewidth}
\includegraphics[width=1\textwidth]{{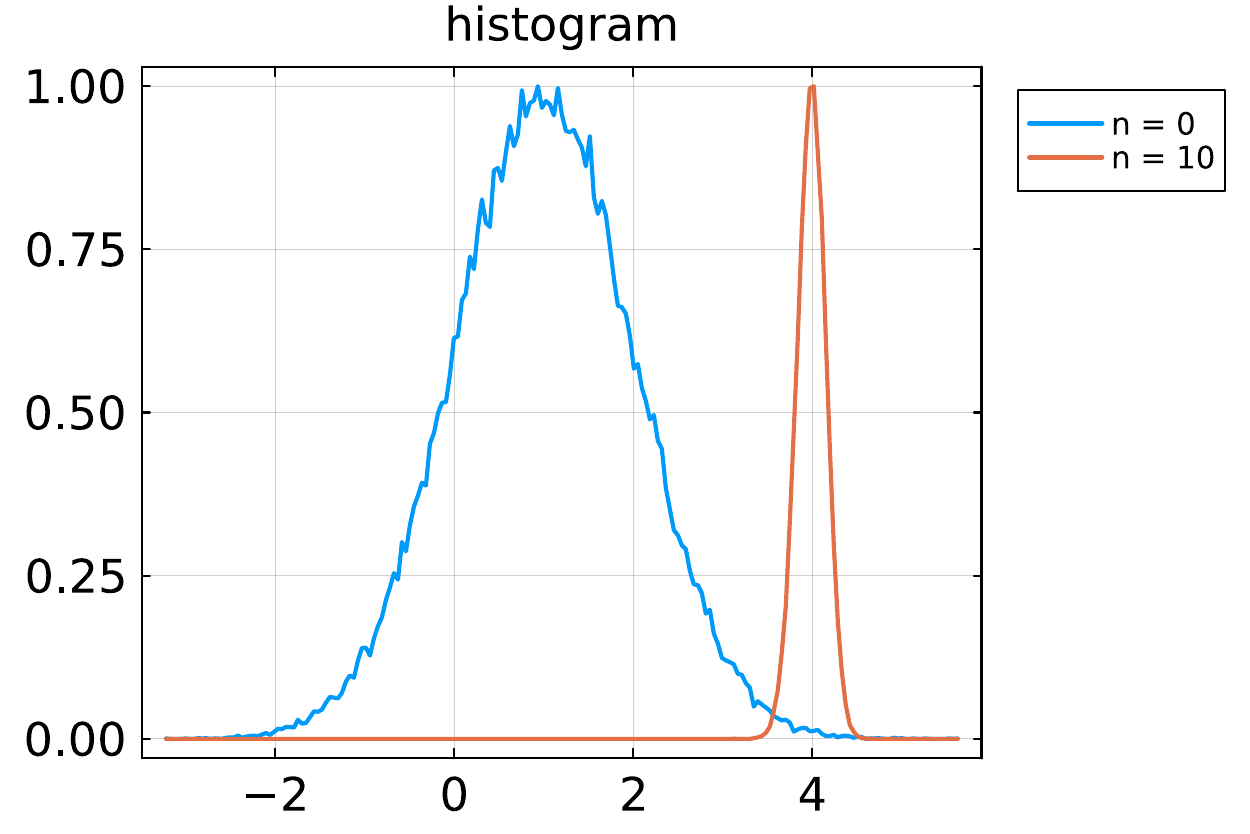}}%
	\caption{Setup 5}
  \end{subfigure}
\end{figure}

\paragraph{Parameters of the McKean-Vlasov SDE}
In \cref{table_overview_settings_bimodal_SDE}, we give an overview of the tests setups (parameters) used in our numerical simulations.
In \cref{fig_overview_all_tests} we provide an overview of the histograms in these setups at time $t=0$ and $t=T$.
In setup 2,  at $t=T$, there are a very small number of particles in the left-sided well of the potential.
In setup 5, the mean-field coupling makes the McKean-Vlasov SDE as a unimodal one.

\begin{table}[H]
\centering
\begin{tabular}{c||c|c|c|c||c|c||c|c|c}
name	
	& $a$
	& $\gamma$
	& $\beta$ 
	& $\sigma$
	& $m_0$
	& $\sigma_0$
	& $T$
	& property
	& example
\\ \hline
Setup 1 
	& $1/4$
	& $1/2$	
	& 0
	& 0.5
	& 1
	& 1
	& 20
	& 
	& \ref{Bimodal_SDE} \\
Setup 2 
	& $1/4$
	& $1/2$	
	& 0
	& 0.5
	& 1
	& 0
	& 10
	& initialisation in one well 
	& \ref{Bimodal_SDE_unimodal_start} \\
Setup 3 
	& $1/4$
	& $1/2$	
	& 0
	& 1
	& 1
	& 0
	& 10
	& idem, with more diffusion 
	& \ref{Bimodal_SDE_setup_3} \\	
Setup 4 
	& $1/4$
	& $1/2$	
	& 0
	& 2 
	& 1
	& 0
	& 50
	& weakly bimodal 
	& \ref{Bimodal_SDE_6} \\
Setup 5
	& $1/4$
	& $1/2$	
	& 1
	& 0.5 
	& 1
	& 1
	& 10
	& non-zero mean-field coupling 
	& \ref{Bimodal_SDE_MEAN_FIELD} \\	
\end{tabular}
\caption{Bimodal SDE: overview of test setups.}
\label{table_overview_settings_bimodal_SDE}
\end{table}

\begin{figure}[h!]
\begin{subfigure}[t]{.33\linewidth}
\includegraphics[width=1\textwidth]{{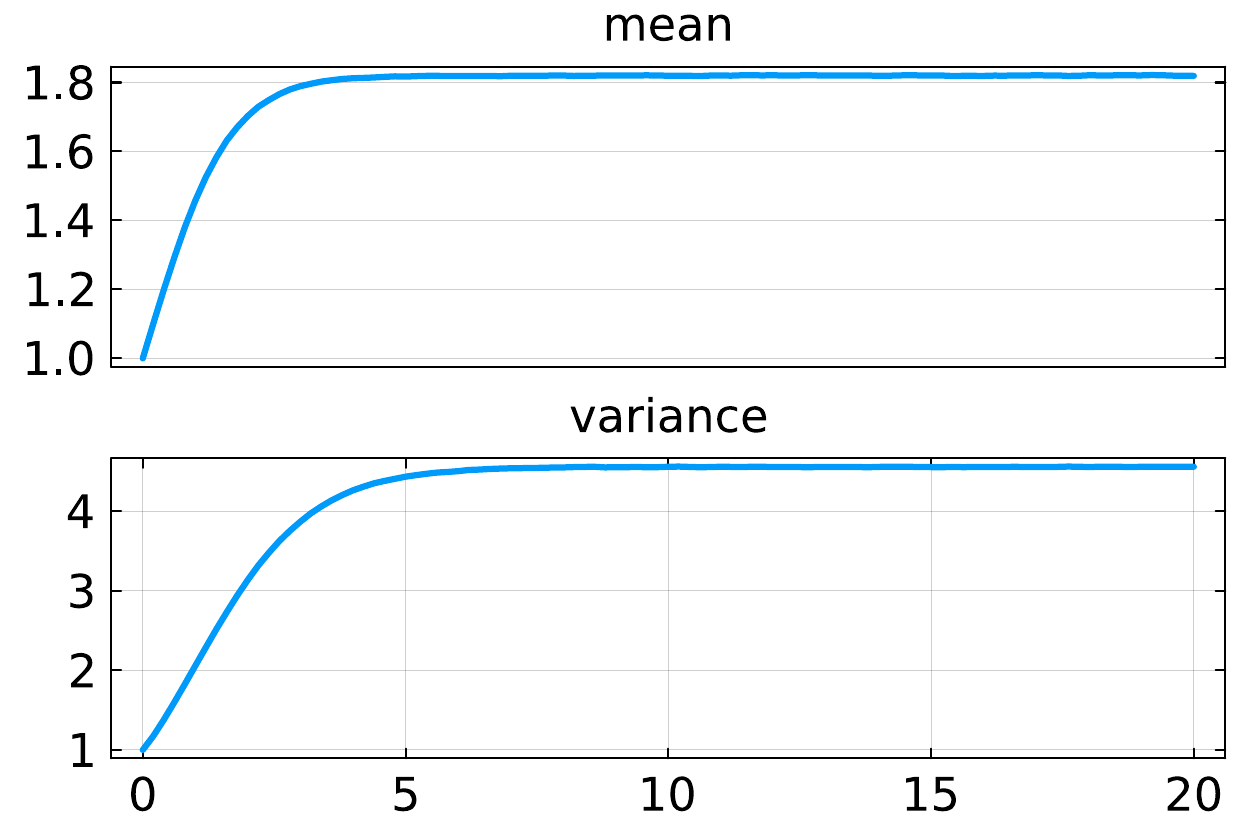}}%
	\caption{Setup 1}
  \end{subfigure}
  \begin{subfigure}[t]{.33\linewidth}
\includegraphics[width=1\textwidth]{{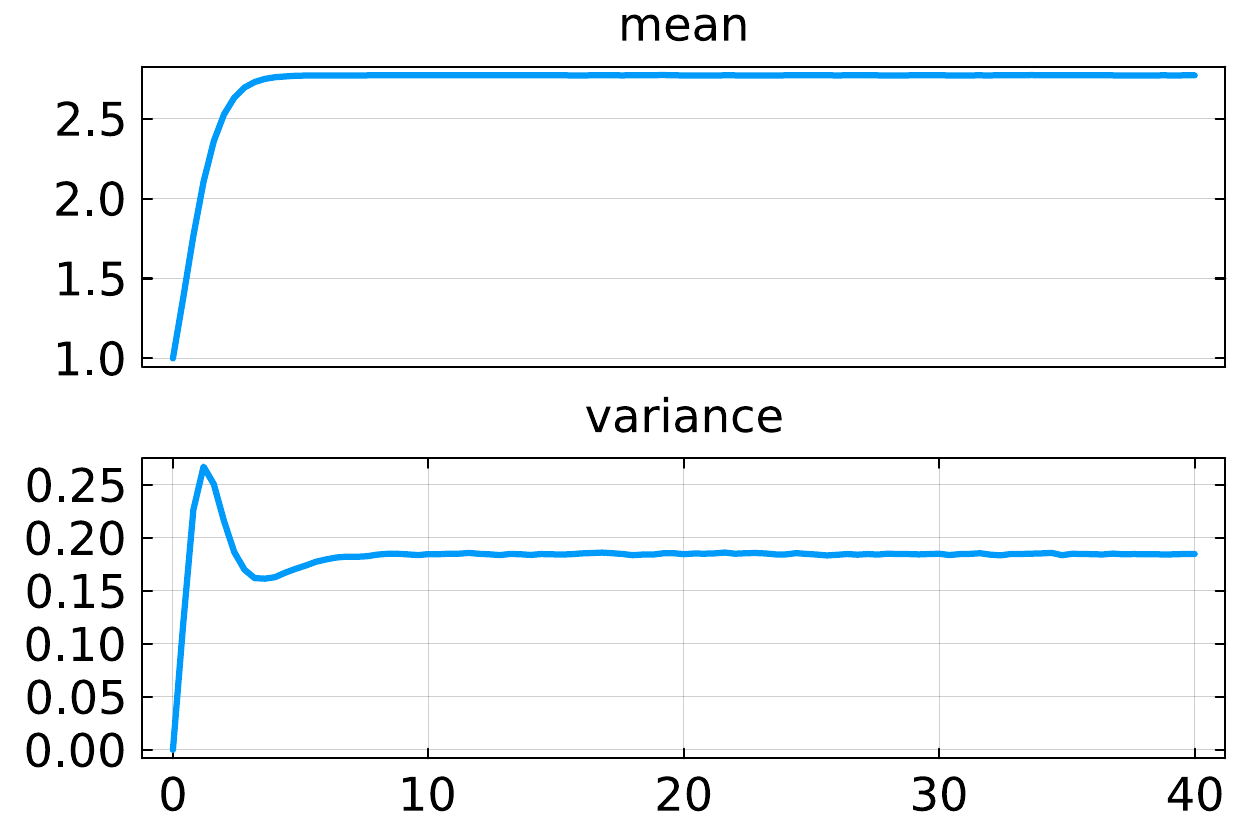}}%
	\caption{Setup 2}
  \end{subfigure}
  \begin{subfigure}[t]{.33\linewidth}
\includegraphics[width=1\textwidth]{{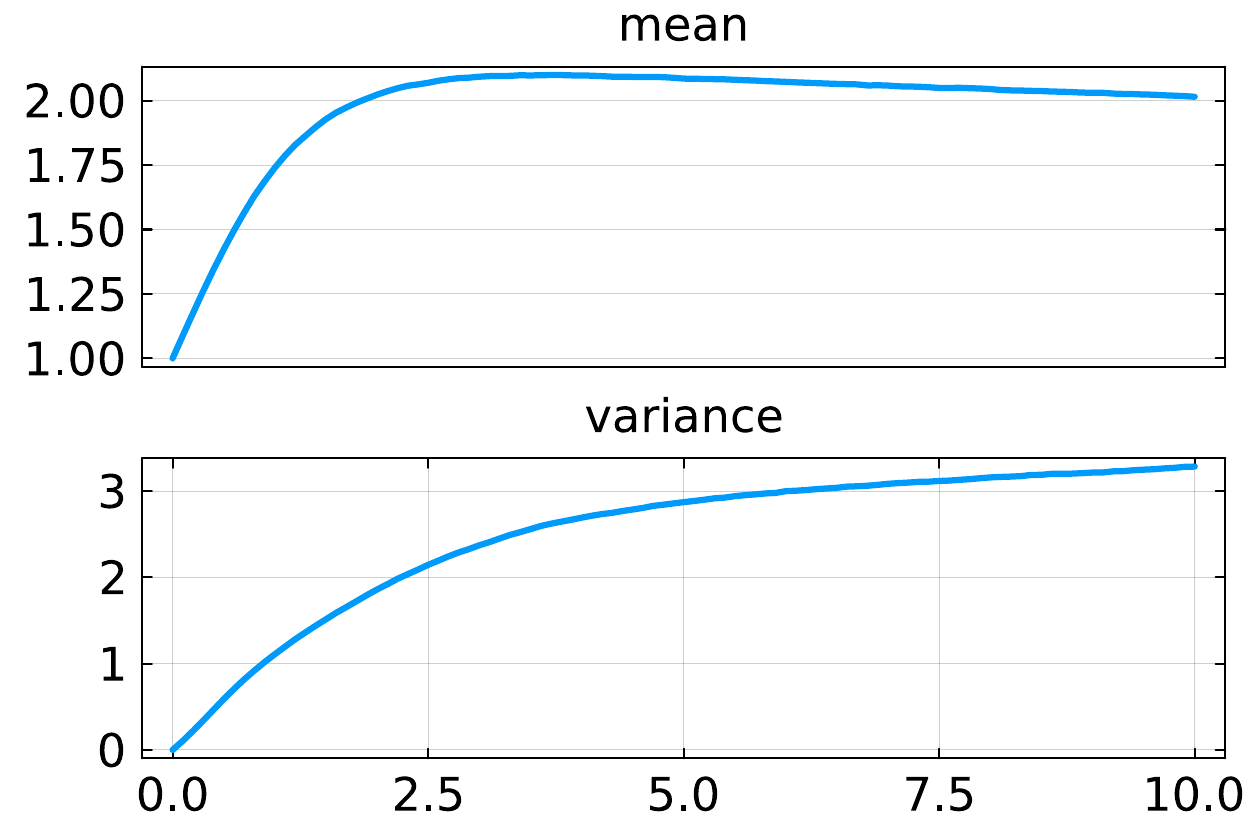}}%
	\caption{Setup 3}
  \end{subfigure}
    \begin{minipage}{.25\linewidth}
    \vspace*{-2cm}
  \caption{Overview of bimodal test setups: 
mean and variance in function of time.}
\label{fig_overview_evolution_mean_var}
    \end{minipage}
    \hspace*{0.06\linewidth}
    \begin{subfigure}[t]{.33\linewidth}
\includegraphics[width=1\textwidth]{{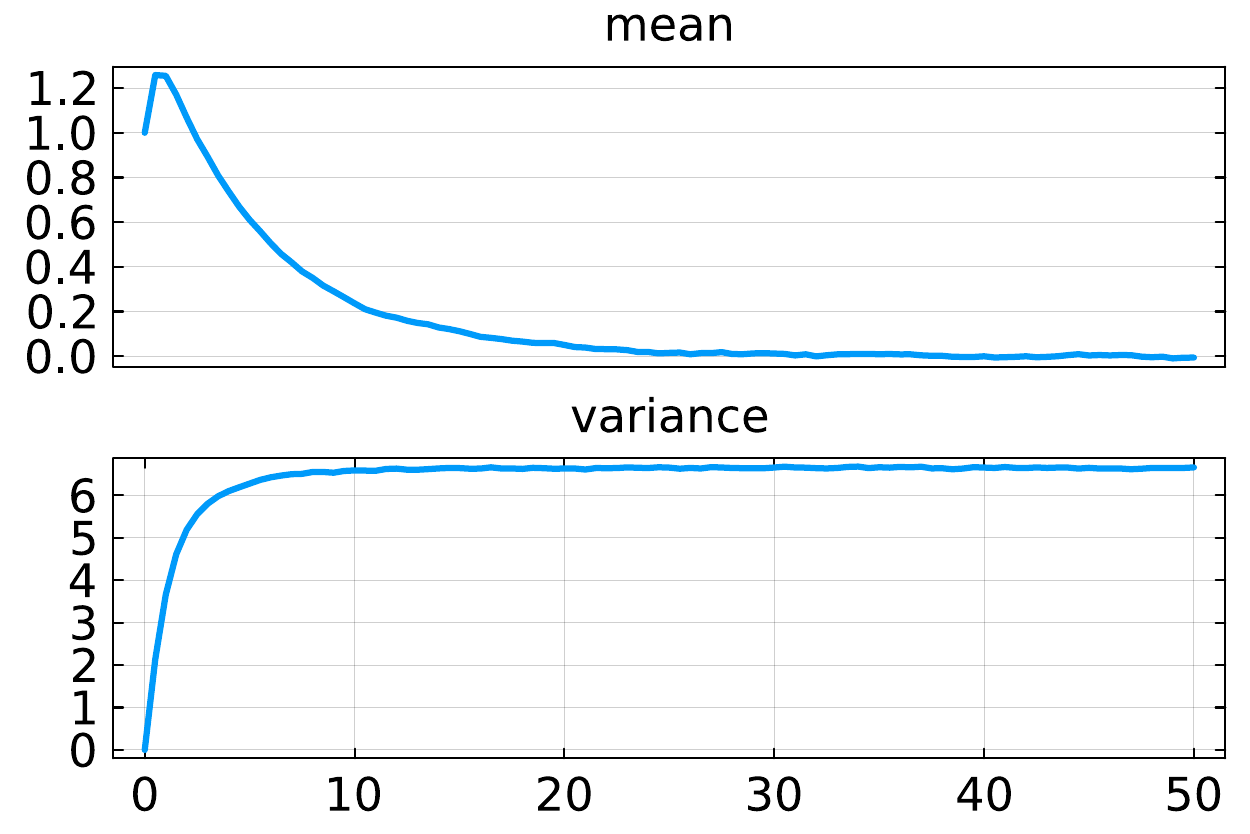}}%
	\caption{Setup 4}
  \end{subfigure}
  \begin{subfigure}[t]{.33\linewidth}
\includegraphics[width=1\textwidth]{{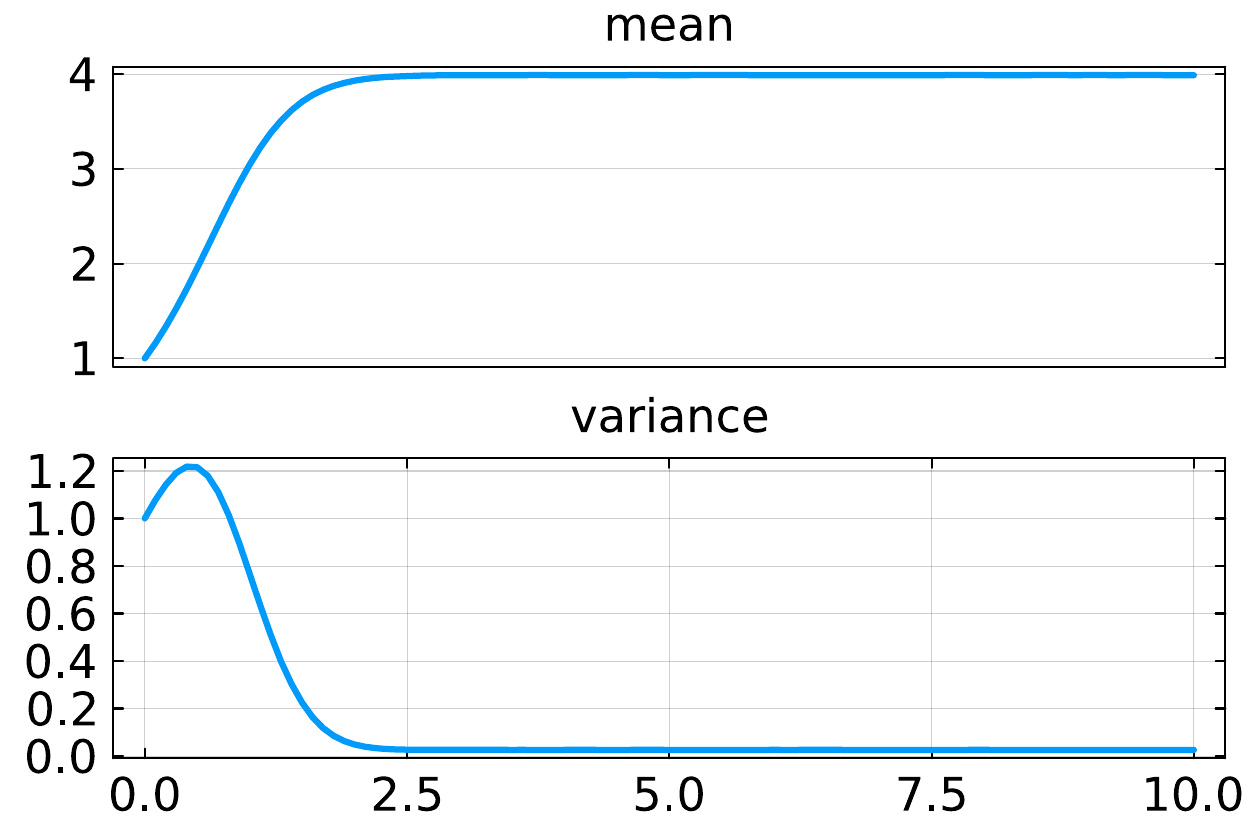}}%
	\caption{Setup 5}
  \end{subfigure}
\end{figure}

\paragraph{Variants of MC-moments Parareal}
We use these four variants of MC-moments Parareal:
\begin{enumerate}
\item Variant 1: MC-moments Parareal with a unimodal coarse propagator from \cref{definition_operators_unimodal}. For the coarse solver, the Gaussian moment model is used with three sigma-points.

\item Variant 2: MC-moments with for bimodal SDEs from \cref{MC_moments_Parareal_for_multimodal}.

\item Variant 4: Learning-based MC-moments Parareal from \cref{definition_MC_moments_Parareal_with_learning}.

\item Variant 3: MC-moments with exact particle fractions, using $M^k_{i,n}$ and $\Sigma^k_{i,n}$ but the exact particle fractions $\mathcal P_{\mathcal D_i,n}$ instead of the Parareal approximations $\mathcal P_{\mathcal D_i,n}^k$.
In practice, this is not implementable, it rather serves as an ideal reference algorithm for comparison purposes.
\end{enumerate}

\begin{figure}[H]
\begin{subfigure}[t]{.33\linewidth}
\includegraphics[width=1\textwidth]{{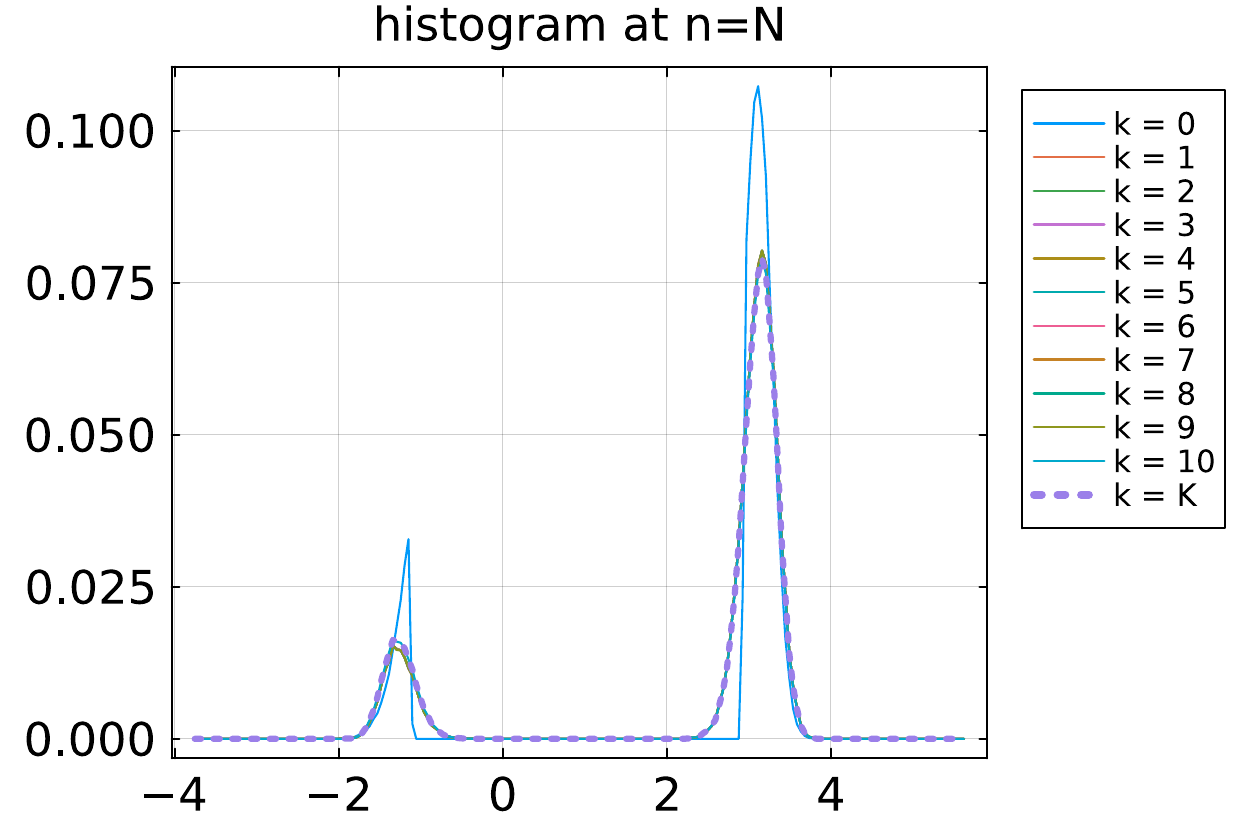}}%
	\caption{Setup 1}
  \end{subfigure}
  \begin{subfigure}[t]{.33\linewidth}
\includegraphics[width=1\textwidth]{{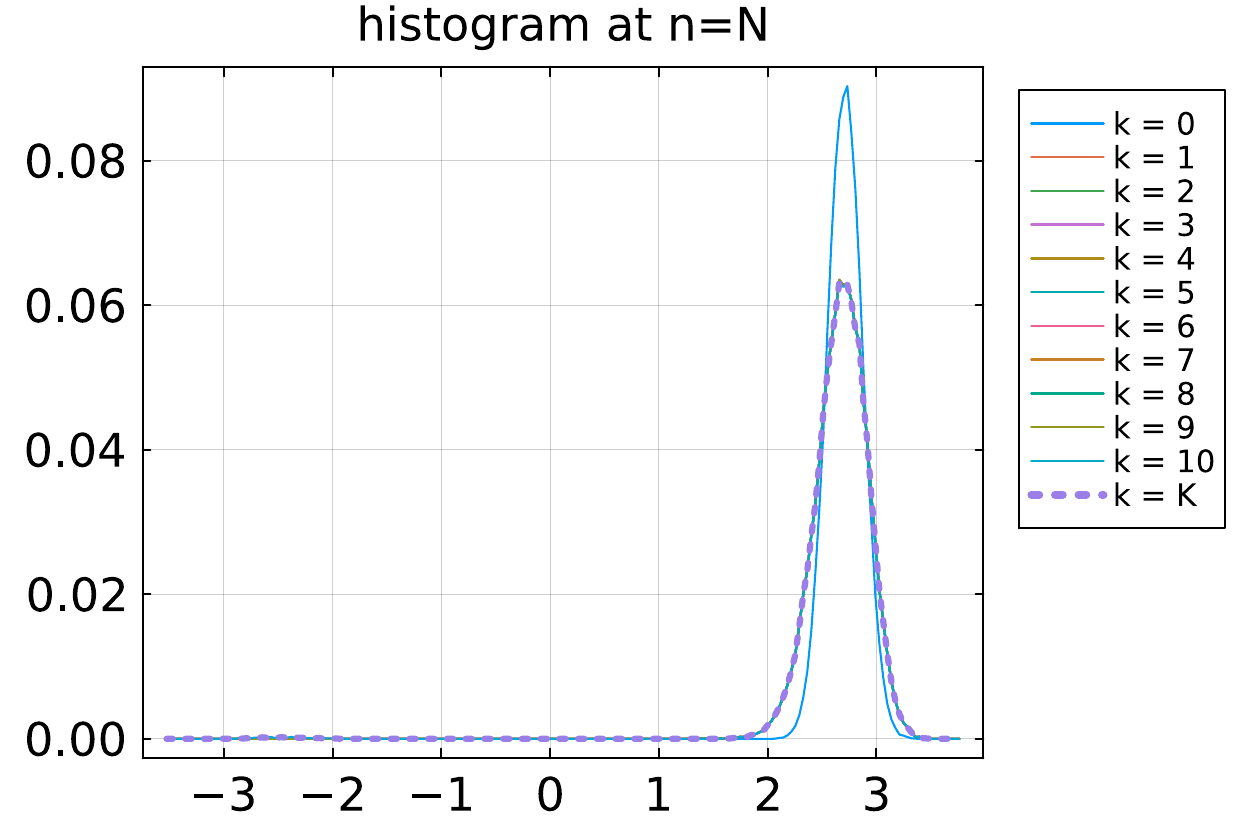}}%
	\caption{Setup 2}
  \end{subfigure}
  \begin{subfigure}[t]{.33\linewidth}
\includegraphics[width=1\textwidth]{{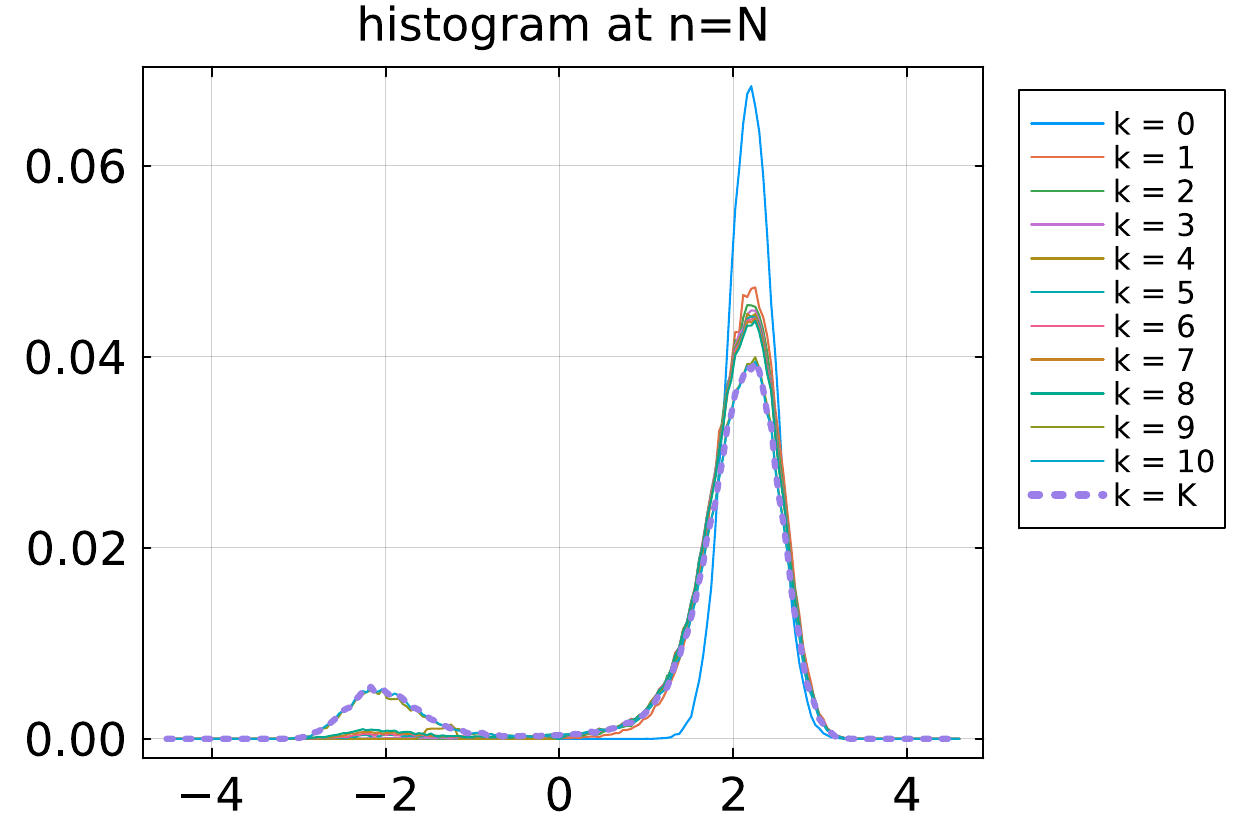}}%
	\caption{Setup 3}
  \end{subfigure}
    \begin{minipage}{.25\linewidth}
    \vspace*{-2cm}
  \caption{Convergence of the histogram at $t=T$ for the different bimodal test setups.}
  \label{fig_overview_convergence_histograms}
    \end{minipage}
    \hspace*{0.06\linewidth}
    \begin{subfigure}[t]{.33\linewidth}
\includegraphics[width=1\textwidth]{{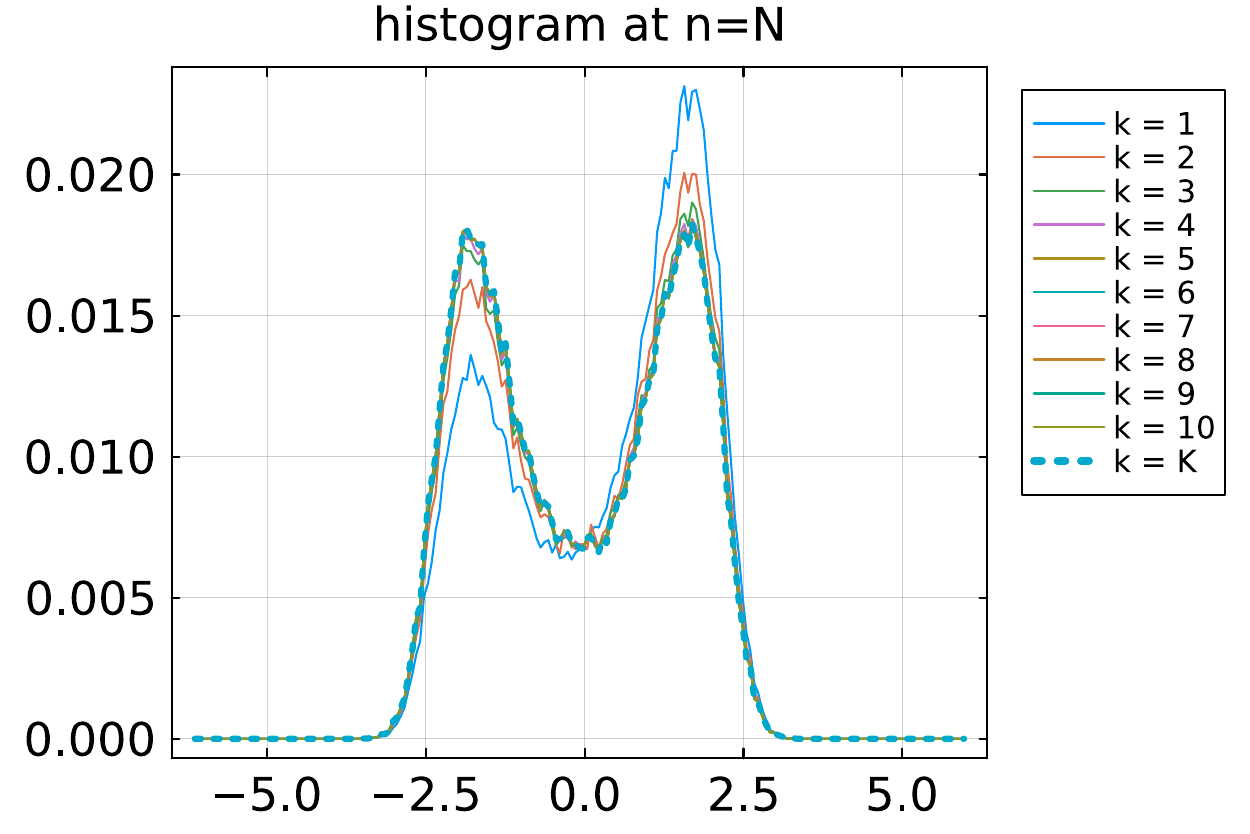}}%
	\caption{Setup 4}
  \end{subfigure}
  \begin{subfigure}[t]{.33\linewidth}
\includegraphics[width=1\textwidth]{{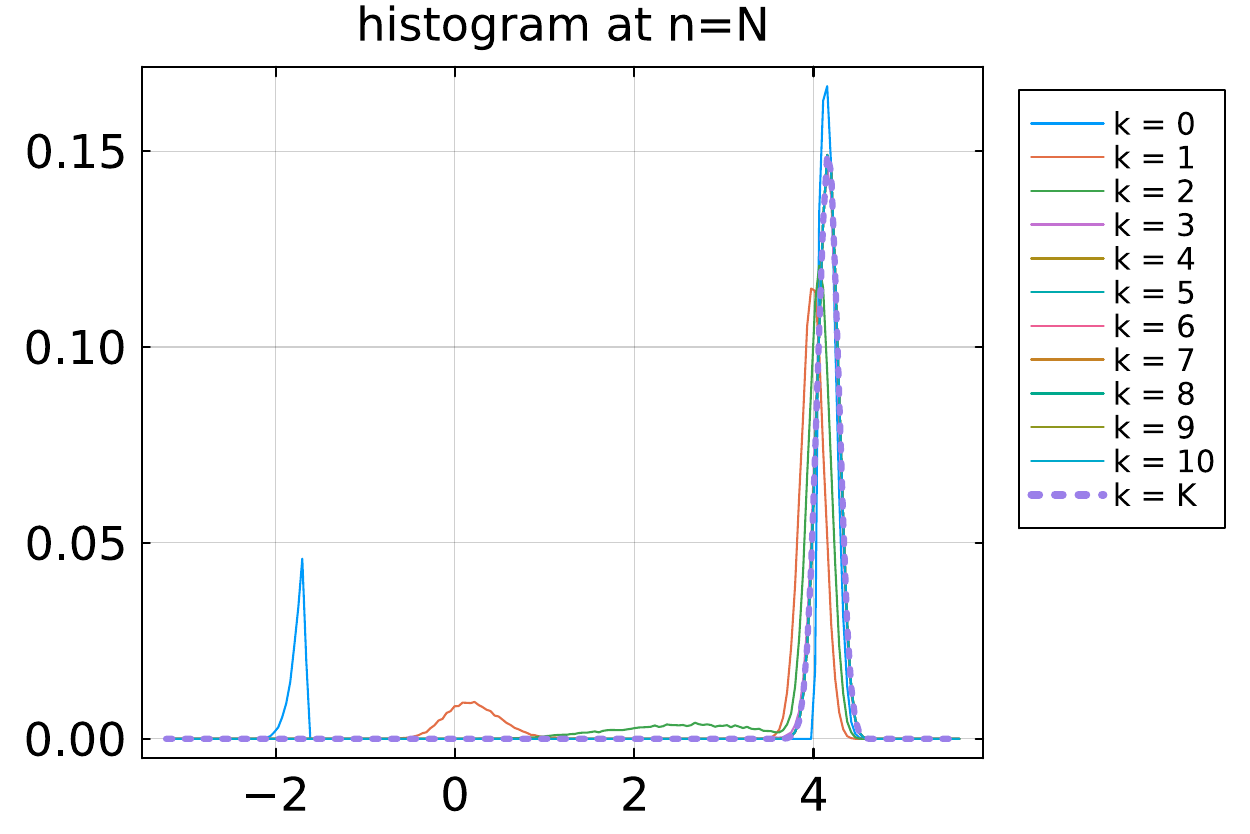}}%
	\caption{Setup 5}
  \end{subfigure}
\end{figure}

In order to judge the computational difficulty of the setups, we show in \cref{fig_overview_evolution_mean_var} the mean and variance of each setup as a function of time. (In \cref{fig_overview_evolution_PFs} (\cref{extra_figures_pf}) we plot the evolution of the particle fractions.)
We then present the convergence of the mean and variance for all these setups.
In \cref{fig_overview_convergence_histograms}, the histograms of the particle ensemble at $t=T$, for different iteration numbers of MC-moments Parareal (variant 2) is shown for all setups.

\newcommand{\FigureIntroConv}{Weak convergence of different variants of MC-moments Parareal. }

\begin{example}[Bimodal SDE, Setup 1] 
\label{Bimodal_SDE}
In this example, the initial condition contains particles in both regions of attraction.
In \cref{fig_bimodal_unimdal_coarse} we plot the convergence of different variants of MC-moments Parareal.
Variant 1 (unimodal MC-moments Parareal) converges extremely slowly to the level of statistical error (only in the last iteration).
Variant 2 (bimodal MC-moments Parareal) has a lower initial error, but still converges very slowly.
Variant 3 (learning-based) converges much faster: after a few iterations the statistical error is reached.
Variant 4 (with exact particle fractions) converges very quickly to the level of statistical error.

\begin{figure}[H]
\centering
\includegraphics[width=0.5\textwidth]{{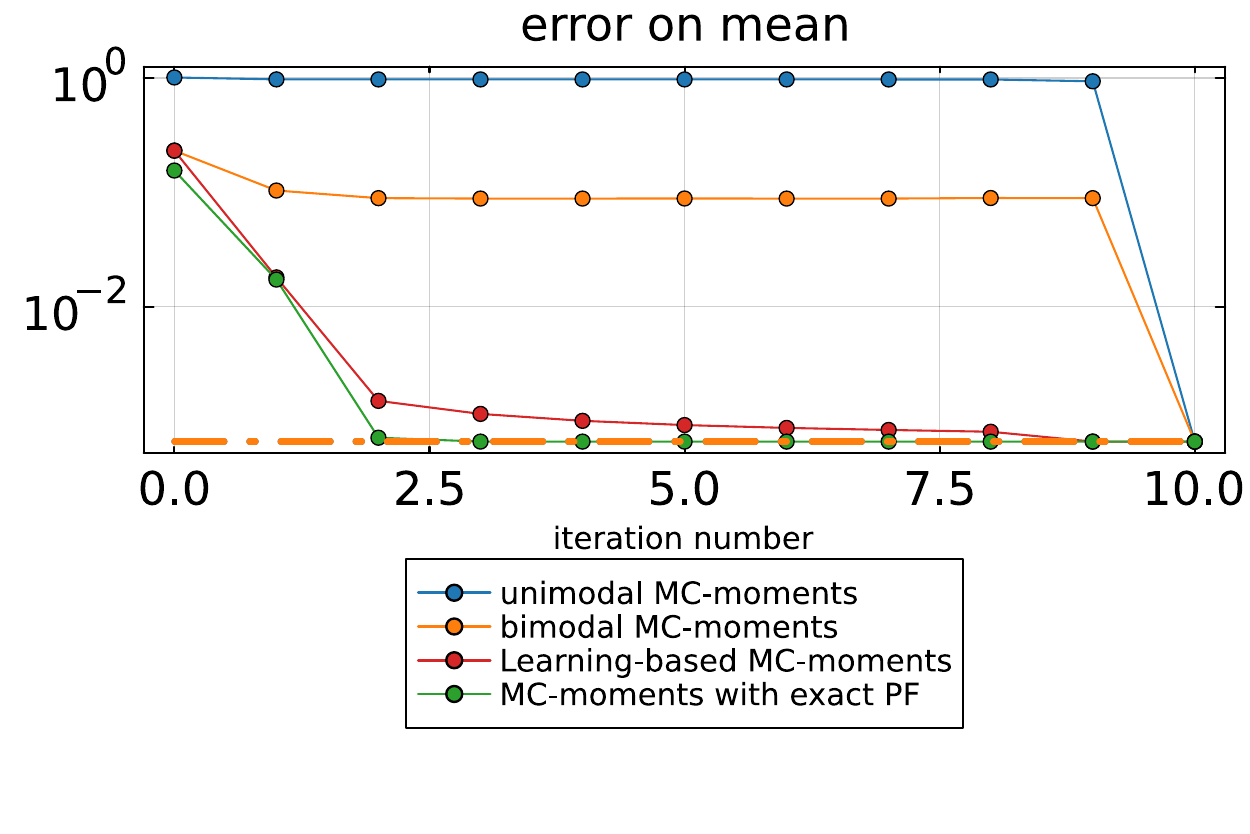}}%
\includegraphics[width=0.5\textwidth]{{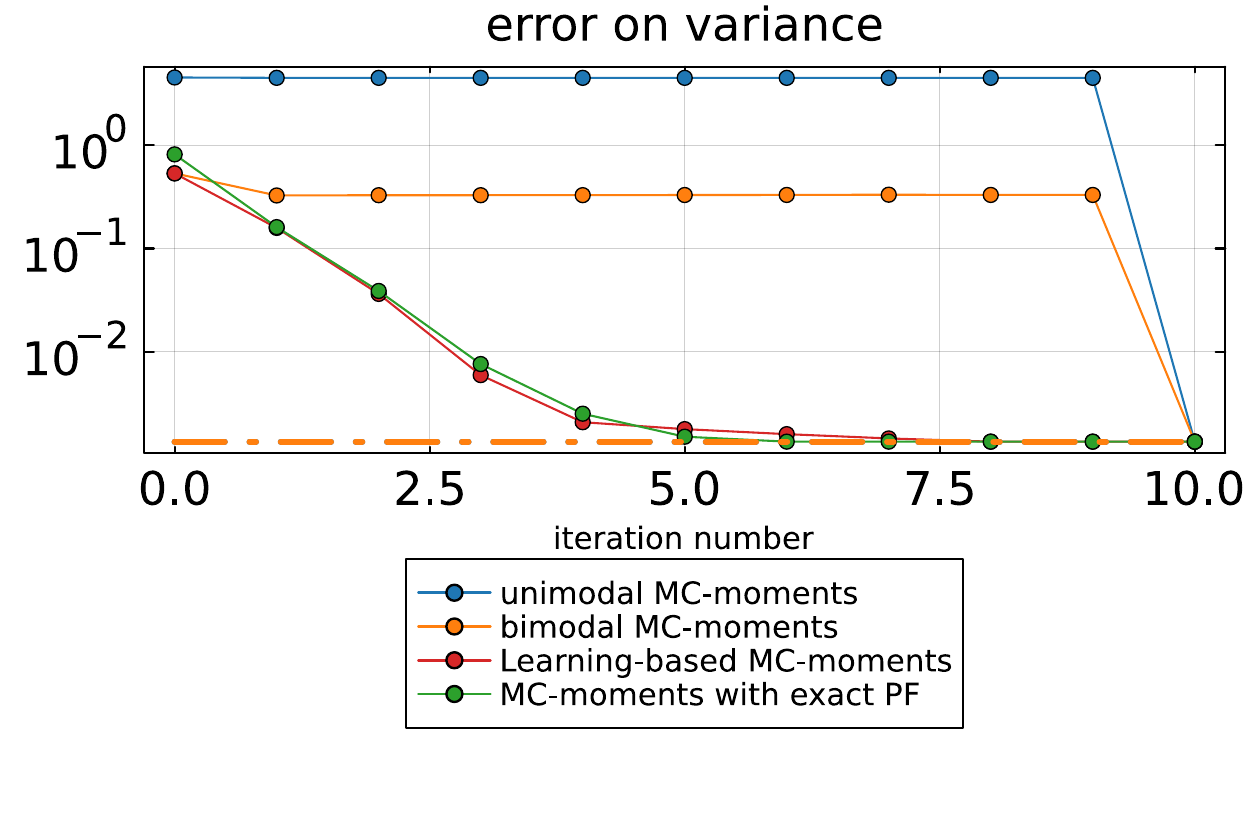}}%
\caption{\FigureIntroConv Setup 1.
}
\label{fig_bimodal_unimdal_coarse}
\end{figure}
\end{example}

\begin{example}[Bimodal SDE, Setup 2]
\label{Bimodal_SDE_unimodal_start}
In this example, the initial condition has all particles in one region of attraction. 
At time $t=T$, the number of particles in the other well is extremely low but nonzero; the influence of these particles on the overall mean and variance is significant.
In comparison to \cref{Bimodal_SDE} (setup 1), the mean and the variance reach their steady-state solution slightly faster in time, suggesting faster Parareal convergence (see \cref{Dahlquist_effect_time_horizon}).
In \cref{fig_bimodal_unimodal_initialisation} we show the convergence of different variants of MC-moments Parareal.
Overall, convergence is slightly faster than in \cref{Bimodal_SDE}, as expected. 
\begin{figure}[h]
\includegraphics[width=0.5\textwidth]{{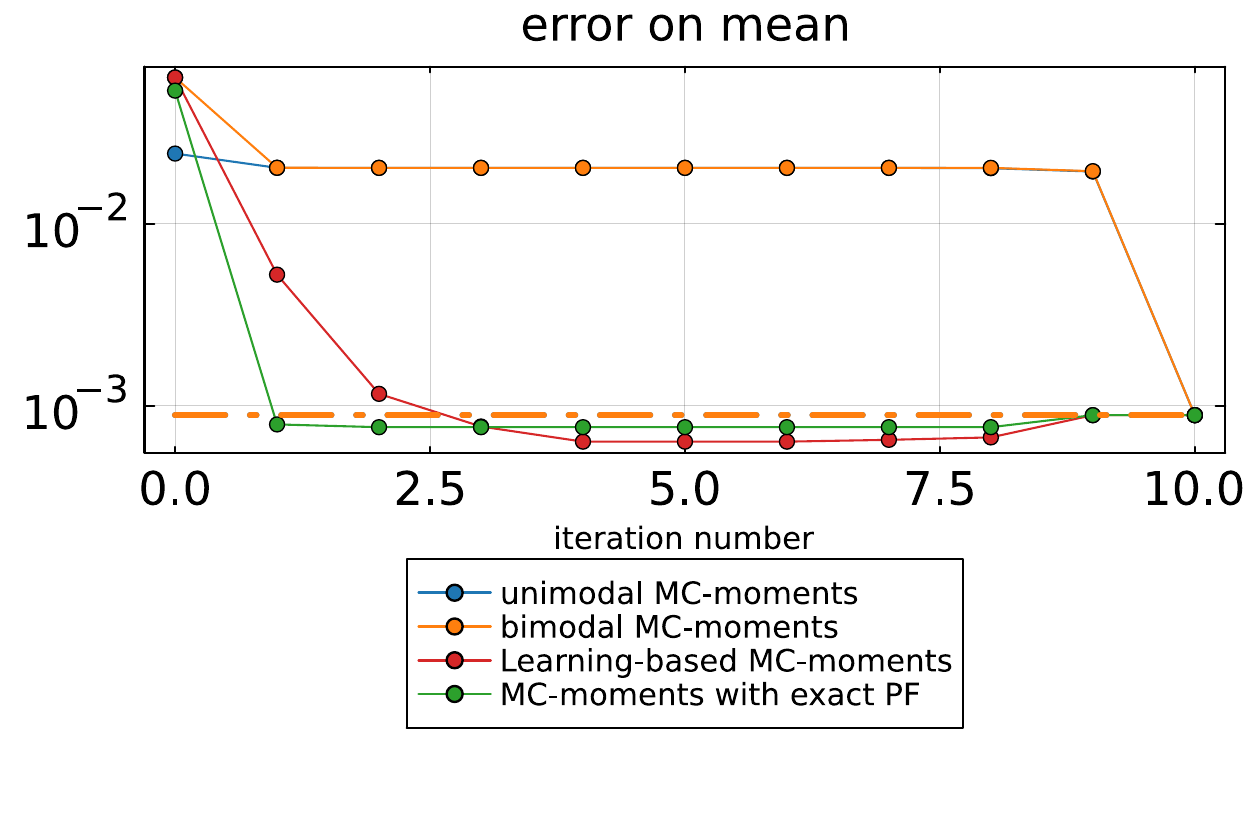}}%
\includegraphics[width=0.5\textwidth]{{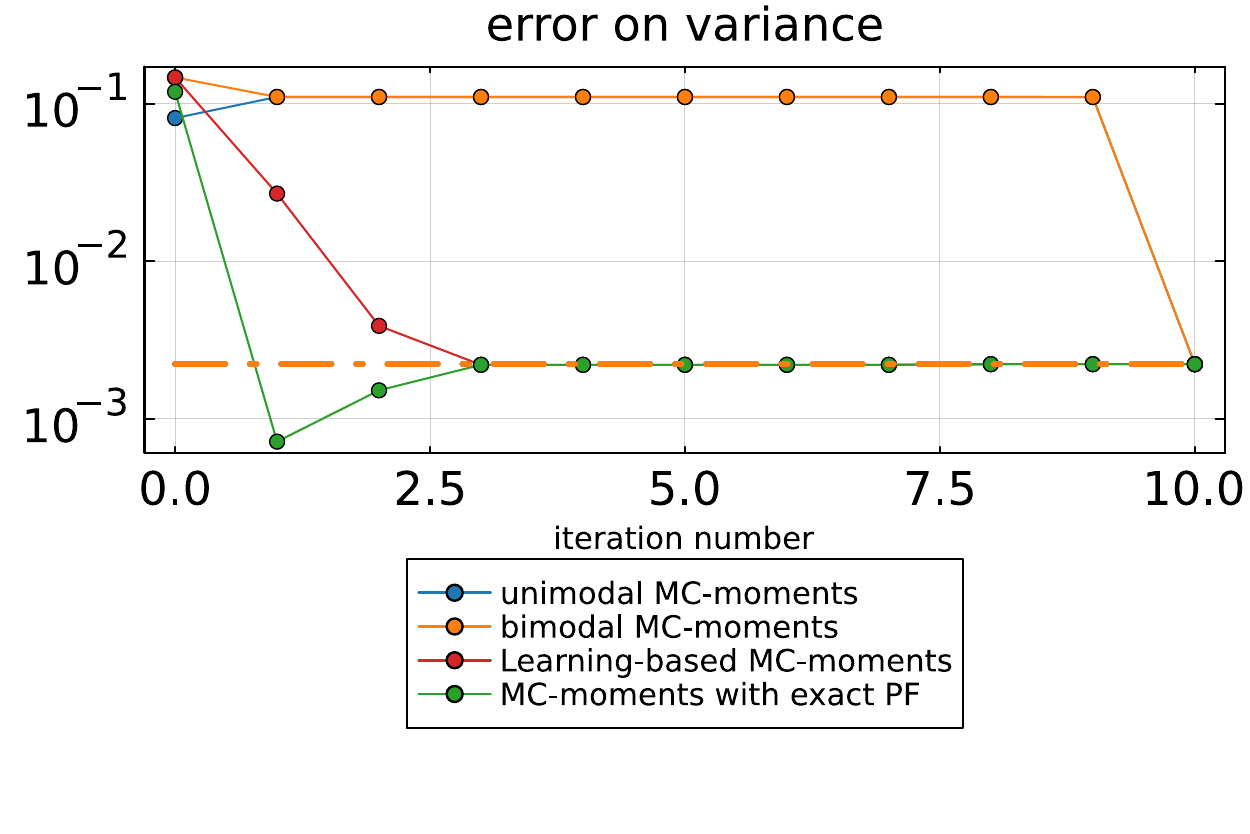}}
\caption{\FigureIntroConv Bimodal SDE setup 2.}
\label{fig_bimodal_unimodal_initialisation}
\end{figure}

The learning-based variant converges slightly slower than variant 4, but still much faster than the other variants.
In experiments not shown here, we observed that the extrapolation on the mean and the variance, based on the coarse propagator alone, proposed in \cref{definition_MC_moments_Parareal_with_learning}, is essential for the convergence of the learning-based variant. 
\end{example}

\begin{example}[Bimodal SDE, Setup 3]
\label{Bimodal_SDE_setup_3}
In this example, the initial condition contains all particles in one region of attraction, but there is a bit more diffusion with respect to setup 2.
This setup is challenging because a steady-state is not reached within the time horizon $[0,T]$, as can be seen in \cref{fig_overview_evolution_mean_var}.
From \cref{Dahlquist_effect_time_horizon}, we may expect that convergence will be slower then for setup 2.
The convergence behavior of the mean and variance is shown in \cref{fig_bimodal_setup_3}.
The convergence of all variants is quite slow, and can be explained by the absence of a steady state before time $T$.

\begin{figure}[h]
\includegraphics[width=0.5\textwidth]{{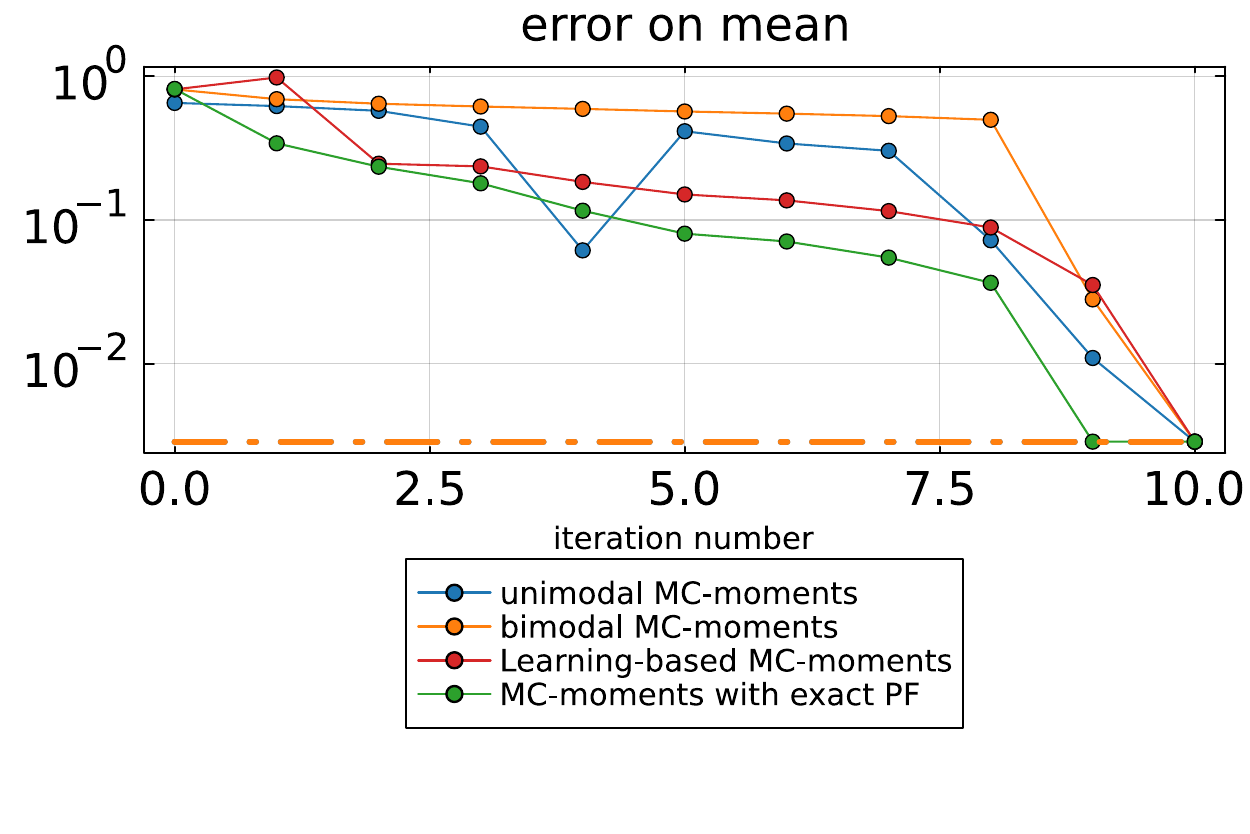}}%
\includegraphics[width=0.5\textwidth]{{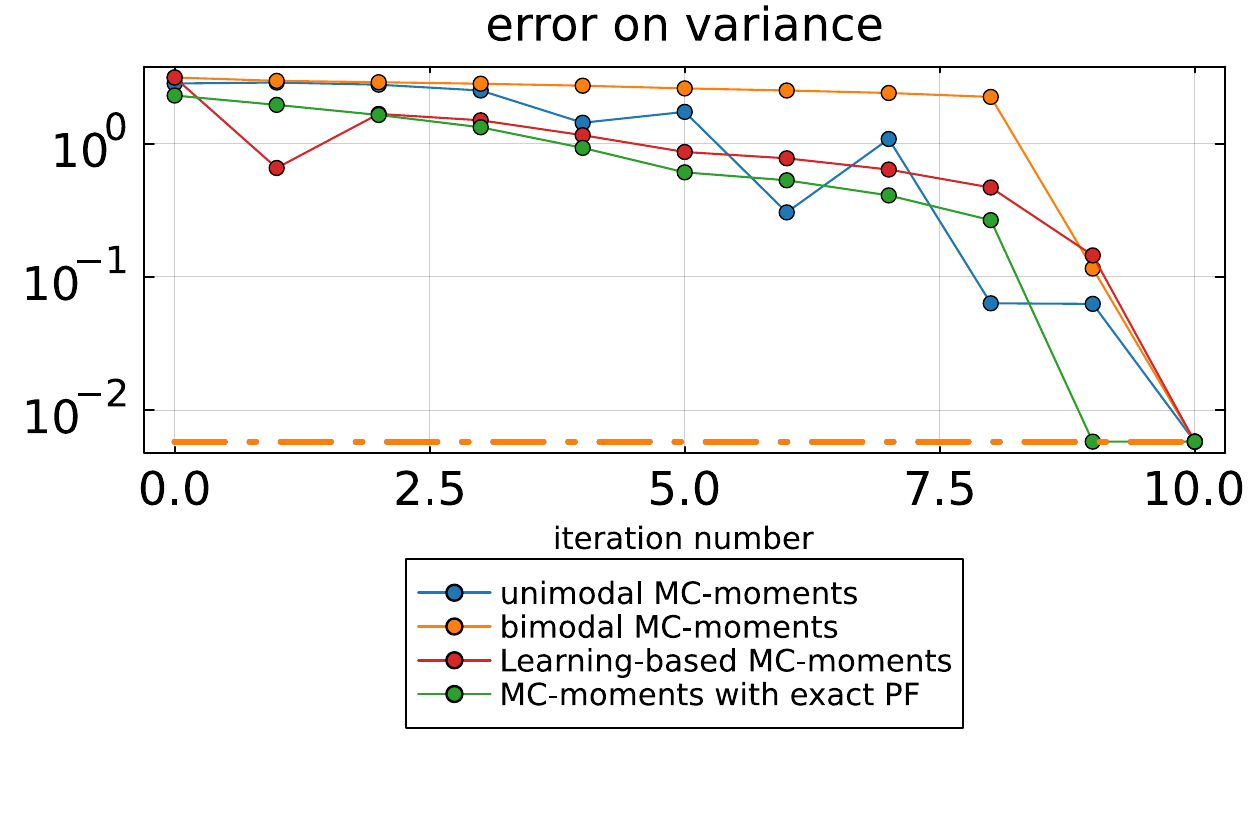}}
\caption{\FigureIntroConv Bimodal SDE setup 3.}
\label{fig_bimodal_setup_3}
\end{figure}
\end{example}

\begin{example}[Bimodal SDE, Setup 4]
\label{Bimodal_SDE_6}
In this example, there is more diffusion, as compared to \cref{Bimodal_SDE_setup_3}. 
This makes the locally unimodal particle distributions overlap at $t=T$, which makes it challenging to extrapolate the particle fractions using a set of simple linear uncoupled ODE (as is done in the learning-based method from \cref{definition_MC_moments_Parareal_with_learning}).
On the other hand, overall convergence is expected to be faster than in \cref{Bimodal_SDE_setup_3}, because the mean and variance reach a steady-state within the considered time-interval $[0,T]$ (see \cref{fig_overview_evolution_mean_var}).

The convergence behavior of the mean and the variance is shown in \cref{fig_bimodal_6}.
The unimodal Gaussian coarse model converges in a few iterations; even faster than the bimodal variant and the variant with learning. 
Indeed, the particle distribution at time $t=T$ is nearly unimodal and therefore is well described by a single Gaussian.

The convergence is faster than in \cref{Bimodal_SDE_unimodal_start} (setup 3), mainly due to the fact that the solution reaches a steady-state well within the time interval $[0, T]$.
Remarkably, for the error on the mean, variant 4 (learning-based) converges slower then variant 2 (where the particle fractions converge as in Parareal without coarse propagator).
An investigation of this effect is left as a subject of future research.

\begin{figure}[H]
\includegraphics[width=0.5\textwidth]{{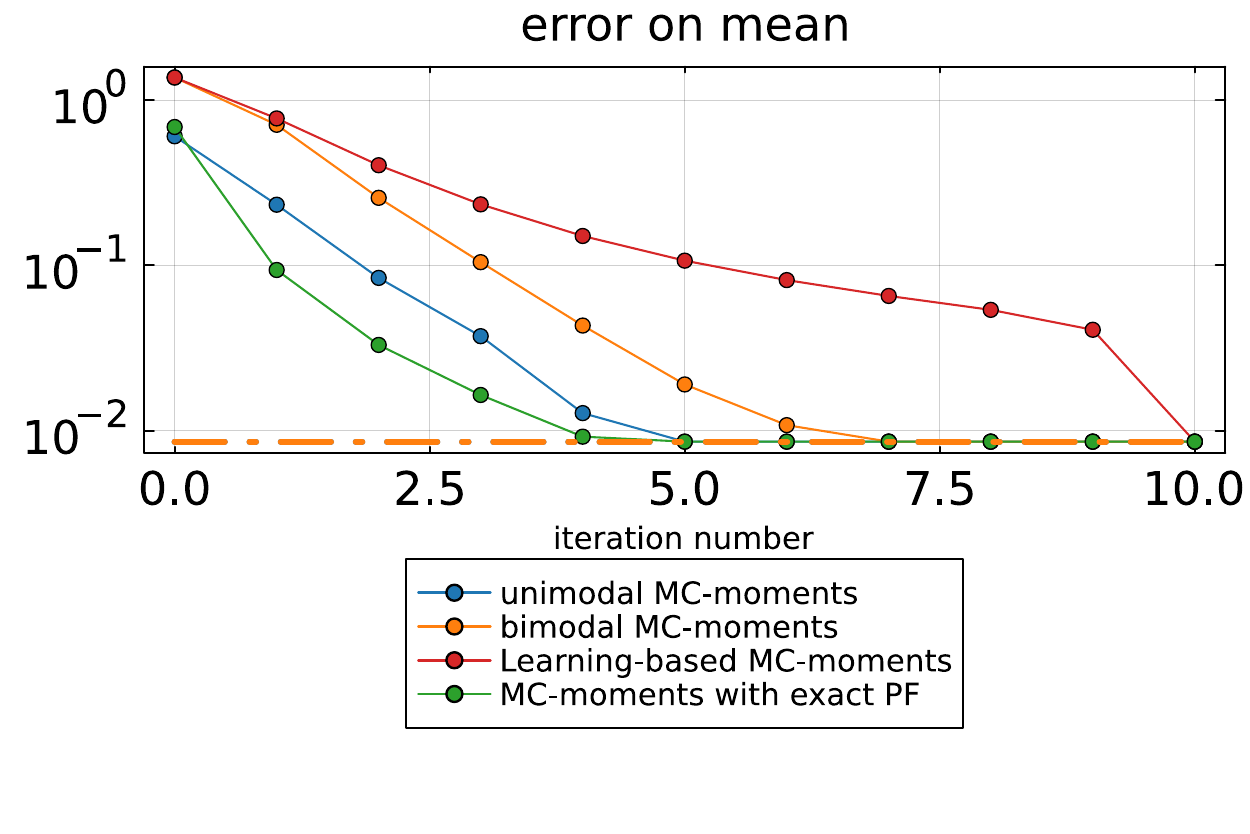}}%
\includegraphics[width=0.5\textwidth]{{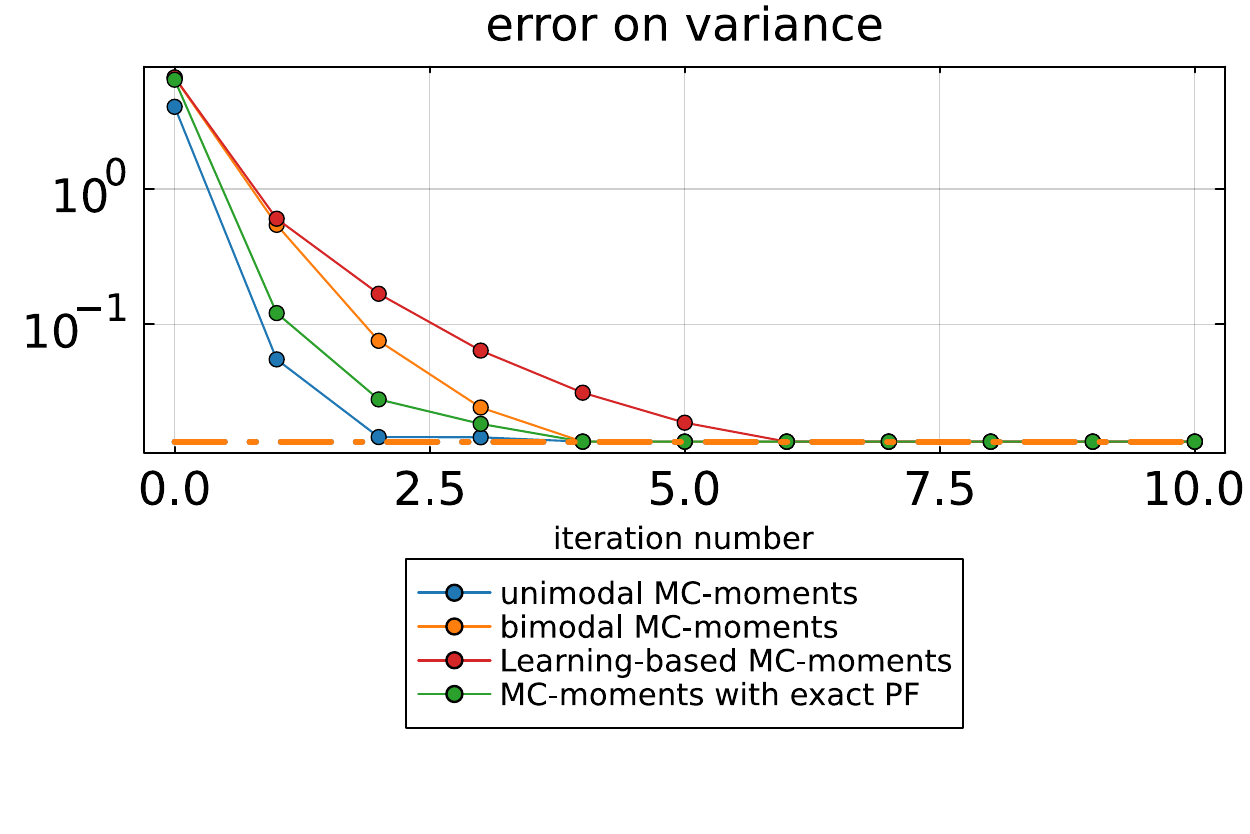}}
\caption{\FigureIntroConv Bimodal SDE setup 4.}
\label{fig_bimodal_6}
\end{figure}
\end{example}

\begin{example}[Bimodal SDE, Setup 5]
\label{Bimodal_SDE_MEAN_FIELD}
In this example, a nonzero mean-field interaction term is added to the SDE from setup 1 ($\beta  = 1$); the solution at $t=T$ is unimodal, and the mean and the variance reach a steady-state before time $T$.
Overall, we thus expect reasonably fast convergence, even with a unimodal coarse moment model.
In \cref{fig_bimodal_MEAN_FIELD} we show the convergence of the mean and the variance.
Variant 1 (unimodal MC-moments Parareal) converges quickly to the level of statistical error.
The other variants converge after one or two more iterations. 

\begin{figure}[h]
\includegraphics[width=0.5\textwidth]{{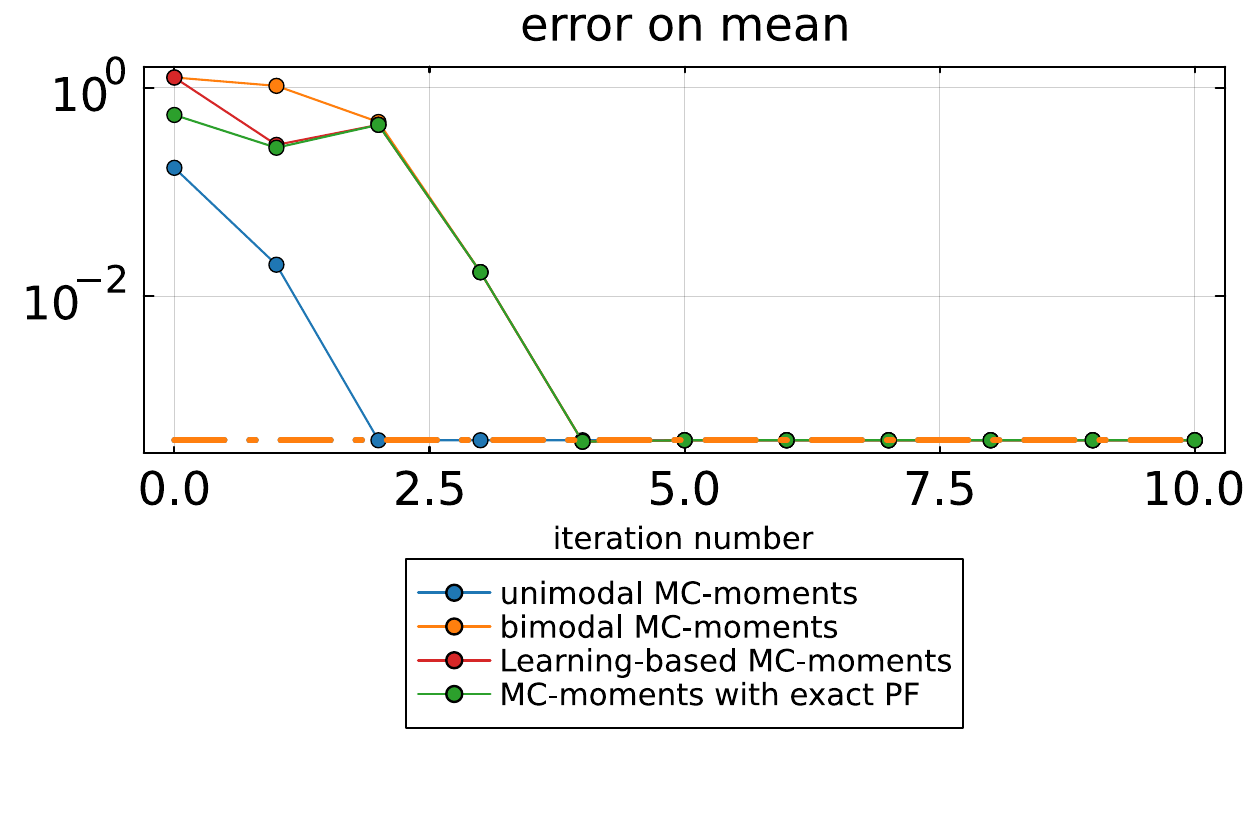}}%
\includegraphics[width=0.5\textwidth]{{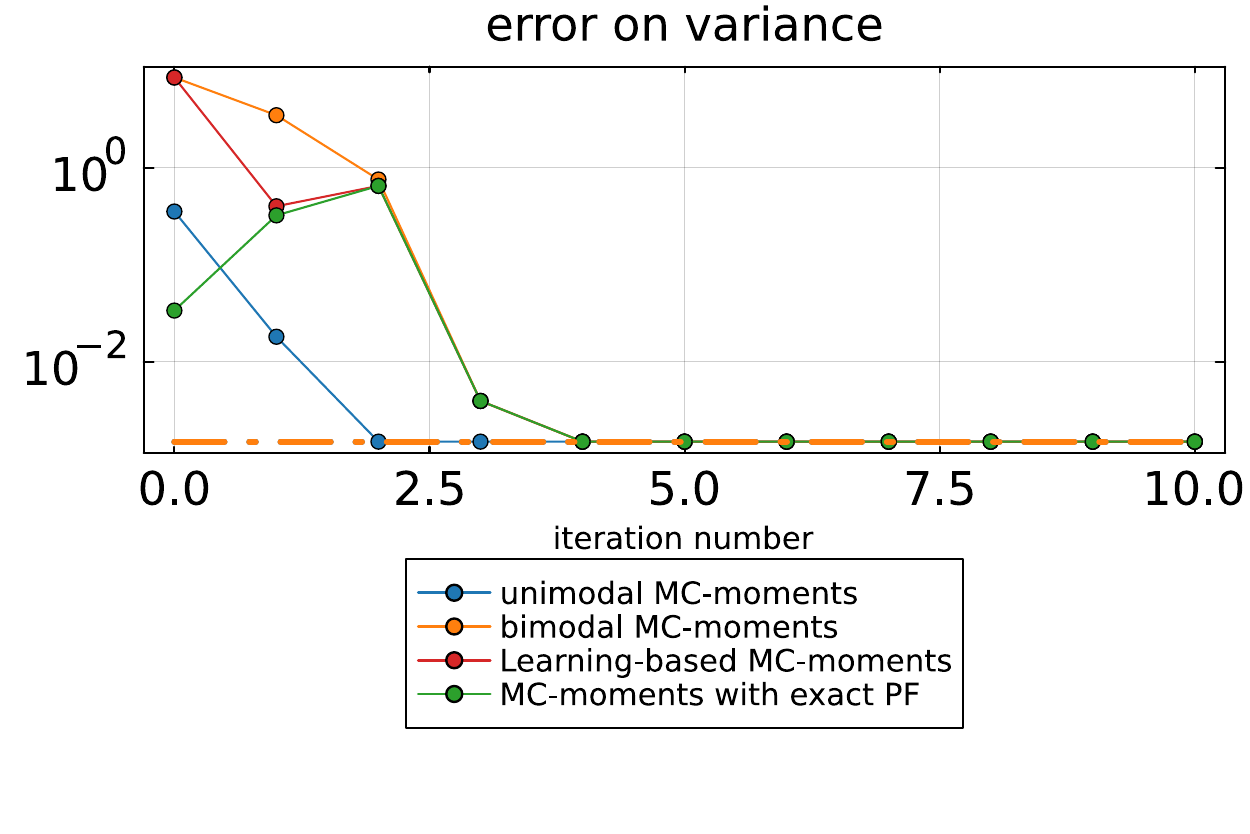}}
\caption{\FigureIntroConv Bimodal SDE setup 5.}
\label{fig_bimodal_MEAN_FIELD}
\end{figure}
\end{example}

\subsection{Effect of time horizon on convergence}
\begin{figure}[h]
\centering
\includegraphics[width=0.33\textwidth]{{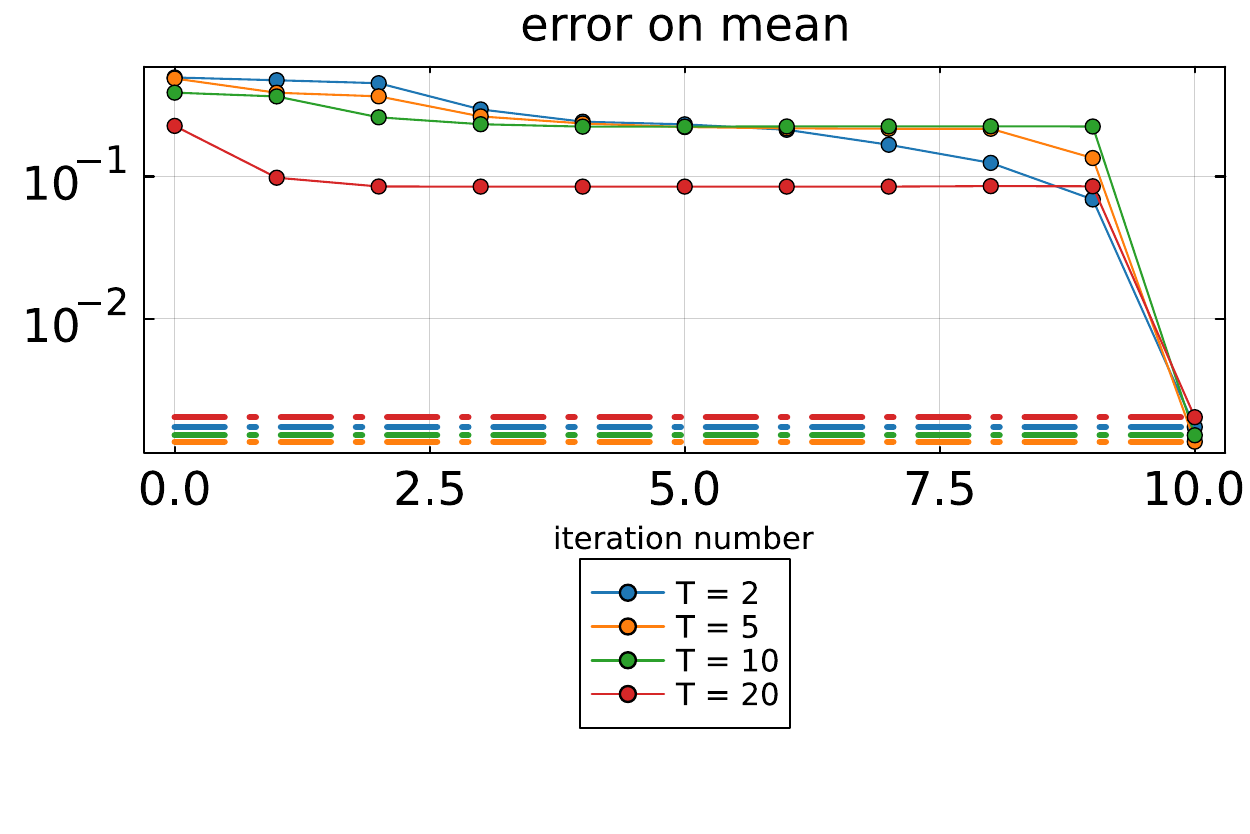}}
\includegraphics[width=0.33\textwidth]{{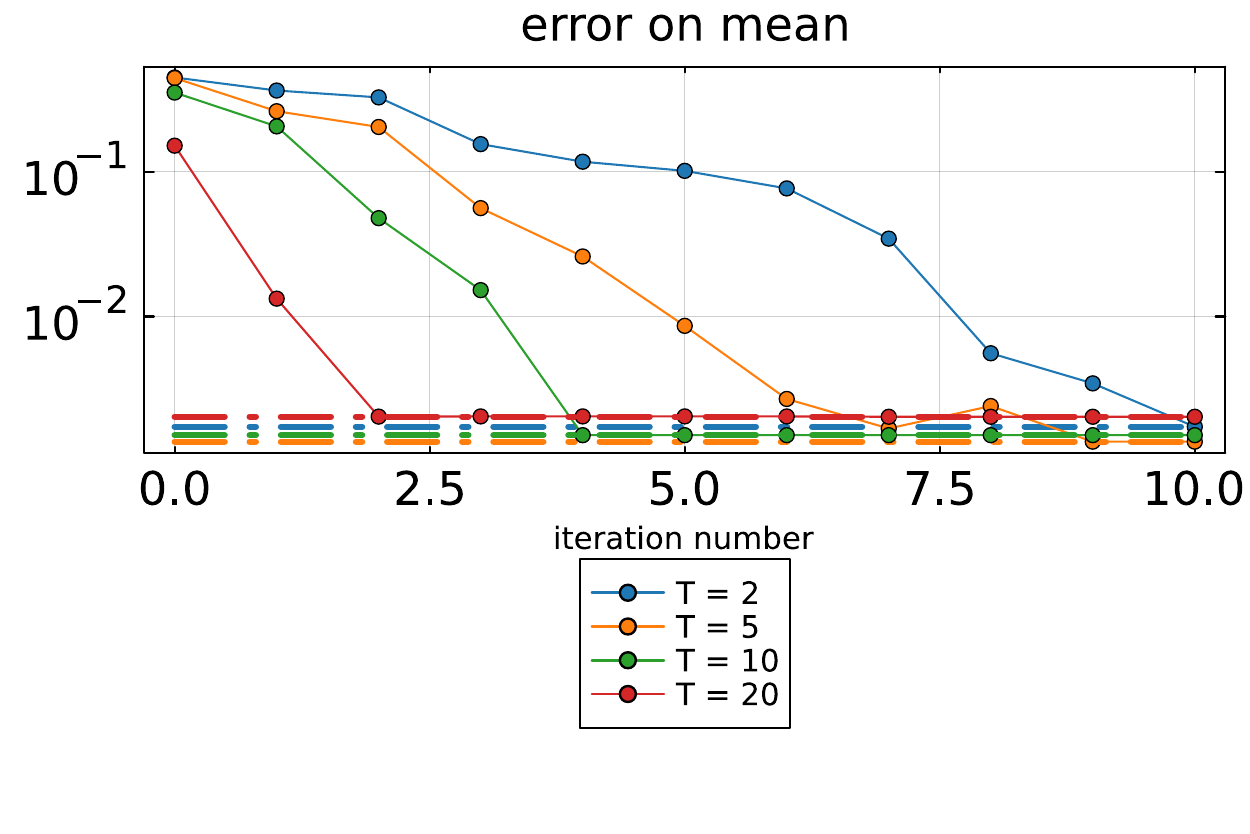}}%
\includegraphics[width=0.33\textwidth]{{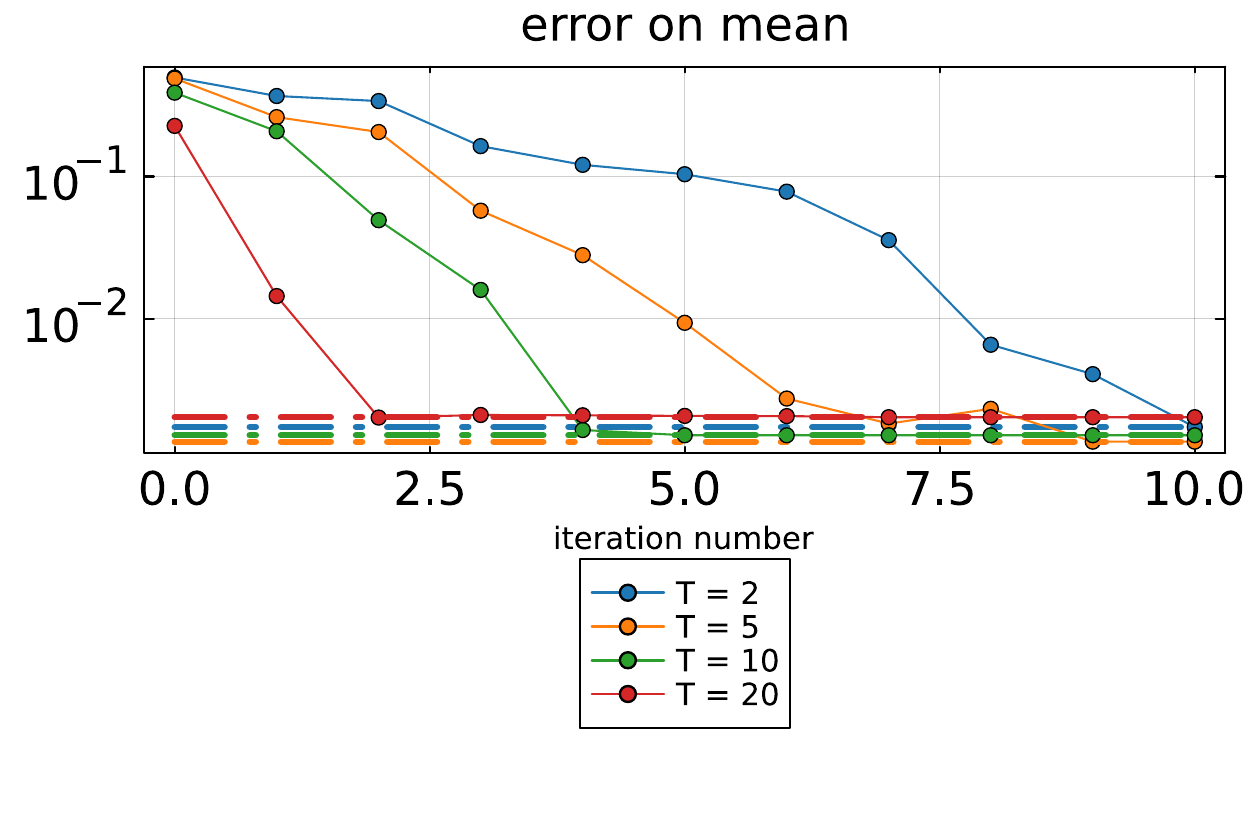}}
\begin{subfigure}[t]{.33\linewidth}
\includegraphics[width=\textwidth]{{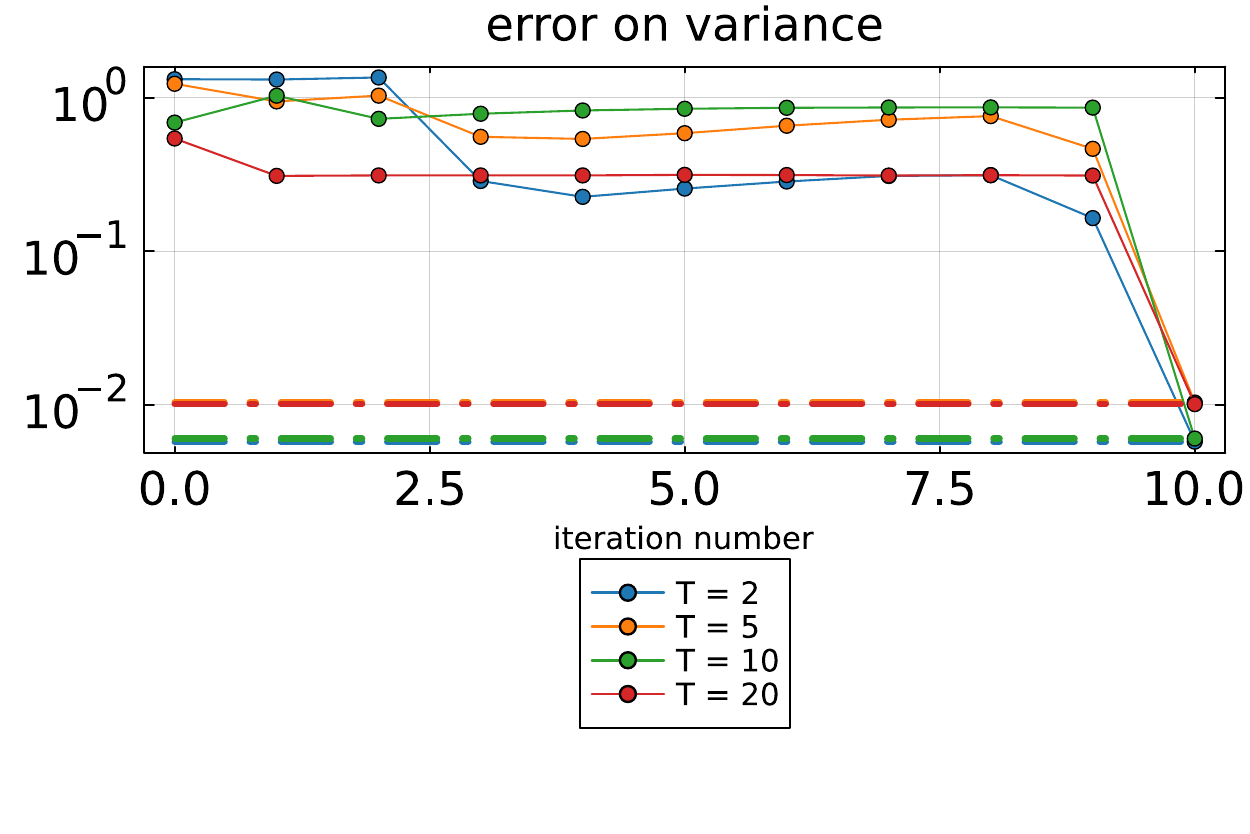}}
\caption{Variant 2 (multimodal)}
\end{subfigure}
\begin{subfigure}[t]{.33\linewidth}
\includegraphics[width=\textwidth]{{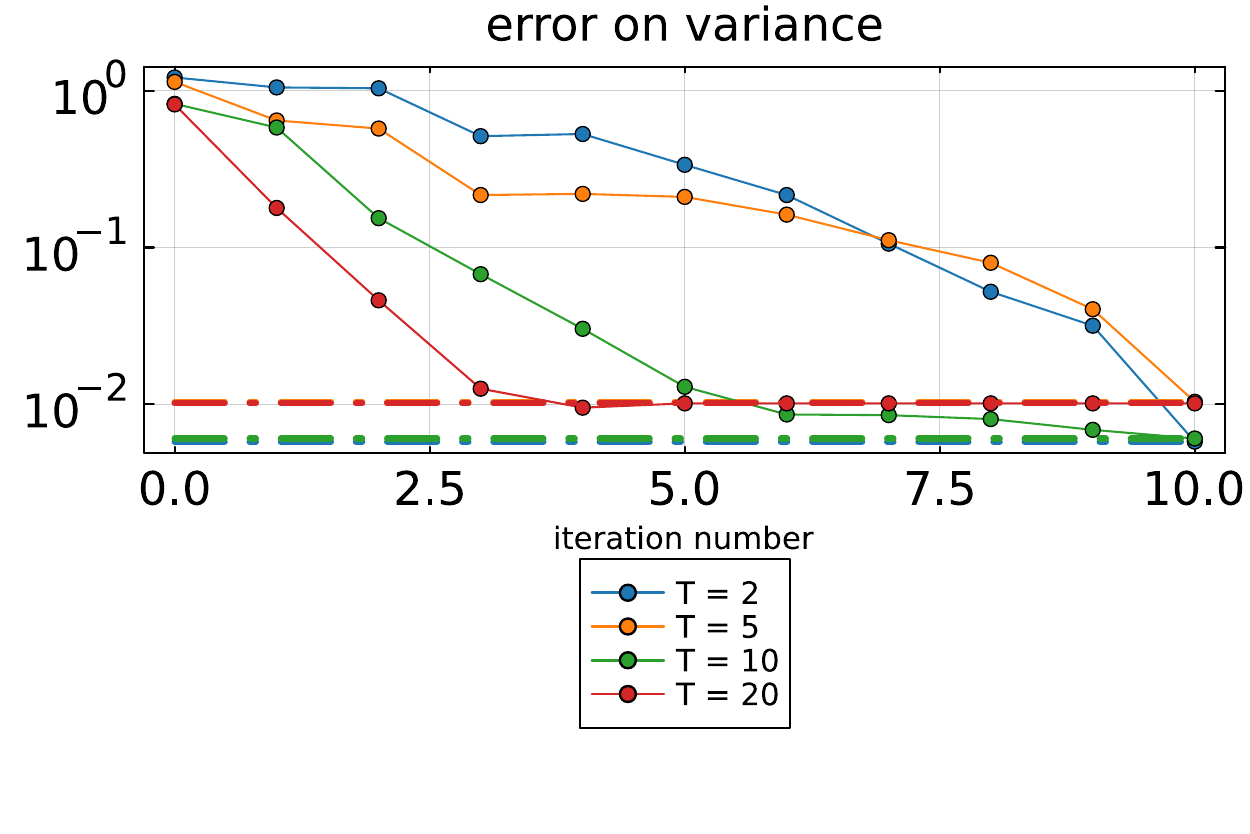}}%
\caption{Variant 3 (learning-based)}
\end{subfigure}%
\begin{subfigure}[t]{.33\linewidth}
\includegraphics[width=\textwidth]{{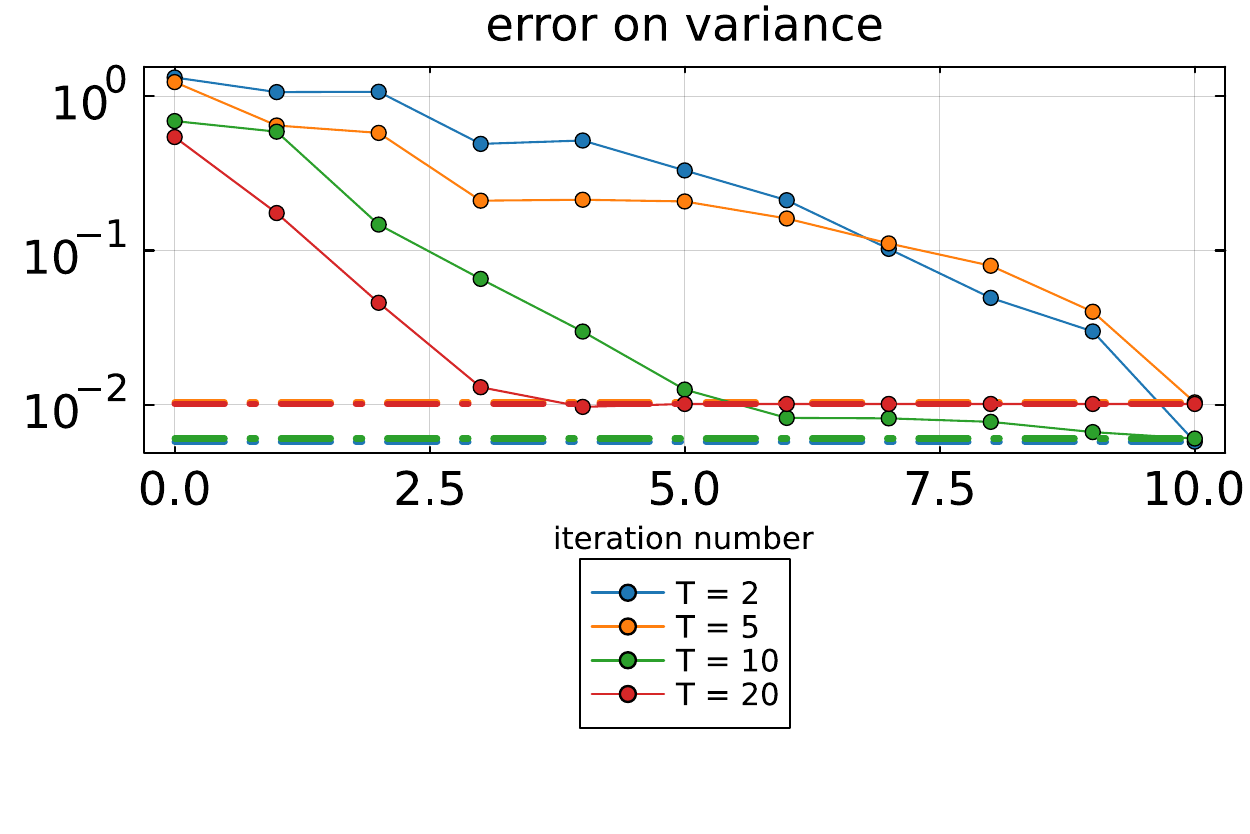}}
\caption{Variant 4 (exact PF)}
\end{subfigure}
\caption{Convergence of the mean (upper row) and of the variance (bottom row) for different MC-moments Parareal variants. Setup 1.}
\label{bimodal_effect_time_horizon}
\end{figure}

\mylabel{numerical_experiments_effect_time_horizon} \,
We aim to study the effect of the time horizon on the convergence of MC-moments Parareal, while keeping the number of subintervals $N$ constant. 
We studied a similar question in \cref{Dahlquist_effect_time_horizon} for the Dahlquist equation.
In \cref{bimodal_effect_time_horizon} we display the convergence for three variants of MC-moments Parareal, with parameters from setup 1, but with various different end times $T$.
An investigation of this effect for other setups is left as a subject for future work.

Variant 2 (bimodal MC-moments Parareal) reaches a plateau and only converges to the level of statistical error in the very last iteration, for all values of $T$.
From a comparison of variant 2 with variant 4 (exact particle fractions), we deduce that the occurence of this plateau is due to the slow convergence of the particle fractions, which are only exact in the last iteration. 

For setup 1, variant 3 (learning-based) converges similarly as variant 4 (exact particle fractions).
Augmenting the time horizon leads to faster convergence, as was also observed for the Dahlquist equation in \cref{Dahlquist_effect_time_horizon}.

\subsection{Discussion}
We numerically studied the convergence of MC-moments Parareal on a few simple scalar McKean-Vlasov SDEs. 
For bimodal SDEs, the key ingredient for the MC-moments Parareal method to work, is the introduction of several moment ODEs instead of a single ODE.
The error on the mean and the variance is strongly influenced by the error on the particle fractions.
Learning-based Parareal can improve the convergence, but is not always able to make MC-moments Parareal converge as quickly as when exact particle fractions are available.

\section{Conclusion and future work}
\label{SECTION_end}
We proposed a new Parareal method for scalar McKean-Vlasov SDEs. 
In the algorithm, the fine Parareal propagator is an Euler-Maruyama simulation of an ensemble of particles, 
while the coarse propagator is a moment model, describing the 
mean and the variance of the particle distribution.
Moment models have a greatly reduced computational cost as compared to stochastic interacting particle simulations, but at the expense of a model error. 

For bimodal SDEs, we use multiple moment models, each describing the mean and variance of the particle distribution in locally unimodal regions of the phase space (\cref{MC_moments_Parareal_for_multimodal}).
We also developed a variant that converges faster by iteratively learning a model for the distribution of particles in the all regions of the phase space (through a least-squares procedure) as the Parareal iterations progress (\cref{definition_MC_moments_Parareal_with_learning}).
\subsection{Future work}
\label{section_future_work}
We finish with some suggestions for future research. 

A next step would be to develop an adaptive variant of the MC-moments Parareal algorithm where multimodality is discovered on-the-fly, and where moment ODEs are then dynamically added or removed. 
This could alleviate the requirement of the availability of a priori estimations of the regions of attraction of the McKean-Vlasov SDE.

We limited the MC-moments Parareal algorithm to scalar McKean-Vlasov SDEs. 
A generalisation for (unimodal) SDEs in higher dimensions is available in \cite{bossuyt_new}, but it would be interesting to carry out more numerical tests.
Based on our results for scalar SDEs, we expect the convergence in high dimensions to crucially depend on the availability of a good coarse model and its coarse approximation.

With respect to implementation, it would be interesting to implement the algorithm in a parallel environment and test it on a massively parallel machine. 

Another extension of the MC-moments Parareal method is to allow the particles to have weights. 
This extension could be useful to apply the MC-moments Parareal algorithm to more general McKean-Vlasov SDEs, such as the Burgers equation with arbitrary initial condition (\cref{remark_about_weighted_particle_ensemble}). 
Also, weights could possibly be useful in the design of new matching operators and iterator functions.

\section*{Acknowledgments}
We sincerely thank the anonymous reviewers for their careful and thorough work.
Their suggestions have greatly ameliorated this text.
We also thank the following colleagues for discussions, suggestions, ideas and advice related to the content of this manuscript: Arne Bouillon, Tom Kaiser, Vince Maes,
Ausra Pogozelskyte (at the PinT workshop 2022), and
Hannes Vandecasteele.

\appendix
\section{Proof of \cref{lemma_what_happens_with_noise} (Convergence of MC-moments Parareal with a noisy fine solver)}
\label{proof_of_main_lemma}
\begin{proof}
Starting from \eqref{equation_starting_point_first_proof}, we write in vector notation, defining 
$\mathbf e^k = \left[ \begin{matrix} e^k_0 & \hdots & e^k_N \end{matrix} \right]$ 
and similarly 
$\bm \epsilon^{(k)}
= \left[ \begin{matrix} \epsilon^k_0 & \hdots & \epsilon^k_N \end{matrix} \right]$,
\begin{equation}
(I-G) \mathbf e^{k+1} = (F-G) \mathbf e^k + \bm \epsilon^{(k)}.
\end{equation}
Let the matrices $M \in \mathbb R^{N \times N}$ and $C \in \mathbb R^{N \times N}$ be defined as
\begin{equation}
M = 
\left[ 
\begin{matrix}
I & \hdots & \\
-G & I & \hdots & \\
\vdots & &  \vdots & \\
& & & -G &  I & \\
0 & & & & -G &  I 
\end{matrix} 
\right];
\qquad
C = 
\left[ \begin{matrix}
0 & \hdots & & \\
(F-G) & \hdots & \\
\vdots & &  \vdots & \\
\hdots & &  \hdots & & 0 \\
0 & & & (F-G) & 0 
\end{matrix} 
\right] 
\end{equation}
Let the matrices $H \in \mathbb R^{N \times N}$ and $D \in \mathbb R^{N \times N}$ be defined as follows:
\begin{equation}
H = \left[ 
\begin{matrix} 
0 & 0 & \hdots & 0 & 0 \\ 
I & 0 & \hdots & 0 & 0 \\
\hdots \\
G^{N-1} & \hdots & I & 0 & 0 \\
G^{N-2} & G^{N-1} & \hdots & I & 0 \\
\end{matrix}
\right];
\qquad 
D = \left[ 
\begin{matrix}
F-G & 0 & \hdots & 0 & 0 \\ 
0 & F-G & \hdots & 0 & 0 \\
0 & & \hdots & 0 & 0 \\
\hdots & & & & \hdots \\
0 & 0 & \hdots & 0 & F-G \\
\end{matrix}
\right].
\end{equation}
Then, it holds that $M^{-1}C = HD$ and $CM^{-1}=DH$.

Thus $M \mathbf e^{k+1} = C \mathbf e^k + \bm \epsilon^{(k)}$, and therefore $\mathbf e^{k+1} = (M^{-1}C)^k \mathbf e^k + M^{-1} \bm \epsilon^{(k)}$.
Then it holds, by induction, that
\begin{equation}
\begin{aligned}
\mathbf e^{k} 
&= (M^{-1}C)^k 
\mathbf  e^0 
+ 
\sum_{j=0}^{k-1} (M^{-1}C)^{j} M^{-1} 
\mathbf \mathbf \bm \epsilon^{(j)} \\
&= 
(HD)^k 
\mathbf e^0 + \sum_{j=0}^{k-1} (HD)^{j} M^{^1} \mathbf \mathbf  \bm \epsilon^{(j)}.
\end{aligned}
\end{equation}

Let us now take the absolute value elementwise on both sides, and apply the expectation operator on both sides (making use of the triangle inequality in the right-hand side), then we obtain for $\xi^k = \mathbb E[|\mathbf e^k|]$, where the expectation and the absolute value are applied elementwise:
\begin{equation}
\mathbf{\xi}^k
\leq 
(HD)^k \xi^0 
+
\sum_{j=0}^{k-1} (HD)^j M^{-1} \mathbb E [|\bm \epsilon^{(j)}|].
\end{equation}
Now we use the triangle inequality and we bound the norms of the matrices $H$ and $D$ and $M^{-1}$ from above.
We also use that $\norm{\mathbb E [|\bm \epsilon^{(j)}|]} \leq \varepsilon$ for all $j$.
We thus have
\begin{equation}
\norm{
\mathbf \xi^k }_{\infty} 
	\leq
\norm{(HD)^k }
\norm{ \mathbf \xi^{0} }
+ 
\sum_{j=0}^{k-1} \norm{ (HD)^{j} M^{-1} }
\varepsilon.
\end{equation}
Now we derive a superlinear and a linear bound based on results from \cite{gander_analysis_2007}.
\begin{itemize}
\item For the linear bound, we write
$\| (HD)^k \| 
\leq
\norm{H^k} \norm{D^k} 
\leq \left( \frac{|F-G|}{1-|G|}\right)^k$ 
and 
$\| (HD)^{j} M^{-1} \| 
\leq 
\| (HD)^{j} \| \| M^{-1} \| 
\leq 
\left( \frac{|F-G|}{1-|G|}\right)^{j} \frac{1}{1-|G|}$.

\item For the superlinear bound, 
$\norm{(HD)^k} \leq (F-G)^k \binom{N-1}{k}$ 
and 
$\norm{(HD)^{j} M^{-1}} \leq (F-G)^{j} \binom{N-1}{j} \frac{1-|G|^{N-1}}{1-|G|}$
\end{itemize}
\end{proof}

\noindent
This leads to equations \eqref{superlinear_bound_with_noise} and \eqref{linear_bound_with_noise}.

\section{Proof of \cref{lemma_exactness_moment_model} (Exactness of moment model for linear SDEs)} \label{proof_exactness_linear_McK_V_SDE}
\begin{proof}
The proof is based on a similar procedure as mentioned in \cite{kloeden_gauss_quadrature_2017} for a linear Ornstein-Uhlenbeck SDE. See also \cite[Section 8.4 and 8.5]{arnold_stochastic_1974}. 
Writing It\^{o}'s lemma for the SDE \eqref{SDE_class_1} and the test functions $f(x)=x$ and $f(x)=(x-M)^2$, where $M = \mathbb{E}[X]$, and then taking expectations leads to the desired result, namely the moment model \eqref{moment_model_class_1}
with $a(X,\mathbb{E}[X],t) = A(t) X + A_E \mathbb{E}[X] + A_0(t)$, $b(X,\mathbb{E}[X],t) = B(t) X + B_E \mathbb{E}[X] + B_0(t)$, thus $a_X = A(t)$ and $b_X = B(t)$ and $b_{XX} = 0$.
We first write It\^{o}'s lemma with the test function $f(x) = x$:
\begin{equation}
dX 
= \left(A(t)X + A_E(t)\mathbb{E}[X],t + a(t) \right) dt + \left( B(t) X + B_E(t) \mathbb{E}[X] + b(t) \right) dW.
\end{equation}
Taking expectation immediately leads to the evolution ODE for $M$ in \eqref{moment_model_class_1}. 
For the test function $f(x)=(x-M)^2$ we write It\^{o}'s lemma as follows:
\begin{equation}
\begin{aligned}
d(X-M)^2 
&= \left[ 2(X-M) 
\left( A(t)X + A_E(t)\mathbb{E}[X] + a(t) \right) 
+  \left( B(t) X + B_E(t) \mathbb{E}[X] + b(t) \right)^2 \right] dt \\
&+ 2(X-M) \left( B(t) X + B_E(t) \mathbb{E}[X] + b(t) \right) dW.
\end{aligned}
\label{tussenresultaat_proof}
\end{equation}
Equivalently, after rewriting \eqref{tussenresultaat_proof} and by using the fact that $(BX + Q)^2 = B^2(X-M)^2 + Q^2 + 2BXQ + B^2(X-M)M + B^2M X$, we obtain
\begin{equation}
\begin{aligned}
d(X-M)^2 
&= \left[ 2(X-M)^2 A(t) + 2M(X-M)A(t) + 2(X-M)(A_E(t)\mathbb{E}[X] + A_0(t)) \right. \\
&\quad 
+ B(t)^2(X-M)^2 
+ (B_E(t) \mathbb{E}[X] + B_0(t))^2 
+ 2B(t) (B_E(t) \mathbb{E}[X] + B_0(t))^2 X  \\
&\quad \left. +  
(X-M)B(t)^2M 
+ B(t)^2MX \right] dt   \\
&\quad + 2(X-M) \left( B(t) X + B_E(t) \mathbb{E}[X] + b(t) \right) dW.
\end{aligned}
\label{ander_tussenresultaat}
\end{equation}
Then, upon taking expectations of \eqref{ander_tussenresultaat} and applying the martingale property and the fact that $\mathbb{E}[\kappa (X-M)]=0$, where $\kappa$ is a constant, the exact moment ODE for the variance $\Sigma$ is obtained:
\begin{equation*}
\begin{aligned}
d \mathbb{E}[(X-M)^2] 
&= 2 A(t) \Sigma 
+ B(t)^2 \Sigma 
+ (B_E(t) M + B_0(t))^2
+ 2B(t) (B_E(t) M + B_0(t))^2 M
+ B(t)^2M^2 \\
&= (2 A(t)  + B(t)^2) \Sigma 
+ (B(t)M + B_E(t)M + B_0(t))^2.
\end{aligned}
\end{equation*}
This corresponds with the ODE for $\Sigma$ in \eqref{moment_model_class_1} and concludes the proof.
\end{proof}

\section{Proof of exactness of moment model for polynomial drift model}
\label{appendix_exactness_polynomial_drift_model}
\begin{example}[Polynomial drift]
\label{example_polynomial_drift}
The polynomial drift SDE is taken from \cite{son_doan_mean_square_2015}
\begin{equation}
\begin{aligned}
dX &= (\alpha X +  \mathbb{E}[X] - X \mathbb{E}[X^2])dt + X dW \\
X(0) &= 1.
\end{aligned}
\end{equation}

In fact, the model \eqref{moment_model_class_1} corresponds to the exact moment ODEs derived in \cite[equations (15) and (16)]{son_doan_mean_square_2015} (see appendix \ref{appendix_exactness_polynomial_drift_model}).
\end{example}

\begin{proof}
The moment model \eqref{moment_model_class_1} can be derived as follows: 
$a(M,f(M) = (\alpha + 1)M - M \mathbb{E}[X^2] = (\alpha + 1)M - MS$, $a_X = \alpha-\mathbb{E}[X^2] = \alpha- S^2 = \alpha - (\Sigma + M^2)$, $a_{XX} = 0$, $b(M,f(M)) = M$, $b_X = 1$, and $b_{XX} = 0$.
Thus we have 
\begin{equation}
\begin{aligned}
\frac{dM}{dt} &= (\alpha+1)M - M(\Sigma + M^2) \\
\frac{d\Sigma}{dt} &= \left[2(\alpha-(\Sigma+M^2)) + 1 \right] \Sigma + M^2.
\end{aligned}
\end{equation}
Now we can use the relation 
\begin{equation}
\frac{dS}{dt} 
= \frac{d}{dt} \left( \Sigma + M^2 \right) 
= \frac{d\Sigma}{dt} + 2M \frac{dM}{dt}.
\end{equation}
Thus 
\begin{equation}
\begin{aligned}
\frac{dS}{dt} &= \underbrace{\left[2(\alpha-S) + 1 \right] (S- M^2) + M^2}_{d\Sigma/dt} 
+ 2M \left( (\alpha+1)M - MS \right) \\
&= \left[ 2 \alpha + 1 \right] (S-M^2) - 2S(S-M^2) + M^2 + 2(\alpha+1)M^2 - 2M^2S \\
&= \left[ 2 \alpha + 1 \right] (S-M^2) - 2S(S-M^2) + M^2 + \left[ 2\alpha+1 \right] M^2 + M^2 - 2M^2S \\
&=  \left[ 2 \alpha + 1 \right] S - 2 S(S-M^2) + 2M^2 - 2M^2S \\
&= \left[2\alpha+1 \right] S  + 2M^2 - 2 S ^2.
\end{aligned}
\end{equation}

The moment equations (for the central moment) with $m=\mathbb E[X]$ and $S=\mathbb E[X^2]$ are \cite{son_doan_mean_square_2015}
\begin{equation}
\begin{aligned}
\frac{dM}{dt} &= (\alpha+1)M - M S \\
\frac{dS}{dt} &= (2\alpha+1) S  + 2 M^2 - 2 S ^2.
\end{aligned}
\end{equation}
This proves that the moment model with the moment model is exact.
\end{proof}

\section{Derivation of the moment model}
\label{appendix_derivation_moment_model}
We fist write It\^{o}'s lemma (see, for instance, \cite{higham_kloeden_introduction_2021}).
If the function $\Phi(x,t)$ is applied on the SDE \eqref{SDE_class_1} (here without mean-field interaction), we obtain
\begin{equation}
d \Phi(X,t) = \left( \Phi_t(X,t) 
+ \Phi_{x}(X,t) a(X) 
+
\frac{1}{2}
\Phi_{xx}(X,t)  b(X)^2 \right)dt 
+ \Phi_x(X,t)b(X,t) dW.
\end{equation}

\paragraph{For the mean}
For the mean, we obtain, using that $\Phi(x,t) = x$ and the martingale property
\begin{equation}
d X = a(X) dt + b(X) dW.
\end{equation}
Then, using the martingale property, writing a Taylor expansion and taking the expected value,
\begin{equation}
\begin{aligned}
\frac{dM}{dt}
&\approx	 
\mathbb E 
\left[
a(M) 
+ (X-M) a'(M) 
+ \frac{(X-M)^2}{2} a''(M)
+ \mathcal O((X-M)^3)
\right] \\
&\approx
a(M) + \frac{\Sigma^2}{2} a''(M).
\end{aligned}
\end{equation}

\paragraph{For the variance}
For the variance, we obtain, using It\^{o}'s rule with $\Phi(x,t) = (x-M)^2$
\begin{equation}
d (X-M)^2 = 
\left( 
2(X-M) a(M) 
+ b(X)^2 \right)dt 
+ 2(X-M) b(X,t) dW,
\end{equation}
and then, using the martingale property, writing a Taylor series and taking the expected value
\begin{equation}
\begin{aligned}
\frac{dM}{dt}
&\approx
\mathbb E 
\left[
2(X-M) 
\left( 
a(M) 
+ (X-M) a'(M) 
+ \frac{(X-M)^2}{2} a''(M)
+ \mathcal O((X-M)^3) 
\right) \right. \\
& \quad + \left.
\left( 
b(M) + (X-M) b'(M) + \frac{(X-M)^2}{2} b''(M) 
+
\mathcal{O}\left( (X-M)^3 \right)
\right)^2  \right] \\
&\approx
\left( 2 a'(M) + b'(M)^2 \right) 
\Sigma
+ b(M)^2.
\end{aligned}
\end{equation}

\section{Expectation of a nonlinear quantity of interest}
\label{appendix_estimating_nonlinear_QoI}
We write a Taylor series around $M = \mathbb E[X]$, and variance $\Sigma = \mathbb E[(X-M)^2]$ and then taking the expectation.
\begin{equation}
\Phi(X) 
\approx
\Phi(M) 
+ (X-M)\Phi'(M)
+ \frac{(X-M)^2}{2} \Phi''(M).
\end{equation}

\section{Proof of \cref{property_of_multimodal_linear_moment_model} (Exactness of multimodal Taylor series-based moment model)}
\label{nog_een_bewijsken}
Let us consider a linear ODE of the form
\begin{equation}
\frac{du}{dt} = Au + b,
\qquad u(0) = 0.
\label{moment_model_for_proof}
\end{equation}
It holds that, if $A^{-1}$ exists,
\begin{equation}
u(t) = e^{At} \left[ u_0 + A^{-1}b \right] - A^{-1}b.
\label{jaja}
\end{equation}
Here we need to prove that, when the model \eqref{moment_model_for_proof} is applied on the mean (and variance), the result is equivalent to applying the model on an ensemble of local means $M_i$ and variances $\Sigma_i$  and then combining these results with \eqref{mixture_of_distributions_mean_variance}.
\newline
\textbf{For the mean}, we have that, on the one hand, using a unimodal exact moment model
\begin{equation}
M(t) = e^{At} \left[ M_0 + A^{-1}b \right] - A^{-1}b.
\end{equation}
On the other hand, using equation \eqref{mixture_of_distributions_mean_variance} on the local means $M_i$, 
\begin{equation}
\begin{aligned}
\hat M(t) 
= 
\sum_{i=1}^I 
M_i(t) \mathcal P_{\mathcal D_i}(t)
&= 
\sum_{i=1}^I 
\left(
e^{At} \left[ M_{i,0} + A^{-1}b \right] - A^{-1}b  
\right) \mathcal P_{\mathcal D_i}(t) \\
&= 
\left(
e^{At} 
\sum_{i=1}^I M_{i,0} 
\mathcal P_{\mathcal D_i}(t) 
+ A^{-1}b \right) 
- A^{-1}b.
\end{aligned}
\end{equation}
Now $\hat M(t) = M(t)$ if $\sum_{i=1}^I M_{i,0} \mathcal P_{\mathcal D_i}(t) = u_0$, which is true if the weights stay constant and if the initial means $M_i(0)$ are consistent with the global mean $M(0)$.

\noindent
\textbf{For the variance}, we have, on the one hand, using a unimodal exact moment model
\begin{equation}
\Sigma(t) = e^{At} \left[ \Sigma_0 + A^{-1}b \right] - A^{-1}b.
\end{equation}
On the other hand, 
\begin{equation}
\begin{aligned}
\hat \Sigma (t)
= 
\sum_{i=1}^I 
(M_i(t)^2 + \Sigma_i^2) \mathcal P_{\mathcal D_i}(t) - M^2 
&= 
\sum_{i=1}^I \Sigma_i^2 \mathcal P_{\mathcal D_i}(t)
+
\sum_{i=1}^I  M_i(t)^2 \mathcal P_{\mathcal D_i} + 
 - M^2 
\\
&= 
\sum_{i=1}^I \Sigma_i^2 \mathcal P_{\mathcal D_i}(t) \\
&= 
\left(
e^{At} 
\sum_{i=1}^I \Sigma_{i,0} 
\mathcal P_{\mathcal D_i}(t) 
+ A^{-1}b \right) 
- A^{-1}b  
\end{aligned}.
\end{equation}
Now $\hat \Sigma(t) = \Sigma(t)$ if $\sum_{i=1}^I \Sigma_{i,0} \mathcal P_{\mathcal D_i}(t) = u_0$, which is true if the weights stay constant and if the initial variances $\Sigma_i(0)$ are consistent with the global mean $\Sigma(0)$.

\section{Proof of \cref{consistency_of_operators_S_and_T} Cconsistency property of the operators $\mathcal T$ and $\mathcal S$)}
\begin{proof}
For (i), we first write 
$\mathcal T \left( 
\mathcal S(u),  
u \right) 
= 
\mathcal T \left( 
(\mathbb E[u], \mathbb V[u]),  
u \right)
$. If $\mathbb E[u] = 0$, then $\sigma = 0$ thus (i) holds. If $\mathbb E[u] \neq 0$ then (i) also holds because, in that case,
\begin{equation}
\mathcal T \left( 
(\mathbb E[u], \mathbb V[u]),  
u \right)
= 
\sqrt{ \frac{ \mathbb V[u]}{\mathbb V[u]} }
\left(
u - \mathbb E [u]
\right) + \mathbb E [u]
= u
\end{equation}

For the proof of (ii), we have 
\begin{equation}
\begin{aligned}
\mathbb E \left[
\sqrt{ \frac{ \sigma}{\mathbb V[u]} }
\left(
u - \mathbb E [u]
\right) + \mu 
\right ] 
&= \mu,
 \\
\mathbb V \left[
\sqrt{ \frac{ \sigma}{\mathbb V[u]} }
\left(
u - \mathbb E [u]
\right) + \mu 
\right ]
&= 
\frac{ \sigma}{\mathbb V[u]} \mathbb V [u - \mathbb E [u]]
= \sigma.
\end{aligned}
\end{equation}
\end{proof}

\section{Model for particle fractions in learning-based MC-moments Parareal}
\label{appendix_model_particle_fractions}
\paragraph{Bimodal SDE}
In \cite[Chapter XI, equation 1.4]{kampen_stochastic_1981} this model is presented for the evolution of particle fractions $\mathcal  P_{\mathcal D_1}$ and $\mathcal P_{\mathcal D_2}$ in a bimodal SDE:
\begin{equation}
\frac{d \mathcal P_{\mathcal D_1}}{dt}
 = 
- \frac{d \mathcal P_{\mathcal D_2}}{dt}
= 
\frac{\mathcal P_{\mathcal D_2}}{\tau_{21}} 
- 
\frac{\mathcal P_{\mathcal D_1}}{\tau_{12}},
\label{equation_van_kampen}
\end{equation}
with $\mathcal P_{\mathcal D_1}(0) = \mathcal P_{\mathcal D_1, 0}$ 
and 
$\mathcal P_{\mathcal D_2}(0) = 1 - \mathcal P_{\mathcal D_1, 0}$.
In this model $\tau_{12}$ and $\tau_{21}$ are the exit times from region $1$ to region $2$, and reversely. 
These quantities are difficult to obtain in general without (sequential) time simulation of sample paths.
For a bimodal SDE, it holds that $\mathcal P_{\mathcal D_1}(t) = 1 - \mathcal P_{\mathcal D_2}(t)$.
Thus \eqref{equation_van_kampen} can be rewritten as follows
\begin{equation}
\begin{aligned}
\frac{d \mathcal P_{1}}{dt}
&=
\frac{1 - \mathcal P_{\mathcal D_1}}{\tau_{21}} 
-
\frac{\mathcal P_{\mathcal D_1}}{\tau_{12}}  
&=
-\mathcal P_{\mathcal D_1} 
\left(
\frac{1}{\tau_{21}} + 
\frac{1}{\tau_{12}} 
\right)
+
\frac{1}{\tau_{12}},
\end{aligned}
\end{equation}
with $\mathcal P_{\mathcal D_1}(0) = \mathcal P_{\mathcal D_1, 0}$.
This ODE has the exact solution
\begin{equation}
\mathcal P_{\mathcal D_1}(t) 
= 
\exp 
\left(
\underbrace{
-
\frac{\tau_{21} + \tau_{12}}{\tau_{21} \tau_{12}}
}_{\alpha} t
\right) 
\left[ 
\mathcal P_{\mathcal D_1,0} 
-
\underbrace{
\frac{\tau_{21}}{\tau_{21} + \tau_{12}}}_{\beta_1}
\right]
+
\underbrace{
\frac{\tau_{21}}{\tau_{21} + \tau_{12}}}_{\beta_1}.
\end{equation}
It is possible to derive a similar model for $\mathcal P_{\mathcal D_2}(t)$.
\paragraph{Multimodal SDEs}
For multimodal SDEs, the model \eqref{equation_van_kampen} can be generalised as follows:
\begin{equation}
\frac{d \mathcal P_{\mathcal D_i}}{dt}
= 
\sum_{j \neq i} \frac{\mathcal P_{\mathcal D_j}}{\tau_{ji}}
-
\sum_{j \neq i} \frac{\mathcal P_{\mathcal D_i}}{\tau_{ij}}.
\label{generalised_Van_Kampen}
\end{equation}
Because no particles disappear or emanate during the evolution McKean-Vlasov dynamics \eqref{general_equation}, 
it must hold that $\sum_{i=1}^I \mathcal P_{\mathcal D_i}(t) = 1$.
Let $\omega_{ij} = 1/\tau_{ij}$, then it is possible to reorganise equation \eqref{generalised_Van_Kampen}
in matrix form:
\begin{equation}
\frac{d}{dt} 
\left[ 
\begin{matrix}
\mathcal P_{\mathcal D_1} \\
\mathcal P_{\mathcal D_2} \\
\vdots  \\
\mathcal P_{\mathcal D_I}
\end{matrix}
\right]
=
\underbrace{
\left[ 
\begin{matrix}
-\sum_{j \neq 1} \omega_{1j} 
	& \omega_{21} 
	& 
	& \hdots
	& \omega_{I1} \\ 
\omega_{12} 
	& -\sum_{j \neq 2} \omega_{2j} 
	& 
	& \hdots
	& \omega_{I2} \\ 
\hdots & \hdots & & & \vdots \\	
\omega_{1I} 
	& \omega_{2I}
	& 
	& \hdots
	& - \sum_{j \neq I} \omega_{Ij}  \\ 
\end{matrix}
\right]}_{A}
\left[ 
\begin{matrix}
\mathcal P_{\mathcal D_1} \\
\mathcal P_{\mathcal D_2} \\
\vdots  \\
\mathcal P_{\mathcal D_I}
\end{matrix}
\right]
\label{generalised_model_multimodal_PFs}
\end{equation}
with the initial condition $
\left[ 
\begin{matrix}
\mathcal P_{\mathcal D_1}(0) &
\mathcal P_{\mathcal D_2}(0) &
\hdots  &
\mathcal P_{\mathcal D_I}(0)
\end{matrix}
\right]
= 
\left[ 
\begin{matrix}
\mathcal P_{\mathcal D_1, 0} &
\mathcal P_{\mathcal D_2, 0} &
\hdots  &
\mathcal P_{\mathcal D_I, 0}
\end{matrix}
\right]
$
with $\sum_i \mathcal P_{\mathcal D_i,0} = 1$.
\begin{property}[Conservation (invariance) property of the model \eqref{generalised_model_multimodal_PFs}]
It holds that $\sum_i \mathcal P_{\mathcal D_i}(t) = 1$ for all $t \geq 0$.
\end{property}
\begin{proof}
First remark that the sum of the rows of the matrix $A$ in \eqref{generalised_model_multimodal_PFs} equals zero.
The ODE \eqref{generalised_model_multimodal_PFs} thus has in invariant. 
Indeed, an ODE for 
$\mathcal P_{\mathrm{tot}}(t) 
= 
\sum_i \mathcal P_{\mathcal D_i}(t)$ can be constructed by summing the left-hand side and the right-hand side of \eqref{generalised_model_multimodal_PFs}: $\frac{d}{dt}\mathcal P_{\mathrm{tot}} = 0$.
Thus we have that $\mathcal P_{\mathrm{tot}}$ equals its value at the initial condition, namely $\mathcal P_{\mathrm{tot}} = \sum_i \mathcal P_{\mathcal D_i, 0} = 1$.
\end{proof}

\section{Extra figure} 
\label{extra_figures_pf}
\,
\begin{figure}[H]
\begin{subfigure}[t]{.33\linewidth}
\includegraphics[width=1\textwidth]{{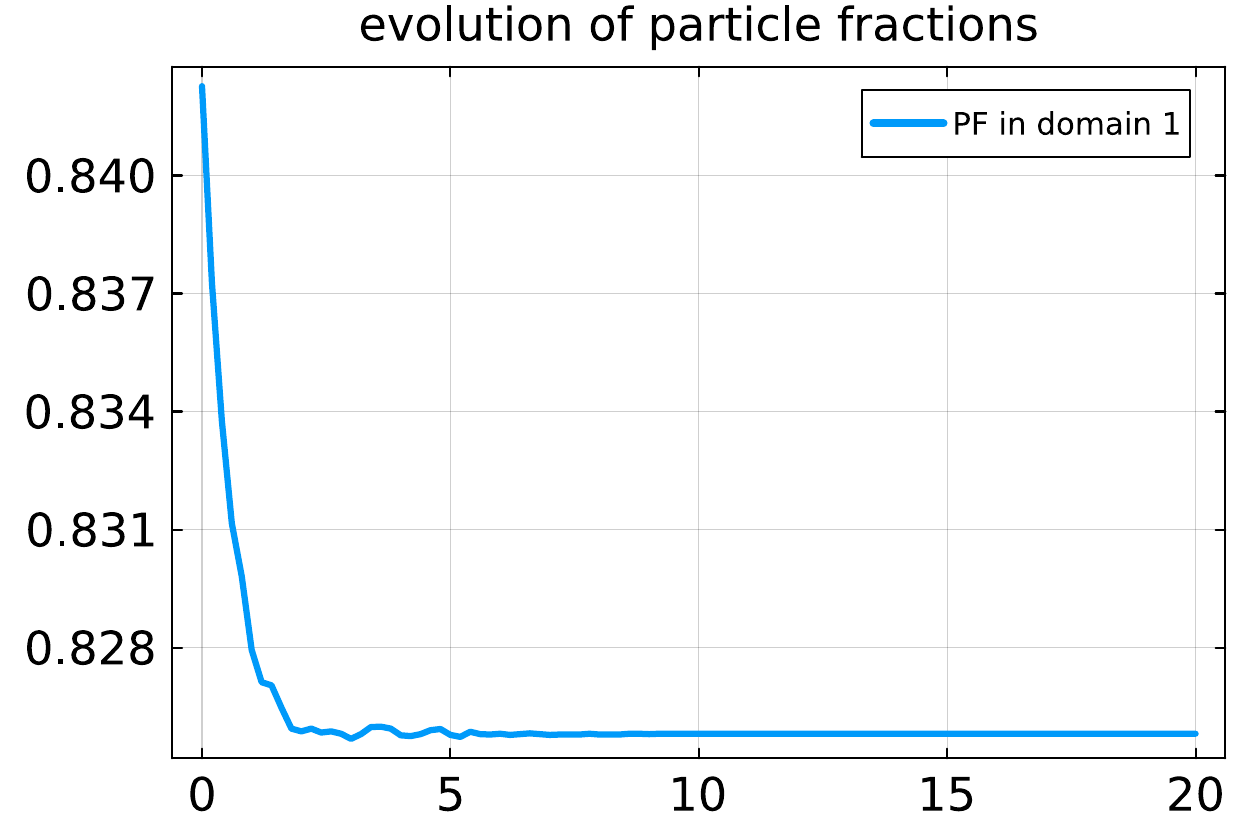}}%
	\caption{Setup 1}
  \end{subfigure}
  \begin{subfigure}[t]{.33\linewidth}
\includegraphics[width=1\textwidth]{{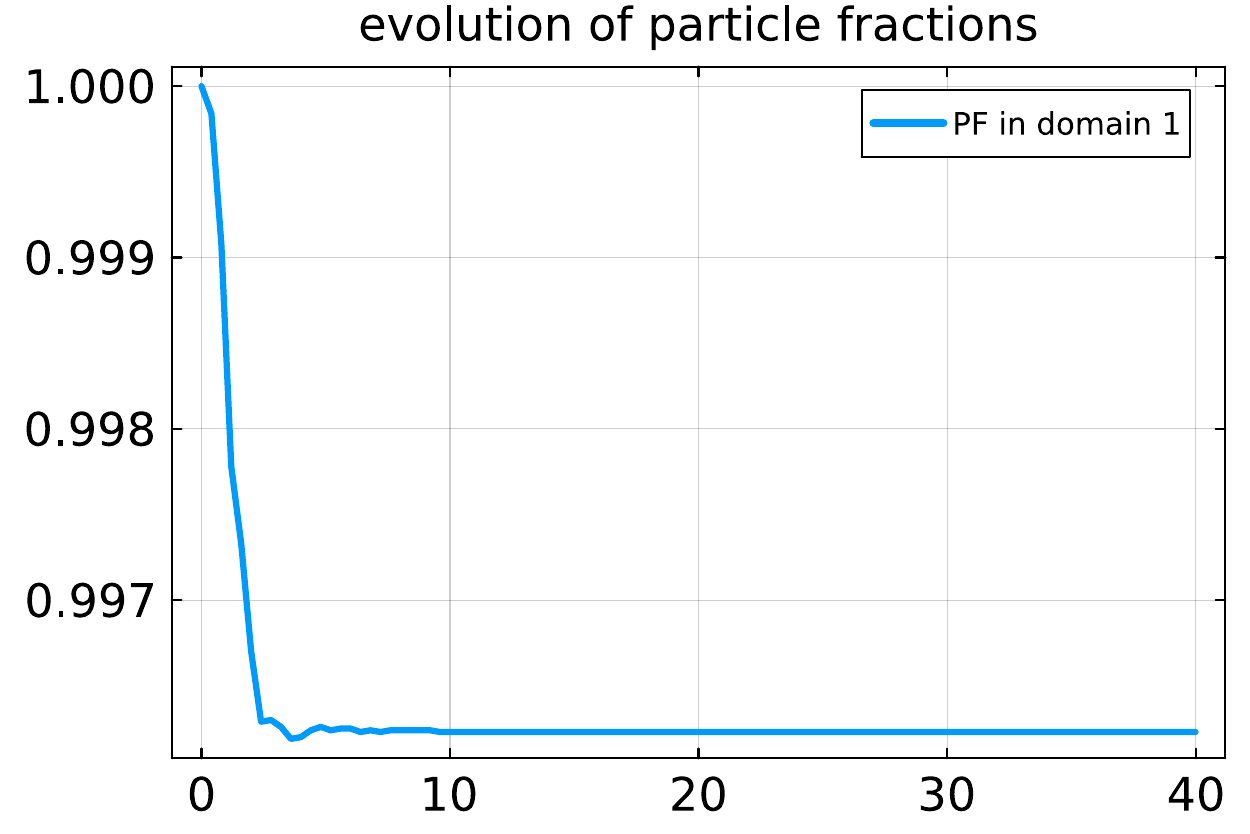}}%
	\caption{Setup 2}
  \end{subfigure}
  \begin{subfigure}[t]{.33\linewidth}
\includegraphics[width=1\textwidth]{{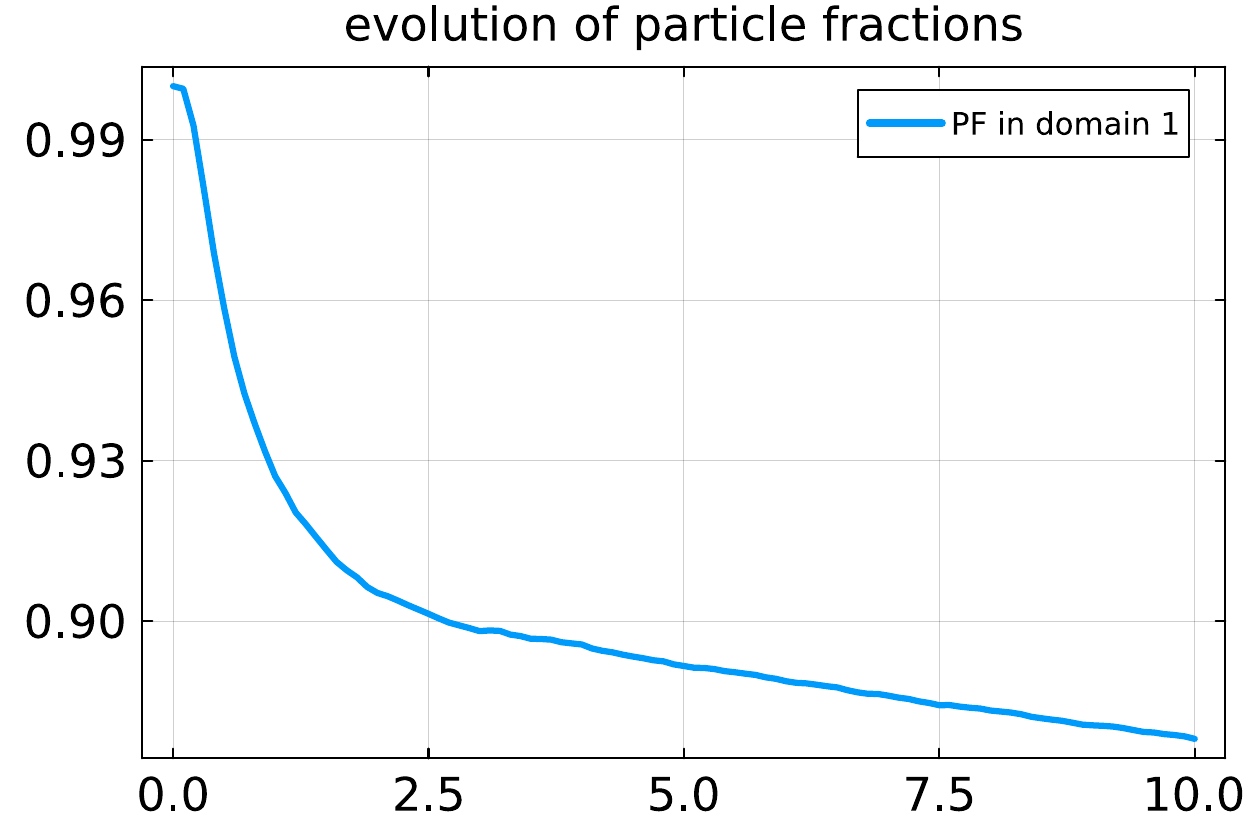}}%
	\caption{Setup 3}
  \end{subfigure}
    \begin{minipage}{.25\linewidth}
    \vspace*{-2cm}
  \caption{Overview of test systems for the bimodal SDE: 
evolution of particle fractions in function of time.}
\label{fig_overview_evolution_PFs}
    \end{minipage}
    \hspace*{0.06\linewidth}
    \begin{subfigure}[t]{.33\linewidth}
\includegraphics[width=1\textwidth]{{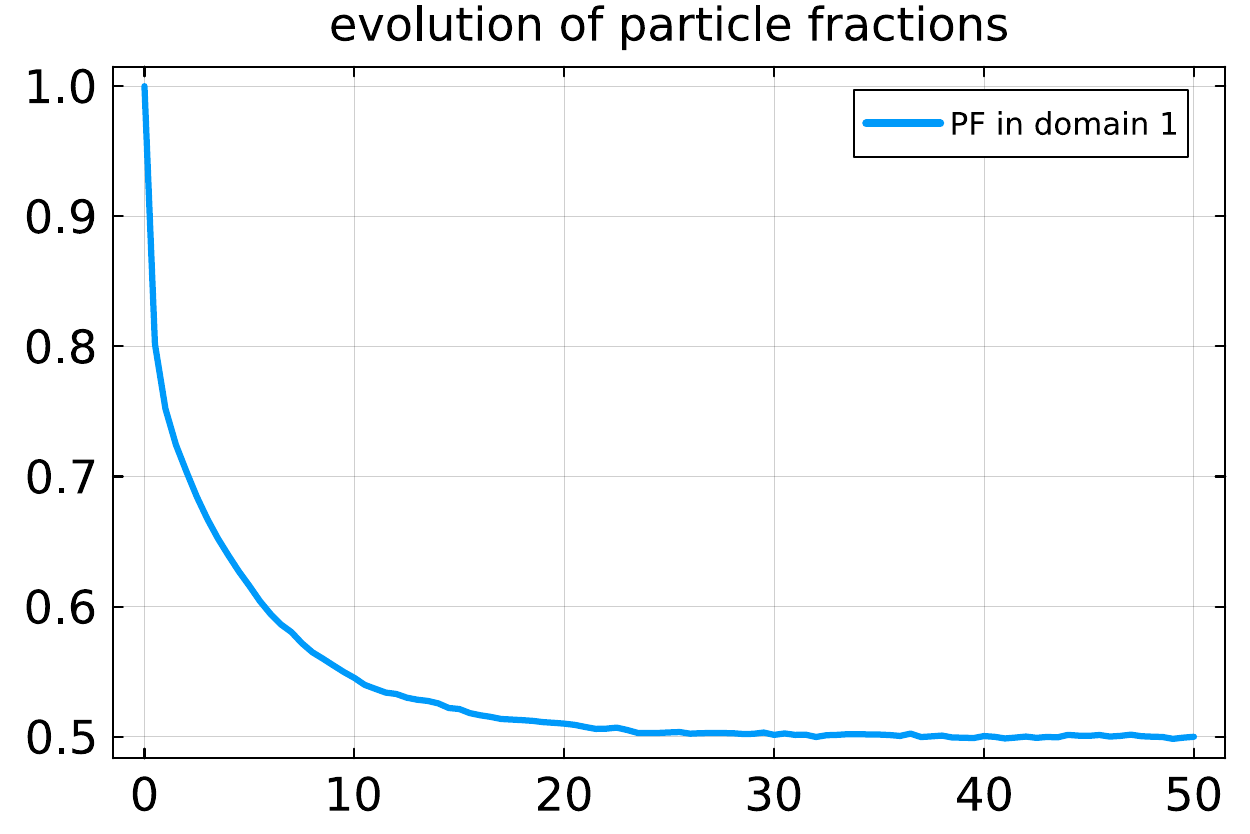}}%
	\caption{Setup 4}
  \end{subfigure}
  \begin{subfigure}[t]{.33\linewidth}
\includegraphics[width=1\textwidth]{{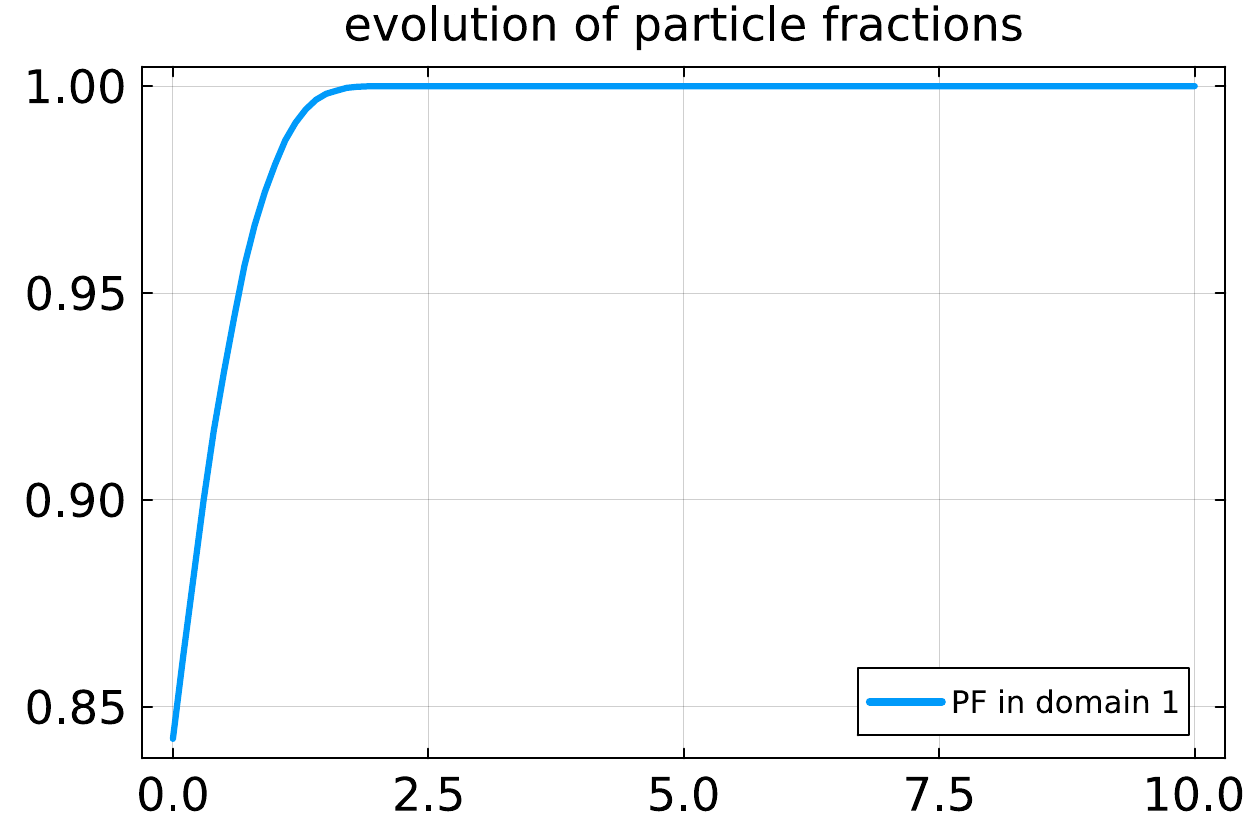}}%
	\caption{Setup 5}
  \end{subfigure}
\end{figure}

\bibliographystyle{siamplain}
\bibliography{all}

\begin{thebibliography}{10}

\bibitem{alspach_nonlinear_1972}
{\sc D.~Alspach and H.~Sorenson}, {\em Nonlinear {Bayesian} estimation using
  {Gaussian} sum approximations}, IEEE Transactions on Automatic Control, 17
  (1972), pp.~439--448, \url{https://doi.org/10.1109/TAC.1972.1100034},
  \url{http://ieeexplore.ieee.org/document/1100034/} (accessed 2024-10-23).

\bibitem{arasaratnam_cubature_2009}
{\sc I.~Arasaratnam and S.~Haykin}, {\em Cubature {Kalman} {Filters}}, IEEE
  Transactions on Automatic Control, 54 (2009), pp.~1254--1269,
  \url{https://doi.org/10.1109/TAC.2009.2019800},
  \url{http://ieeexplore.ieee.org/document/4982682/} (accessed 2024-10-23).

\bibitem{archambeau07a}
{\sc C.~Archambeau, D.~Cornford, M.~Opper, and J.~Shawe-Taylor}, {\em Gaussian
  process approximations of stochastic differential equations}, in Gaussian
  Processes in Practice, N.~D. Lawrence, A.~Schwaighofer, and
  J.~Quiñonero~Candela, eds., vol.~1 of Proceedings of Machine Learning
  Research, Bletchley Park, UK, 12--13 Jun 2007, PMLR, pp.~1--16,
  \url{https://proceedings.mlr.press/v1/archambeau07a.html}.

\bibitem{arnold_stochastic_1974}
{\sc L.~Arnold}, {\em Stochastic differential equations: theory and
  applications}, Wiley, New York, 1974.

\bibitem{bal_parallelization_2003}
{\sc G.~Bal}, {\em Parallelization in time of (stochastic) ordinary
  differential equations},  (2003), pp.~1--23,
  \url{https://galton.uchicago.edu/~guillaumebal/PAPERS/paralleltime.pdf}.

\bibitem{baladron_mean_field_2012}
{\sc J.~Baladron, D.~Fasoli, O.~Faugeras, and J.~Touboul}, {\em Mean-field
  description and propagation of chaos in networks of {Hodgkin}-{Huxley} and
  {FitzHugh}-{Nagumo} neurons}, The Journal of Mathematical Neuroscience, 2
  (2012), p.~10, \url{https://doi.org/10.1186/2190-8567-2-10},
  \url{https://mathematical-neuroscience.springeropen.com/articles/10.1186/2190-8567-2-10}.

\bibitem{jourdain_2019}
{\sc O.~Bencheikh and B.~Jourdain}, {\em Bias behaviour and antithetic sampling
  in mean-field particle approximations of {SDEs} nonlinear in the sense of
  {McKean}}, ESAIM: Proceedings and Surveys, 65 (2019), pp.~219--235,
  \url{https://doi.org/10.1051/proc/201965219},
  \url{https://www.esaim-proc.org/10.1051/proc/201965219}.

\bibitem{blouza_parallel_2010}
{\sc A.~Blouza, L.~Boudin, and S.~M. Kaber}, {\em Parallel in time algorithms
  with reduction methods for solving chemical kinetics}, Communications in
  Applied Mathematics and Computational Science, 5 (2010), pp.~241--263,
  \url{https://doi.org/10.2140/camcos.2010.5.241},
  \url{http://msp.org/camcos/2010/5-2/p04.xhtml}.

\bibitem{ignace_software_paper}
{\sc I.~Bossuyt}, {\em {MC-moments-Parareal}, implementation in {Julia}}, 2024,
  \url{https://gitlab.kuleuven.be/numa/public/mc-moments-parareal}.

\bibitem{bossuyt_new}
{\sc I.~Bossuyt, S.~Vandewalle, and G.~Samaey}, {\em {Micro-macro Parareal,
  from ODEs to SDEs and back again}}, Oct. 2023,
  \url{http://arxiv.org/abs/2310.11365} (accessed 2023-11-08).
\newblock arXiv:2310.11365 [math.NA, physics, stat].

\bibitem{bossy_stochastic_2005}
{\sc M.~Bossy}, {\em Some stochastic particle methods for nonlinear parabolic
  {PDEs}}, ESAIM: Proceedings, 15 (2005), pp.~18--57,
  \url{https://doi.org/10.1051/proc:2005019},
  \url{http://www.esaim-proc.org/10.1051/proc:2005019}.

\bibitem{bossy_synchronization_2019}
{\sc M.~Bossy, J.~Fontbona, and H.~Olivero}, {\em Synchronization of stochastic
  mean field networks of {Hodgkin}–{Huxley} neurons with noisy channels},
  Journal of Mathematical Biology, 78 (2019), pp.~1771--1820,
  \url{https://doi.org/10.1007/s00285-019-01326-7},
  \url{http://link.springer.com/10.1007/s00285-019-01326-7}.

\bibitem{bossy_stochastic_1997}
{\sc M.~Bossy and D.~Talay}, {\em A stochastic particle method for the
  {McKean}-{Vlasov} and the {Burgers} equation}, Mathematics of Computation, 66
  (1997), pp.~157--193, \url{https://doi.org/10.1090/S0025-5718-97-00776-X},
  \url{http://www.ams.org/journal-getitem?pii=S0025-5718-97-00776-X}.

\bibitem{caflisch_monte_1998}
{\sc R.~E. Caflisch}, {\em Monte {Carlo} and quasi-{Monte} {Carlo} methods},
  Acta Numerica, 7 (1998), pp.~1--49,
  \url{https://doi.org/10.1017/S0962492900002804},
  \url{https://www.cambridge.org/core/product/identifier/S0962492900002804/type/journal_article}.

\bibitem{carrel_low-rank_2023}
{\sc B.~Carrel, M.~J. Gander, and B.~Vandereycken}, {\em Low-rank {Parareal}: a
  low-rank parallel-in-time integrator}, BIT Numerical Mathematics, 63 (2023),
  p.~13, \url{https://doi.org/10.1007/s10543-023-00953-3},
  \url{https://link.springer.com/10.1007/s10543-023-00953-3} (accessed
  2024-10-25).

\bibitem{Chen2014}
{\sc F.~Chen, J.~S. Hesthaven, and X.~Zhu}, {\em On the Use of Reduced Basis
  Methods to Accelerate and Stabilize the Parareal Method}, Springer
  International Publishing, Cham, 2014, pp.~187--214,
  \url{https://doi.org/10.1007/978-3-319-02090-7_7},
  \url{https://doi.org/10.1007/978-3-319-02090-7_7}.

\bibitem{dabaghi_hybrid_2023}
{\sc J.~Dabaghi, Y.~Maday, and A.~Zoia}, {\em A hybrid parareal {Monte} {Carlo}
  algorithm for parabolic problems}, Journal of Computational and Applied
  Mathematics, 420 (2023), p.~114800,
  \url{https://doi.org/10.1016/j.cam.2022.114800},
  \url{https://linkinghub.elsevier.com/retrieve/pii/S0377042722004071}.

\bibitem{del_moral_stability_2018}
{\sc P.~Del~Moral and J.~Tugaut}, {\em On the stability and the uniform
  propagation of chaos properties of {Ensemble} {Kalman}–{Bucy} filters}, The
  Annals of Applied Probability, 28 (2018),
  \url{https://doi.org/10.1214/17-AAP1317},
  \url{https://projecteuclid.org/journals/annals-of-applied-probability/volume-28/issue-2/On-the-stability-and-the-uniform-propagation-of-chaos-properties/10.1214/17-AAP1317.full}.

\bibitem{emmett_toward_2012}
{\sc M.~Emmett and M.~L. Minion}, {\em Toward an efficient parallel in time
  method for partial differential equations}, Communications in Applied
  Mathematics and Computational Science, 7 (2012),
  \url{https://doi.org/10.2140/camcos.2012.7.105},
  \url{https://projecteuclid.org/euclid.camcos/1513732042}.
\newblock ISBN: 0001404105.

\bibitem{engblom_parallel_2009}
{\sc S.~Engblom}, {\em Parallel in {Time} {Simulation} of {Multiscale}
  {Stochastic} {Chemical} {Kinetics}}, Multiscale Modeling \& Simulation, 8
  (2009), pp.~46--68, \url{https://doi.org/10.1137/080733723},
  \url{http://epubs.siam.org/doi/10.1137/080733723}.

\bibitem{falgout_parallel_2014}
{\sc R.~D. Falgout, S.~Friedhoff, T.~V. Kolev, S.~P. Maclachan, and J.~B.
  Schroder}, {\em Parallel {Time} {Integration} with {Multigrid}}, SIAM Journal
  of Scientific Computing., 36 (2014), pp.~C635--C661,
  \url{https://doi.org/10.1137/130944230}.

\bibitem{fruhwirth-schnatter_finite_2006}
{\sc S.~Frühwirth-Schnatter}, {\em Finite mixture and {Markov} switching
  models}, Springer series in statistics, Springer, New York, 2006.

\bibitem{gander_analysis_2008}
{\sc M.~Gander and M.~Petcu}, {\em Analysis of a {Krylov} subspace enhanced
  parareal algorithm for linear problems}, ESAIM: Proceedings, 25 (2008),
  pp.~114--129, \url{https://doi.org/10.1051/proc:082508}.

\bibitem{Gander2015}
{\sc M.~J. Gander}, {\em 50 {Years} of {Time} {Parallel} {Time} {Integration}},
  in Multiple {Shooting} and {Time} {Domain} {Decomposition} {Methods},
  T.~Carraro, M.~Geiger, S.~Körkel, and R.~Rannacher, eds., vol.~9, Springer
  International Publishing, Cham, 2015, pp.~69--113,
  \url{https://doi.org/10.1007/978-3-319-23321-5_3},
  \url{https://link.springer.com/10.1007/978-3-319-23321-5_3}.
\newblock Series Title: Contributions in Mathematical and Computational
  Sciences.

\bibitem{gander_parareal_2024}
{\sc M.~J. Gander, M.~Ohlberger, and S.~Rave}, {\em A {Parareal} algorithm
  without {Coarse} {Propagator}?}, 2024,
  \url{https://doi.org/10.48550/ARXIV.2409.02673},
  \url{https://arxiv.org/abs/2409.02673} (accessed 2024-10-25).
\newblock Version Number: 1.

\bibitem{gander_analysis_2007}
{\sc M.~J. Gander and S.~Vandewalle}, {\em Analysis of the parareal
  time-parallel time-integration method}, SIAM Journal on Scientific Computing,
  29 (2007), pp.~556--578, \url{https://doi.org/10.1137/05064607X},
  \url{https://epubs.siam.org/doi/abs/10.1137/05064607X}.

\bibitem{haji_ali_multilevel_2018}
{\sc A.-L. Haji-Ali and R.~Tempone}, {\em Multilevel and {Multi}-index {Monte}
  {Carlo} methods for the {McKean}–{Vlasov} equation}, Statistics and
  Computing, 28 (2018), pp.~923--935,
  \url{https://doi.org/10.1007/s11222-017-9771-5},
  \url{http://link.springer.com/10.1007/s11222-017-9771-5}.

\bibitem{higham_kloeden_introduction_2021}
{\sc D.~J. Higham and P.~E. Kloeden}, {\em An introduction to the numerical
  simulation of stochastic differential equations}, no.~169 in Other titles in
  applied mathematics, Society for Industrial and Applied Mathematics,
  Philadelphia, 2021.

\bibitem{jazwinski_stochastic_1970}
{\sc A.~H. Jazwinski}, {\em Stochastic processes and filtering theory}, no.~v.
  64 in Mathematics in science and engineering, Academic Press, New York, 1970.

\bibitem{julier_new_2000}
{\sc S.~Julier, J.~Uhlmann, and H.~Durrant-Whyte}, {\em A new method for the
  nonlinear transformation of means and covariances in filters and estimators},
  IEEE Transactions on Automatic Control, 45 (2000), pp.~477--482,
  \url{https://doi.org/10.1109/9.847726},
  \url{http://ieeexplore.ieee.org/document/847726/} (accessed 2024-10-23).

\bibitem{kalai_disentangling_2012}
{\sc A.~T. Kalai, A.~Moitra, and G.~Valiant}, {\em Disentangling {Gaussians}},
  Communications of the ACM, 55 (2012), pp.~113--120,
  \url{https://doi.org/10.1145/2076450.2076474},
  \url{https://dl.acm.org/doi/10.1145/2076450.2076474} (accessed 2024-10-28).

\bibitem{kampen_stochastic_1981}
{\sc N.~G.~v. Kampen}, {\em Stochastic processes in physics and chemistry},
  North-{Holland} personal library, Elsevier, Amsterdam ; Boston, 3rd ed~ed.,
  2007.
\newblock OCLC: ocm81453662.

\bibitem{kloeden_gauss_quadrature_2017}
{\sc P.~Kloeden and T.~Shardlow}, {\em Gauss-quadrature method for
  one-dimensional mean-field sdes}, SIAM Journal on Scientific Computing, 39
  (2017), pp.~A2784--A2807, \url{https://doi.org/10.1137/16M1095688},
  \url{https://epubs.siam.org/doi/10.1137/16M1095688}.

\bibitem{kloeden_platen_1999}
{\sc P.~E. Kloeden and E.~Platen}, {\em Numerical solution of stochastic
  differential equations}, no.~23 in Applications of mathematics, Springer,
  Berlin Heidelberg, 1999.

\bibitem{kostur_nonequilibrium_2002}
{\sc M.~Kostur, J.~Łuczka, and L.~Schimansky-Geier}, {\em Nonequilibrium
  coupled {Brownian} phase oscillators}, Physical Review E, 65 (2002),
  p.~051115, \url{https://doi.org/10.1103/PhysRevE.65.051115},
  \url{https://link.aps.org/doi/10.1103/PhysRevE.65.051115}.

\bibitem{lambert_chewi_bach_bonnabel_rigollet}
{\sc M.~Lambert, S.~Chewi, F.~Bach, S.~Bonnabel, and P.~Rigollet}, {\em
  Variational inference via wasserstein gradient flows}, in Advances in Neural
  Information Processing Systems, S.~Koyejo, S.~Mohamed, A.~Agarwal,
  D.~Belgrave, K.~Cho, and A.~Oh, eds., vol.~35, Curran Associates, Inc., 2022,
  pp.~14434--14447,
  \url{https://proceedings.neurips.cc/paper_files/paper/2022/file/5d087955ee13fe9a7402eedec879b9c3-Paper-Conference.pdf}.

\bibitem{Legoll2020}
{\sc F.~Legoll, T.~Leli{\`{e}}vre, K.~Myerscough, and G.~Samaey}, {\em
  {Parareal computation of stochastic differential equations with time-scale
  separation: a numerical convergence study}}, Computing and Visualization in
  Science, 23 (2020), \url{https://doi.org/10.1007/s00791-020-00329-y},
  \url{https://doi.org/10.1007/s00791-020-00329-y}.

\bibitem{legoll_adaptive_Parareal}
{\sc F.~Legoll, T.~Leli\`{e}vre, and U.~Sharma}, {\em An adaptive parareal
  algorithm: Application to the simulation of molecular dynamics trajectories},
  SIAM Journal on Scientific Computing, 44 (2022), pp.~B146--B176,
  \url{https://doi.org/10.1137/21M1412979},
  \url{https://doi.org/10.1137/21M1412979},
  \url{https://arxiv.org/abs/https://doi.org/10.1137/21M1412979}.

\bibitem{Legoll2013}
{\sc F.~Legoll, T.~Lelièvre, and G.~Samaey}, {\em {A micro-macro parareal
  algorithm: application to singularly perturbed differential equations}}, SIAM
  Journal on Scientific Computing, 2013-01, 35 (2013), pp.~p.A1951--A1986,
  \url{https://doi.org/10.1137/120872681}.

\bibitem{li_numerical_2021}
{\sc L.~Li, J.~Lu, J.~C. Mattingly, and L.~Wang}, {\em Numerical methods for
  stochastic differential equations based on {Gaussian} mixture},
  Communications in Mathematical Sciences, 19 (2021), pp.~1549--1577,
  \url{https://doi.org/10.4310/CMS.2021.v19.n6.a5},
  \url{https://link.intlpress.com/JDetail/1806261779988697090} (accessed
  2024-10-25).

\bibitem{lions_resolution_2001}
{\sc J.-L. Lions, Y.~Maday, and G.~Turinici}, {\em Résolution d'edp par un
  schéma en temps "pararéel"}, Comptes Rendus de l'Académie des Sciences -
  Series I - Mathematics, 332 (2001), pp.~661--668,
  \url{https://doi.org/https://doi.org/10.1016/S0764-4442(00)01793-6}.

\bibitem{najafi_fast_2016}
{\sc E.~Najafi, R.~Babuška, and G.~A.~D. Lopes}, {\em A fast sampling method
  for estimating the domain of attraction}, Nonlinear Dynamics, 86 (2016),
  pp.~823--834, \url{https://doi.org/10.1007/s11071-016-2926-7},
  \url{http://link.springer.com/10.1007/s11071-016-2926-7}.

\bibitem{thalhammer_MLMC_space_time_multigrid}
{\sc M.~Neumüller and A.~Thalhammer}, {\em Combining {Space}-{Time}
  {Multigrid} {Techniques} with {Multilevel} {Monte} {Carlo} {Methods} for
  {SDEs}}, in Domain {Decomposition} {Methods} in {Science} and {Engineering}
  {XXIV}, P.~E. Bjørstad, S.~C. Brenner, L.~Halpern, H.~H. Kim, R.~Kornhuber,
  T.~Rahman, and O.~B. Widlund, eds., vol.~125, Springer International
  Publishing, Cham, 2018, pp.~493--501,
  \url{https://doi.org/10.1007/978-3-319-93873-8_47},
  \url{http://link.springer.com/10.1007/978-3-319-93873-8_47}.
\newblock Series Title: Lecture Notes in Computational Science and Engineering.

\bibitem{ong_applications_2020}
{\sc B.~W. Ong and J.~B. Schroder}, {\em Applications of time parallelization},
  Computing and Visualization in Science, 23 (2020), p.~11,
  \url{https://doi.org/10.1007/s00791-020-00331-4},
  \url{https://link.springer.com/10.1007/s00791-020-00331-4}.

\bibitem{pamela_neural-parareal_2025}
{\sc S.~J.~P. Pamela, N.~Carey, J.~Brandstetter, R.~Akers, L.~Zanisi,
  J.~Buchanan, V.~Gopakumar, M.~Hoelzl, G.~Huijsmans, K.~Pentland, T.~James,
  G.~Antonucci, and t.~J. Team}, {\em Neural-{Parareal}: {Dynamically}
  {Training} {Neural} {Operators} as {Coarse} {Solvers} for
  {Time}-{Parallelisation} of {Fusion} {MHD} {Simulations}}, Computer Physics
  Communications, 307 (2025), p.~109391,
  \url{https://doi.org/10.1016/j.cpc.2024.109391},
  \url{http://arxiv.org/abs/2405.01355} (accessed 2024-10-27).
\newblock arXiv:2405.01355 [physics].

\bibitem{pentland_stochastic_2022}
{\sc K.~Pentland, M.~Tamborrino, D.~Samaddar, and L.~C. Appel}, {\em Stochastic
  {Parareal}: {An} {Application} of {Probabilistic} {Methods} to
  {Time}-{Parallelization}}, SIAM Journal on Scientific Computing,  (2022),
  pp.~S82--S102, \url{https://doi.org/10.1137/21M1414231},
  \url{https://epubs.siam.org/doi/10.1137/21M1414231}.

\bibitem{rackauckas_differentialequationsjl_2017}
{\sc C.~Rackauckas and Q.~Nie}, {\em {DifferentialEquations}.jl – {A}
  {Performant} and {Feature}-{Rich} {Ecosystem} for {Solving} {Differential}
  {Equations} in {Julia}}, Journal of Open Research Software, 5 (2017), p.~15,
  \url{https://doi.org/10.5334/jors.151},
  \url{https://openresearchsoftware.metajnl.com/article/10.5334/jors.151/}.

\bibitem{reisinger_adaptive_2022}
{\sc C.~Reisinger and W.~Stockinger}, {\em An adaptive {Euler}–{Maruyama}
  scheme for {McKean}–{Vlasov} {SDEs} with super-linear growth and
  application to the mean-field {FitzHugh}–{Nagumo} model}, Journal of
  Computational and Applied Mathematics, 400 (2022), p.~113725,
  \url{https://doi.org/10.1016/j.cam.2021.113725},
  \url{https://linkinghub.elsevier.com/retrieve/pii/S0377042721003472}.

\bibitem{rodriguez_statistical_1996}
{\sc R.~Rodriguez and H.~C. Tuckwell}, {\em Statistical properties of
  stochastic nonlinear dynamical models of single spiking neurons and neural
  networks}, Physical Review E, 54 (1996), pp.~5585--5590,
  \url{https://doi.org/10.1103/PhysRevE.54.5585}.

\bibitem{schnoerr_validity_2014}
{\sc D.~Schnoerr, G.~Sanguinetti, and R.~Grima}, {\em Validity conditions for
  moment closure approximations in stochastic chemical kinetics}, The Journal
  of Chemical Physics, 141 (2014), p.~084103,
  \url{https://doi.org/10.1063/1.4892838},
  \url{https://pubs.aip.org/aip/jcp/article/352569}.

\bibitem{Solin_NEURIPS}
{\sc A.~Solin, E.~Tamir, and P.~Verma}, {\em Scalable inference in sdes by
  direct matching of the fokker\textendash planck\textendash kolmogorov
  equation}, in Advances in Neural Information Processing Systems, M.~Ranzato,
  A.~Beygelzimer, Y.~Dauphin, P.~Liang, and J.~W. Vaughan, eds., vol.~34,
  Curran Associates, Inc., 2021, pp.~417--429,
  \url{https://proceedings.neurips.cc/paper_files/paper/2021/file/03e4d3f831100d4355663f3d425d716b-Paper.pdf}.

\bibitem{son_doan_mean_square_2015}
{\sc T.~Son~Doan, M.~Rasmussen, P.~E.~Kloeden, {,Department of Mathematics,
  Imperial College London, 180 Queen's Gate, London SW7 2AZ}, and {,School of
  Mathematics and Statistics, Huazhong University of Science \& Technology,
  Wuhan 430074}}, {\em The mean-square dichotomy spectrum and a bifurcation to
  a mean-square attractor}, Discrete \& Continuous Dynamical Systems - B, 20
  (2015), pp.~875--887, \url{https://doi.org/10.3934/dcdsb.2015.20.875},
  \url{http://aimsciences.org//article/doi/10.3934/dcdsb.2015.20.875}.

\bibitem{sukys_momentclosurejl_2021}
{\sc A.~Sukys and R.~Grima}, {\em {MomentClosure}.jl: automated moment closure
  approximations in {Julia}}, Bioinformatics, 38 (2021), pp.~289--290,
  \url{https://doi.org/10.1093/bioinformatics/btab469},
  \url{https://academic.oup.com/bioinformatics/article/38/1/289/6309452}.

\bibitem{snitzman_1991}
{\sc A.-S. Sznitman}, {\em Topics in propagation of chaos}, in Ecole d'{Eté}
  de {Probabilités} de {Saint}-{Flour} {XIX} — 1989, P.-L. Hennequin, ed.,
  vol.~1464, Springer Berlin Heidelberg, Berlin, Heidelberg, 1991,
  pp.~165--251, \url{https://doi.org/10.1007/BFb0085169},
  \url{http://link.springer.com/10.1007/BFb0085169}.
\newblock Series Title: Lecture Notes in Mathematics.

\bibitem{saarkkaa_solin_applied_SDEs_2019}
{\sc S.~Särkkä and A.~Solin}, {\em Applied stochastic differential
  equations}, no.~10 in Institute of {Mathematical} {Statistics} textbooks,
  Cambridge University Press, Cambridge ; New York, NY, 2019.

\bibitem{tang_potential_2017}
{\sc Y.~Tang, R.~Yuan, G.~Wang, X.~Zhu, and P.~Ao}, {\em Potential landscape of
  high dimensional nonlinear stochastic dynamics with large noise}, Scientific
  Reports, 7 (2017), p.~15762,
  \url{https://doi.org/10.1038/s41598-017-15889-2},
  \url{http://www.nature.com/articles/s41598-017-15889-2}.

\bibitem{terejanu_adaptive_2011}
{\sc G.~Terejanu, P.~Singla, T.~Singh, and P.~D. Scott}, {\em Adaptive
  {Gaussian} {Sum} {Filter} for {Nonlinear} {Bayesian} {Estimation}}, IEEE
  Transactions on Automatic Control, 56 (2011), pp.~2151--2156,
  \url{https://doi.org/10.1109/TAC.2011.2141550},
  \url{http://ieeexplore.ieee.org/document/5746510/} (accessed 2024-10-25).

\bibitem{thalhammer_phd_thesis}
{\sc A.~Thalhammer}, {\em Numerical methods for stochastic partial differential
  equations: Analysis of stability and efficiency}, PhD thesis, Johannes Kepler
  Universitat Linz, 2017.

\bibitem{zagli_dimension_2023}
{\sc N.~Zagli, G.~A. Pavliotis, V.~Lucarini, and A.~Alecio}, {\em Dimension
  reduction of noisy interacting systems}, Physical Review Research, 5 (2023),
  p.~013078, \url{https://doi.org/10.1103/PhysRevResearch.5.013078},
  \url{https://link.aps.org/doi/10.1103/PhysRevResearch.5.013078}.

\end{thebibliography}

\end{document}